\newif\ifroman
\romanfalse
\newif\ifpdf
\ifx\pdfoutput\undefined
\pdffalse % we are not running PDFLaTeX
\else
\pdfoutput=1 % we are running PDFLaTeX
\pdftrue \fi

\documentclass[reqno,twoside,11pt]{amsart}
\IfFileExists{cite.sty}{\usepackage[nosort]{cite}}{}
\usepackage{amsmath}
\usepackage{amsfonts}
\usepackage{amssymb}
\usepackage{wrapfig}
\usepackage{graphicx}
\usepackage{graphics}
\usepackage[normalem]{ulem}
\usepackage{verbatim}
\usepackage{xr}

\IfFileExists{epsf.def}{\input epsf.def}{\usepackage{epsf}}

\IfFileExists{myowntimes.sty}{\usepackage{myowntimes}}
{\usepackage{times}
\IfFileExists{mathrsfs.sty}
{\usepackage{mathrsfs}}{\newcommand{\mathscr}{\mathcal}}}

\DeclareFontFamily{OT1}{eusb}{} \DeclareFontShape{OT1}{eusb}{m}{n}
{<5> <6> <7> <8> <9> <10> <11> <12> <14.4> eusb10}{}
\DeclareMathAlphabet{\eusb}{OT1}{eusb}{m}{n}

\DeclareFontFamily{OT1}{eusm}{} \DeclareFontShape{OT1}{eusm}{m}{n}
{<5> <6> <7> <8> <9> <10> <11> <12> <14.4> eusm10}{}
\DeclareMathAlphabet{\eusm}{OT1}{eusm}{m}{n}

\DeclareFontFamily{OT1}{eufm}{} \DeclareFontShape{OT1}{eufm}{m}{n}
{<5> <6> <7> <8> <9> <10> <11> <12> <14.4> eufm10}{}
\DeclareMathAlphabet{\mathfrak}{OT1}{eufm}{m}{n}

\DeclareFontFamily{OT1}{fraktura}{}
\DeclareFontShape{OT1}{fraktura}{m}{n} {<5> <6> <7> <8> <9> <10>
<11> <12> <13> <14.4> [1.1] eufm10}{}
\DeclareMathAlphabet{\fraktura}{OT1}{fraktura}{m}{n}

\DeclareFontFamily{OT1}{cmfi}{} \DeclareFontShape{OT1}{cmfi}{m}{n}
{<5> <6> <7> <8> <9> <10> <11> <12> <13> <14.4> [0.9] cmfi10}{}
\DeclareMathAlphabet{\cmfi}{OT1}{cmfi}{b}{n}

\DeclareFontFamily{OT1}{cmss}{} \DeclareFontShape{OT1}{cmss}{m}{n}
{<5> <6> <7> <8> <9> <10> <11> <12> <13> <14.4> cmss10}{}
\DeclareMathAlphabet{\cmss}{OT1}{cmss}{m}{n}

\newtheoremstyle{thm}{1.5ex}{1.5ex}{\itshape\rmfamily}{}
{\bfseries\rmfamily}{}{2ex}{}

\newtheoremstyle{rem}{1.3ex}{1.3ex}{\rmfamily}{}
{\itshape}
{} {1.5ex}{}

\newenvironment{proofsect}[1]
{\vskip0.1cm\noindent{\rmfamily\itshape#1.}}{\qed\vspace{0.15cm}}%{\newline\vspace{0.15cm}}

\theoremstyle{thm}
\newtheorem{theorem}{Theorem}[section]
\newtheorem{lemma}[theorem]{Lemma}
\newtheorem{proposition}[theorem]{Proposition}

\newtheorem*{Main Theorem}{Main Theorem.}
\newtheorem{corollary}[theorem]{Corollary}

\theoremstyle{rem}
\newtheorem{remark}{{\itshape Remark}}[]

\numberwithin{equation}{section}

\renewcommand{\section}{\secdef\sct\sect}
\newcommand{\sct}[2][default]{\refstepcounter{section}
\addcontentsline{toc}{section}
{{\tocsection {}{\thesection}{\!\!\!\!#1\dotfill}}{}}
\vspace{0.7cm}
\centerline{ %\large
\scshape\arabic{section}.\ #1} \nopagebreak \vspace{0.2cm}}
\newcommand{\sect}[1]{
\vspace{0.4cm} \centerline{\large\scshape\rmfamily #1}
\vspace{0.2cm}}

\renewcommand{\subsection}{\secdef\subsct\sbsect}
\newcommand{\subsct}[2][default]{\refstepcounter{subsection}
\addcontentsline{toc}{subsection}
{{\tocsection{\!\!}{\hspace{1.2em}\thesubsection}{\!\!\!\!#1\dotfill}}{}}
\nopagebreak \vspace{0.45\baselineskip} {\flushleft\bf
\arabic{section}.\arabic{subsection}~\bf #1.~}
\\*[3mm]\noindent
\nopagebreak}
\newcommand{\sbsect}[1]{\vspace{0.1cm}\noindent
\textbf{#1.~}\vspace{0.1cm}}

\renewcommand{\subsubsection}{%
\secdef \subsubsect\sbsbsect}
\newcommand{\subsubsect}[2][default]{%
\refstepcounter{subsubsection} 
\addcontentsline{toc}{subsubsection}{{\tocsection{\!\!}
{\hspace{3.05em}\thesubsubsection}{\!\!\!\!#1\dotfill}}{}}
\nopagebreak
\vspace{0.15\baselineskip} \nopagebreak {\flushleft\rmfamily
\itshape\arabic{section}.\arabic{subsection}.\arabic{subsubsection}
\ \rmfamily #1\/.}\ }
\newcommand{\sbsbsect}[1]{\vspace{0.1cm}\noindent
\rmfamily \itshape
\arabic{section}.\arabic{subsection}.\arabic{subsubsection} \
\sffamily #1\/.\ }

\newcommand{\dist}{\operatorname{dist}}

\newcommand{\diam}{\operatorname{diam}}

\newcommand{\Int}{{\text{\rm Int}}}
\newcommand{\Ext}{{\text{\rm Ext}}}

\newcommand{\textd}{\text{\rm d}}

\newcommand{\mini}{{\text{\rm min}}}

\newcommand{\Vol}{{\text{\rm Vol}}}

\renewcommand{\AA}{\mathcal A}
\newcommand{\BB}{\mathcal B}
\newcommand{\CC}{\mathcal C}

\newcommand{\EE}{\mathcal E}
\newcommand{\FF}{\mathcal F}

\newcommand{\LL}{\mathcal L}
\newcommand{\MM}{\mathcal M}

\newcommand{\RR}{\mathcal R}
\newcommand{\CalS}{\mathcal S}

\newcommand{\A}{\mathbb A}

\newcommand{\K}{\mathbb K}

\newcommand{\R}{\mathbb R}

\newcommand{\T}{\mathbb T}

\newcommand{\V}{\mathbb V}

\newcommand{\Z}{\mathbb Z}
\newcommand{\bn}{{\boldsymbol n}}
\newcommand{\bh}{{\boldsymbol h}}
\newcommand{\br}{{\boldsymbol r}}

\newcommand{\bzero}{{\boldsymbol 0}}

\newcommand{\Deltac}{\Delta_{\text{\rm c}}}
\newcommand{\lamc}{\lambda_{\text{\rm c}}}
\newcommand{\betac}{\beta_{\text{\rm c}}}

\newcommand{\ssP}{\cmss P}%{\text{\rm P}}}
\newcommand{\dH}{{\text{\it d}_{\text{\rm H}}}}
\newcommand{\mstar}{m^\star}

\newcommand{\ext}{{\text{\rm ext}}}
\newcommand{\mtchoice}[1]{\mathchoice%
	{\text{\fontsize{9}{9}\selectfont #1}}
	{\text{\fontsize{9}{9}\selectfont #1}}
	{\text{\fontsize{6}{6}\selectfont #1}}
	{\text{\fontsize{6}{6}\selectfont #1}}}

\newcommand\IN{{\mtchoice{\rm INT}}}
\newcommand\OUT{{\mtchoice{\rm EXT}}}

\renewcommand{\eqref}[1]{(\ref{#1})}
\newcommand{\twoeqref}[2]{(\ref{#1}--\ref{#2})}

\newcommand{\oExxt}{\Ext^\circ}
\newcommand{\oInnt}{\Int^\circ}
\newcommand{\Exxt}{\Ext}%{\overline{\Ext}}
\newcommand{\Innt}{\Int}%{\overline{\Int}}

\newcommand{\frakS}{\mathfrak S}
\newcommand{\frakP}{\mathfrak P}
\newcommand{\frakM}{\mathfrak M}
\newcommand{\frakb}{\mathfrak b}

\begin{document}

\title[Droplet formation in the 2D Ising model]%, \version]
{\fontsize{16}{20}\selectfont Critical region for droplet formation\\in the
two-dimensional Ising model}

\smallskip

\author[M.~Biskup, L.~Chayes and R.~Koteck\'y]%, \version]
{Marek~Biskup,${}^1\,$ Lincoln~Chayes${}^1$ and$\,$
Roman~Koteck\'y${}^2$}
\maketitle

\vspace{-2mm} \centerline{${}^1$\textit{Department of Mathematics,
UCLA, Los Angeles, California, USA}} %\vspace{-4mm}
\centerline{${}^2$\textit{Center for Theoretical Study, Charles
University, Prague, Czech Republic}}

\thispagestyle{empty}

\vspace{4mm}
\begin{quote}
\footnotesize {\bf Abstract:} 
We study the formation/dissolution of equilibrium droplets in finite systems at parameters 
corresponding to phase coexistence. Specifically, we consider 
the 2D Ising model in volumes of size~$L^2$, inverse temperature~$\beta>\betac$ 
and overall magnetization  conditioned to take the value 
$\mstar L^2-2\mstar v_L$, where~$\betac^{-1}$ 
is the critical temperature,~$\mstar=\mstar(\beta)$ is the spontaneous magnetization
and $v_L$ is a sequence of positive numbers. 
We find that the critical scaling for droplet formation/dissolution is when~$v_L^{3/2} L^{-2}$ 
tends to a definite limit. Specifically, we identify a dimensionless parameter~$\Delta$, 
proportional to this limit, a non-trivial critical value~$\Deltac$ and a function~$\lambda_\Delta$ 
such that the following holds: For~$\Delta<\Deltac$, there are no droplets beyond~$\log L$ scale, 
while for~$\Delta>\Deltac$, there is a single, Wulff-shaped droplet containing 
a fraction~$\lambda_\Delta\ge\lamc=2/3$ of the magnetization deficit and 
there are no other droplets beyond the scale of~$\log L$. 
Moreover,~$\lambda_\Delta$ and~$\Delta$ are related via a universal equation that apparently 
is independent of the details of the system.
\end{quote}

\bigskip\bigskip
\setcounter{tocdepth}{3}

\footnotesize

\begin{list}{}
{\setlength{\topsep}{0in}\setlength{\leftmargin}{0.34in}\setlength{\rightmargin}{0.5in}}
\item[]
\tableofcontents
\end{list}
\vspace{-1.2cm}

\begin{list}{}
{\setlength{\topsep}{0in}\setlength{\leftmargin}{0.5in}\setlength{\rightmargin}{0.5in}}
\item[]
\hskip-0.01in
\hbox to 2.3cm{\hrulefill}
\item[]
{\fontsize{8.6}{8.6}\selectfont\copyright\,\,\,Copyright rests with the authors. Reproduction
of the entire article for non-commercial purposes
is permitted without charge.\vspace{2mm}}
\end{list}
\normalsize

\section{Introduction}
\vspace{-2mm}
\subsection{Motivation}
\label{sec1.1}\noindent
The connection between microscopic interactions and pure-phase 
(bulk) thermodynamics has been understood at a mathematically sophisticated level for many years.  
However, an analysis of systems at phase coexistence which contain droplets
has begun only recently. Over a century ago, 
Curie~\cite{Curie}, Gibbs~\cite{Gibbs} and 
Wulff~\cite{Wulff} derived from surface-thermodynamical considerations that a single droplet 
of a particular shape---the \emph{Wulff shape}---will appear in systems 
that are forced to exhibit a fixed excess of a minority phase. 
A mathematical proof of this fact starting from a system defined on the microscopic scale 
has been given in the context of percolation and Ising systems, first in dimension 
$d=2$~\cite{ACC,DKS} and, more recently, 
in all dimensions $d\ge3$ \cite{Cerf,Bodineau,Cerf-Pisztora}.
Other topics related to the droplet shape have intensively been studied:
Fluctuations of a contour 
line \cite{Campanino-Chayes2,Campanino-Ioffe,Campanino-Ioffe-Velenik,Dobrushin-Hryniv, 
Alexander,Hryniv-Kotecky},
wetting phenomena~\cite{Pf-Velenik2} and Gaussian
fields near a ``wall'' \cite{BenArous-Deuschel,Bolthausen-Ioffe,Dunlop-Magnen}.
See~\cite{BIV} for a summary of these results and comments on the (recent) 
history of these developments.

The initial stages of the rigorous ``Wulff construction'' program have focused 
on systems in which the droplet subsumes a finite fraction of the available volume.
Of no less interest is the situation when the excess represents only a vanishing fraction 
of the total volume. In~\cite{DS}, substantial progress has been made on these questions 
in the context of the Ising model at low temperatures. 
Subsequent developments~\cite{Pfister,Pf-Velenik, Ioffe1,Ioffe2} have allowed the extension, 
in~$d=2$, of the aforementioned results up to the critical point~\cite{Bob+Tim}. 
Specifically, what has so far been shown is as follows: For two-dimensional volumes~$\Lambda_L$ 
of side~$L$ and~$\delta>0$ arbitrarily small, if the magnetization deficit 
exceeds~$L^{4/3+\delta}$,  then a Wulff droplet accounts, pretty much, for all the deficit, 
while if the magnetization deficit  is bounded by~$L^{4/3-\delta}$, 
there are no droplets beyond the scale of~$\log L$. 
The preceding are of course asymptotic statements that hold 
with probability tending to one as~$L\to\infty$.

The focus of this paper is the intermediate regime, which has not yet received 
appropriate attention. Assuming the magnetization deficit divided by~$L^{4/3}$ 
tends to a definite limit, we define a dimensionless parameter, denoted by~$\Delta$, 
which is proportional to this limit. (A precise definition of~$\Delta$ is provided 
in \eqref{Delta-lim}.) Our principal result is as follows: There is a critical value~$\Deltac$ 
such that for~$\Delta<\Deltac$, there are no large droplets (again, nothing beyond~$\log L$ scale), 
while for~$\Delta>\Deltac$, there is a single, large droplet of 
a diameter of the order~$L^{2/3}$. 
However, in contrast to all situations that have previously been analyzed, 
this large droplet only accounts for a finite fraction,~$\lambda_\Delta<1$, 
of the magnetization deficit, which, in addition, does \emph{not} tend to zero 
as~$\Delta\downarrow\Deltac$! (Indeed,~$\lambda_\Delta\downarrow \lamc$, with~$\lamc=2/3$.) 
Whenever the droplet appears, 
its interior is representative of the minus phase, its shape is close to 
the optimal (Wulff) shape and its volume 
is tuned to contain the~$\lambda_\Delta$-fraction of the deficit magnetization. 
Furthermore, for all values of~$\Delta$, there is at most one droplet of size~$L^{2/3}$ 
and nothing else beyond the scale~$\log L$. 
At~$\Delta=\Deltac$ the situation is not completely resolved. 
However, there are only two possibilities: Either there is one droplet of linear size~$L^{2/3}$ or no droplet at all.

The above transition is the result of a competition between two mechanisms 
for coping with a magnetization deficit in the system: Absorption of the deficit by the ambient 
\emph{fluctuations} or the formation of a \emph{droplet}. The results obtained in~\cite{DKS,DS}
and~\cite{Bob+Tim}  deal with the situations when one of the two mechanisms completely dominates the
other.  As is seen by a simple-minded comparison of the exponential costs of the two mechanisms,~$L^{4/3}$
is the only conceivable scaling of the magnetization deficit where these are able to coexist.  
(This is the core of the heuristic approach outlined
in~\cite{Binder-Kalos,Neuhaus-Hager} and~\cite{Binder}, see also~\cite{BCK-comment,Binder-reply}.) 
However, at the point where the droplets first appear,
one can envision alternate scenarios involving complicated fluctuations and/or a multitude of droplets with
effective interactions ranging across many scales. To rule out such possibilities it is necessary to
demonstrate the absence of these ``intermediate-sized'' droplets and the insignificance---or absence---of
large fluctuations. This was argued on a heuristic level in~\cite{BCK} and will be proven rigorously~here.

Thus, instead of blending into each other through 
a series of intermediate scales, 
the droplet-dominated and the fluctuation-dominated regimes meet---literally---at a single point. 
Furthermore, all essential system dependence is encoded 
into one dimensionless parameter~$\Delta$ and the transition between 
the Gaussian-dominated and the droplet-dominated regimes is thus characterized by a \emph{universal}
constant~$\Deltac$. In addition, the relative fraction~$\lambda_\Delta$ of the deficit ``stored'' 
in the droplet depends on~$\Delta$ via a \emph{universal} equation which is 
apparently independent of the details of the system~\cite{BCK}.
At this point  we would like to stress that, even though the rigorous results presented here 
are restricted to the
case of the two-dimensional Ising model, we expect that their validity can be extended to a much larger
class of models and the universality of the dependence on~$\Delta$ will become the subject of a
 \emph{mathematical} statement.
Notwithstanding the rigorous analysis, this universal setting offers the possibility of fitting
experimental/numerical data from a variety of systems onto a single curve.

A practical understanding of how droplets disappear is by no means an esoteric issue. 
Aside from the traditional, i.e., three-dimensional, setting, 
there are experimental realizations which are effectively 
two-dimensional (see~\cite{Sethna} and references therein). 
Moreover, there are purported applications of Ising systems undergoing 
``fragmentation'' in such diverse areas as nuclear physics and adatom formation~\cite{G}.
From the perspective of statistical physics, perhaps more important are the investigations 
of small systems at parameter values corresponding to a first order transition in the bulk. 
In these situations, non-convexities appear 
in finite-volume thermodynamic functions~\cite{G,PH,Chinese,Kosterlitz}, 
which naturally suggest the appearance of a droplet. 
Several papers have studied the disappearance of droplets and reported intriguing
finite-size characteristics~\cite{PS,PH,MS,Sethna,Binder-Kalos,Neuhaus-Hager,Binder}. 
It is hoped that the results established here will shed some light in these situations.

\subsection{The model}
\label{sec1.2}\noindent 
The primary goal of this paper is a detailed description of the above droplet-formation phenomenon
in the Ising model.
In general dimension, this system is defined by the formal Hamiltonian
\begin{equation}
\label{Ham} \mathscr{H}=-\sum_{\langle x,y\rangle}\sigma_x\sigma_y,
\end{equation}
where~$\langle x,y\rangle$ denotes a nearest-neighbor pair on~$\Z^d$ and where 
$\sigma_x\in\{-1,+1\}$ denotes an Ising spin. 
To define the Hamiltonian in a finite volume~$\Lambda\subset\Z^d$, we use~$\partial\Lambda$ 
to denote the external boundary 
of~$\Lambda$, $\partial\Lambda=\{x\notin \Lambda: \text{ there exists a bond } 
\langle x,y\rangle \text{ with } y\in\Lambda\}$,
 fix a collection of boundary spins 
$\sigma_{\partial\Lambda}=(\sigma_x)_{x\in\partial\Lambda}$ and restrict the sum 
in \eqref{Ham} to bonds~$\langle x,y\rangle$ such that~$\{x,y\}\cap\Lambda\ne\emptyset$. 
We denote this finite-volume Hamiltonian
by~$\mathscr{H}_\Lambda(\sigma_\Lambda,\sigma_{\partial\Lambda})$.  The special choices of the boundary
configurations such that~$\sigma_x=+1$, resp., 
$\sigma_x=-1$ for all~$x\in\partial\Lambda$ will be referred to as plus, resp., 
minus boundary conditions.

The Hamitonian gives rise to the concept of a finite-volume \emph{Gibbs measure}
(also known as \emph{Gibbs state}) which is a measure assigning each configuration 
$\sigma_\Lambda=(\sigma_x)_{x\in\Lambda}\in\{-1,+1\}^\Lambda$ the probability
\begin{equation}
\label{Gibb}
P_\Lambda^{\sigma_{\partial\Lambda},\beta}(\sigma_\Lambda)=
\frac{e^{-\beta \mathscr{H}_\Lambda(\sigma_\Lambda,\sigma_{\partial\Lambda})}}
{Z_\Lambda^{\sigma_{\partial\Lambda}}(\beta)}.
\end{equation}
Here~$\beta\ge0$ denotes the inverse temperature,~$\sigma_{\partial\Lambda}$ is 
an arbitrary boundary configuration and $Z_\Lambda^{\sigma_{\partial\Lambda}}(\beta)$ 
is the partition function. Most of this work will concentrate on squares of~$L\times L$ sites, 
which we will denote by~$\Lambda_L$, and the plus boundary conditions. 
In this case we denote the above probability by~$P_L^{+,\beta}(-)$ and 
the associated expectation by~$\langle-\rangle_L^{+,\beta}$. 
As the choice of the signs in \twoeqref{Ham}{Gibb} 
indicates, the measure~$P_L^{+,\beta}$ with~$\beta>0$ tends to favor alignment of neighboring spins with 
an excess of plus spins over minus spins.

\begin{remark}
As is well known, the Ising model is equivalent to a model of a lattice gas where 
at most one particle is allowed to occupy each site. In our case, 
the sites occupied by a particle are represented by minus spins, 
while the plus spins correspond to the sites with no particles. 
In the particle distribution induced by~$P_L^{+,\beta}$, the total number of particles 
is not fixed; hence, we will occasionally refer to this measure as the ``grand canonical'' ensemble. 
On the other hand, if the number of minus spins is fixed 
(by conditioning on the total magnetization, see Section~\ref{sec1.3}), 
the resulting measure will sometimes be referred to as the ``canonical'' ensemble.
\end{remark}

The Ising model has been studied very extensively by 
mathematical physicists in the last~20-30 years 
and a lot of interesting facts have been rigorously established. 
We proceed by listing the properties of the \emph{two-dimensional} model 
which will ultimately be needed in this paper. For general overviews of various aspects 
mentioned below we refer to, e.g.,~\cite{Georgii,Simon,GHM,BIV}. 
The readers familiar with the background (and the standard notation) should feel free to skip 
the remainder of this section and go directly to Section~\ref{sec1.3} where we discuss 
the main results of the present paper.

\smallskip$\bullet$\ \ \textit{Bulk properties. }
For all~$\beta\ge0$, the measure~$P_L^{+,\beta}$ has a unique infinite volume
 (weak) limit~$P^{+,\beta}$ which is a translation-invariant, ergodic, extremal Gibbs state 
for the interaction \eqref{Ham}. 
Let~$\langle-\rangle^{+,\beta}$  denote the expectation with respect to~$P^{+,\beta}$.
The persistence of the plus-bias in the thermodynamic limit, 
characterized by the~\emph{magnetization}
\begin{equation}
m^\star(\beta)=\langle\sigma_0\rangle^{+,\beta},
\end{equation} 
marks the region of phase coexistence in this model.
Indeed, there is a non-trivial critical value 
$\betac\in(0,\infty)$---known~\cite{Onsager,Kaufman-Onsager,AM,BGJS} 
to satisfy~$e^{2\betac}=1+\sqrt2$---such that for 
$\beta>\betac$, we have~$m^\star(\beta)>0$ and there are 
multiple  infinite-volume Gibbs states, while for~$\beta\le\betac$, 
the magnetization vanishes and there is a unique infinite-volume Gibbs state 
for the interaction~\eqref{Ham}. 
Further,  using~$\langle A;B\rangle^{+,\beta}$ to denote the truncated correlation function 
$\langle AB\rangle^{+,\beta}-\langle A\rangle^{+,\beta}\langle B\rangle^{+,\beta}$,
the magnetic \emph{susceptibility}, defined by
\begin{equation}
\chi(\beta)=\sum_{x\in\Z^2}\langle\sigma_0;\sigma_x\rangle^{+,\beta},
\end{equation}
is finite for all~$\beta>\betac$, 
see~\cite{CCS,Schonmann-Shlosman}. By the GHS or FKG inequalities, 
we have~$\chi(\beta)\ge1-m^\star(\beta)^2>0$ for all~$\beta\in[0,\infty)$.

\begin{figure}[t]
\refstepcounter{obrazek}
\ifpdf \centerline{\includegraphics[width=2.7truein]{contours.pdf}}
\else
\centerline{\epsfxsize=2.7truein \epsffile{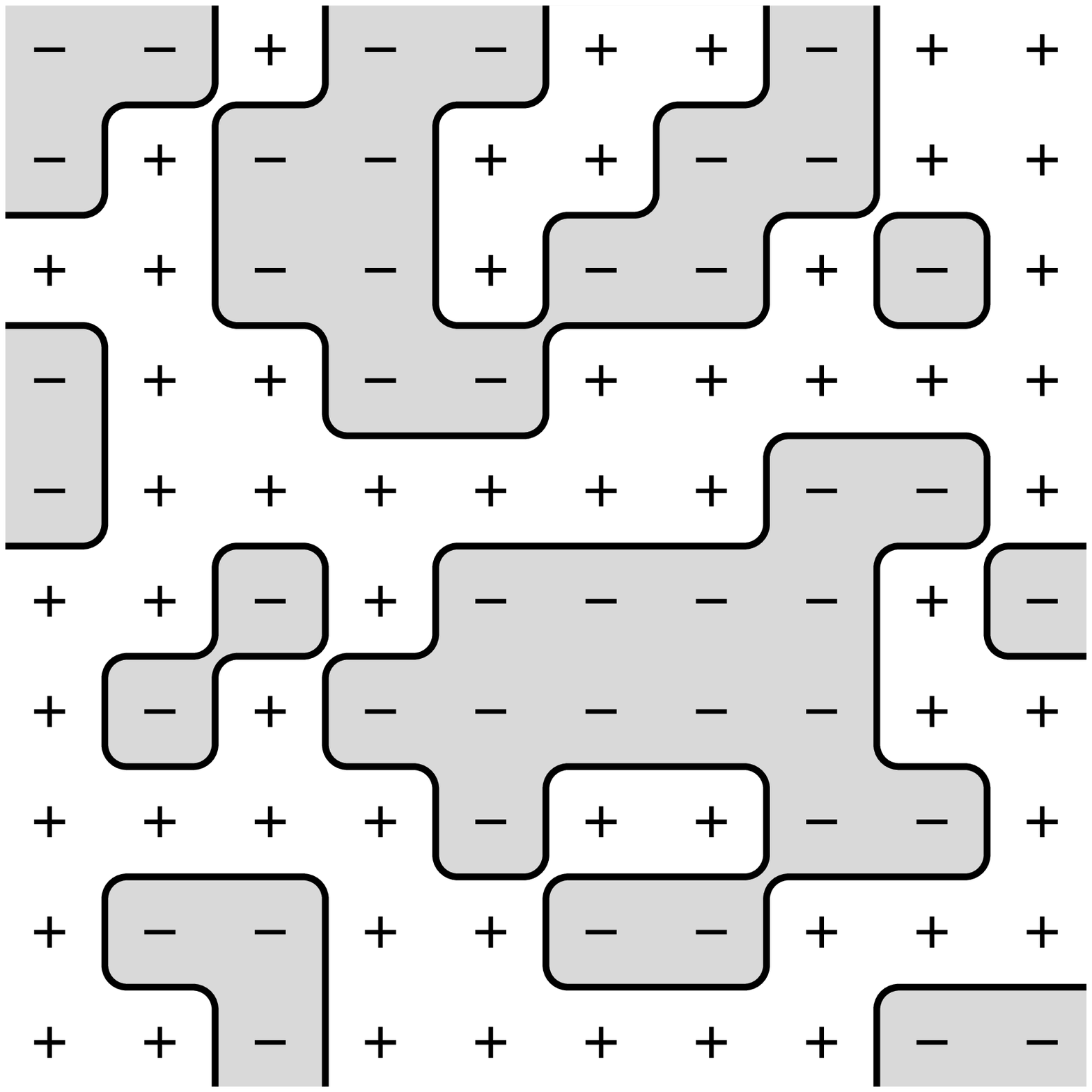}} 
\fi
\bigskip
\begin{quote}
%\small
\fontsize{9.5}{7}\selectfont
{\sc Figure~\theobrazek.\ }
\label{fig1}
An example of an Ising spin configuration and its associated Peierls' contours. In general, 
a contour consists of a string of dual lattice bonds that bisect a direct bond between 
a plus spin and a minus spin. When four such dual bonds meet at a single (dual) lattice site, 
an ambiguity is resolved by applying the south-east/north-west rounding rule. 
(The remaining corners are rounded just for \ae sthetic reasons.) 
The shaded areas correspond to the part of~$V(\gamma)$ occupied by the minus spins.
\normalsize
\end{quote}
\end{figure}

\smallskip$\bullet$\ \ \textit{Peierls' contours. }
Our next requisite item is a description of the Ising configurations in terms of
\textrm{Peierls' contours}.
Given an Ising configuration in~$\Lambda$ 
with plus boundary conditions, we consider the set of dual bonds intersecting direct bonds 
that connect a plus spin with a minus spin. 
These dual bonds will be assembled into contours as follows: 
First we note that only an even number of dual bonds meet at each site of the dual lattice. 
When two bonds meet at a single dual site, we simply connect them. 
When four bonds are incident with one dual lattice site, we apply the 
rounding rule  ``south-east/north-west'' to resolve the ``cross'' into two curves ``bouncing'' off each other 
(see, e.g.,~\cite{DKS,Pf-Velenik} or Figure~\ref{fig1}). 
Using these rules consistently,  the aforementioned set of dual bonds decomposes 
into a set of non self-intersecting polygons with rounded corners. 
These are our \emph{contours}. 

Each contour~$\gamma$ is a boundary of a bounded subset of~$\R^2$, 
which we denote by~$V(\gamma)$. 
We will also need a symbol for the set of sites in the interior of~$\gamma$; 
we let~$\V(\gamma)=V(\gamma)\cap\Z^2$. The \emph{diameter} 
of a contour~$\gamma$ is defined as the diameter of the 
set~$V(\gamma)$  in  the~$\ell_2$-metric on~$\R^2$. 
In the thermodynamic interpretation used 
in Section~\ref{sec1.1}, contours represent microscopic boundaries of droplets.
The advantage of the contour language is that it permits the identification of a sharp boundary 
between two phases; the disadvantage is that, in order to study the typical shape 
(and other properties) of large droplets, one has to first resum over small fluctuations of this boundary.

%%%% COUNTER FOR FIGURES
\newcounter{obrazek}

\smallskip$\bullet$\ \ \textit{Surface tension. }
In order to study droplet equilibrium, we need to introduce the concept 
of microscopic surface tension. Following~\cite{ACC,Pfister}, on~$\Z^2$ 
we can conveniently use \emph{duality}. 
Given a~$\beta>\betac$, let~$\beta^*=\frac12\log\coth\beta$ denote the \emph{dual temperature}. 
For any~$(k_1,k_2)\in\Z^2$ 
and~$k=(k_1^2+k_2^2)^{1/2}$, 
let~$\bn=(k_1/k,k_2/k)\in\CalS_1=\{x\in\R^2\colon\Vert x\Vert=1\}$. 
(Here $\Vert x\Vert$ is the Euclidean norm of $x$.)
Then the limit
\begin{equation}
\label{tau-beta-n}
\tau_\beta(\bn)=\lim_{N\to\infty}\frac1{Nk}
\log\langle\sigma_0 \sigma_{Nk\bn}\rangle^{+,\beta^*},
\end{equation}
where~$Nk\bn=(k_1N,k_2N)\in\Z^2$, exists
independently of what integers~$k_1$ and~$k_2$ we chose to represent~$\bn$ and defines a function 
on a dense subset of~$\CalS_1$. 
It turns out that this function can be continuously extended to all~$\bn\in\CalS_1$.
We call the resulting quantity~$\tau_\beta(\bn)$ the \emph{surface tension} in direction~$\bn$ 
at inverse temperature~$\beta$. 
As is well known, 
$\bn\mapsto\tau_\beta(\bn)$ is invariant under rotations of~$\bn$ by integer multiples of~$\frac\pi2$ 
and~$\tau_\mini=\inf_{\bn\in\CalS_1}\tau_\beta(\bn)
>0$ for all~$\beta>\betac$
\cite{Pfister}. 
Informally, the quantity~$\tau_\beta(\bn)N$ represents the statistical-mechanical cost 
of a (fluctuating) contour line connecting two sites 
at distance~$N$ on a straight line with direction (or normal vector)~$\bn$.

\begin{remark}
Our definition of the surface tension differs from the standard  definition by a factor of~$\beta^{-1}$.
In particular, the physical 
units of~$\tau_\beta$ are length$^{-1}$ rather than energy$\times$length$^{-1}$. The present definition
eliminates the need for an explicit occurrence of~$\beta$ in many expressions 
throughout this paper and, as such, is notationally more convenient.
\end{remark}

\smallskip$\bullet$\ \ \textit{Surface properties. }
On the level of macroscopic thermodynamics, it is obvious that when a droplet 
of the minority phase is present in the system, it is pertinent to minimize the total surface cost. 
By our previous discussion, the 
cost per unit length 
is given by the surface tension~$\tau_\beta(\bn)$. 
Thus, one is naturally led to the functional~$\mathscr{W}_\beta(\gamma)$ 
that assigns the number
\begin{equation}
\label{Wbeta}
\mathscr{W}_\beta(\gamma)=\int_\gamma\tau_\beta(\bn_t)
\textd t
\end{equation} 
to each rectifiable, closed curve~$\gamma=(\gamma_t)$ in~$\R^2$. 
Here $\bn_t$ denotes the normal vector at~$\gamma_t$.
The goal of the resulting variational problem is to minimize~$\mathscr{W}_\beta(\partial D)$ 
over all~$D\subset\R^2$ with rectifiable boundary subject to the constraint that the volume 
of~$D$ coincides with that of the droplet. 
The classic solution, due to Wulff~\cite{Wulff}, is that~$\mathscr{W}_\beta(\partial D)$ 
is minimized by the shape
\begin{equation}
\label{1.7}
D_W=\bigl\{\br\in\R^2\colon \br\cdot\bn\le\tau_\beta(\bn),\,\bn\in\CalS_1\bigr\}
\end{equation}
rescaled to contain the appropriate volume. (Here~$\br\cdot\bn$ denotes the dot product in~$\R^2$.) 
We will use~$W$ to denote the shape~$D_W$ 
scaled to have a \emph{unit} (Lebesgue) volume.
It follows from \eqref{1.7} that~$W$ is a convex set in~$\R^2$. 
We define
\begin{equation}
w_1(\beta)=\mathscr{W}_\beta(\partial W)
\end{equation}
and note that~$w_1(\beta)>0$ once~$\beta>\betac$.

\smallskip
Our preliminary arsenal is now complete and we are prepared to discuss the main results.

\subsection{Main results}
\label{sec1.3}\noindent
Recall the notation~$\Lambda_L$ for a square of~$L\times L$ sites in~$\Z^2$.
Consider the Ising model in volume~$\Lambda_L$ with plus boundary condition 
and inverse temperature~$\beta$. Let us define the total magnetization 
(of a configuration $\sigma$) in~$\Lambda_L$ by the formula
\begin{equation}
\label{cond-prob}
M_L=\sum_{x\in\Lambda_L}\sigma_x.
\end{equation}
Let~$(v_L)_{L\ge1}$ be a sequence of positive numbers, with~$v_L\to\infty$ as~$L\to\infty$, 
such that~$\mstar\,|\Lambda_L|-2\mstar\, v_L$ is an allowed value of~$M_L$ for all~$L\ge1$. 
Our first result concerns the decay rate of the probability that 
$M_L=\mstar\,|\Lambda_L|-2\mstar\, v_L$ in the ``grand canonical'' ensemble~$P_L^{+,\beta}$:

\begin{theorem}
\label{LDP-thm} Let~$\beta>\betac$ and let~$\mstar\,=\mstar\,(\beta)$,~$\chi=\chi(\beta)$, 
and~$w_1=w_1(\beta)$ be as above. Suppose that the limit
\begin{equation}
\label{Delta-lim} \Delta=2\frac{(\mstar)^2}{\chi w_1}
\,\lim_{L\to\infty}\frac{\,v_L^{3/2}}{|\Lambda_L|}
\end{equation}
exists with~$\Delta\in(0,\infty)$. Then
\begin{equation}
\label{LDP} \lim_{L\to\infty}\frac1{\sqrt{v_L}}\,\log
P_L^{+,\beta}\bigl(M_L=\mstar\,|\Lambda_L|-2\mstar\, v_L\bigr)=
-w_1\inf_{0\le \lambda\le1}\Phi_\Delta(\lambda),
\end{equation}
where
\begin{equation}
\label{Phi} 
\Phi_\Delta(\lambda)=\sqrt \lambda+\Delta (1-\lambda)^2, \qquad 0\le
\lambda\le 1.
\end{equation}
\end{theorem}

\smallskip

The proof of Theorem~\ref{LDP-thm} is a direct consequence of Theorems~\ref{lowerbound} 
and~\ref{upperbound}; the actual proof comes in Section~\ref{sec5}. 
We proceed with some remarks:

\begin{remark}
Note that, by our choice of the deviation scale, the term~$\mstar(\beta)|\Lambda_L|$ 
can be replaced by the mean value~$\langle M_L\rangle_L^{+,\beta}$ in all formulas; 
see Lemma~\ref{L:M-fluctuations} below. 
The motivation for introducing the factor~``$2\mstar$'' on the left-hand-side of \eqref{LDP} 
is that then~$v_L$ represents the volume of a droplet that must be created 
in order to achieve the required value of the overall magnetization
(provided the magnetization outside, resp., inside the droplet is $\mstar$, resp., $-\mstar$). 
\end{remark}

\begin{remark}
The quantity~$\lambda$ that appears in \twoeqref{LDP}{Phi} represents the \emph{trial fraction} 
of the deficit magnetization which might go into a large-scale droplet. (So, by our convention, 
the volume of such a droplet is just~$\lambda v_L$.) 
The core of the proof of Theorem~\ref{LDP-thm}, roughly speaking, 
is that the probability of seeing a droplet of this size tends to zero as 
$\exp\{-w_1\sqrt{v_L}\Phi_\Delta(\lambda)\}$. 
Evidently, a large deviation principle for the size of such a droplet is satisfied 
with rate~$L^{2/3}$ and a rate function proportional to~$\Phi_\Delta$. 
However, we will not attempt to make this statement mathematically rigorous.
\end{remark}

Next we shall formulate our main result on the asymptotic form of typical configurations 
in the ``canonical'' ensemble described by the conditional measure 
$P_L^{+,\beta}(\,\cdot\,|M_L=\mstar\,|\Lambda_L|-2\mstar\, v_L)$. 
For any two sets~$A,B\subset\R^2$, let~$\dH(A,B)$ denote the Hausdorff distance 
between~$A$ and~$B$,
\begin{equation}
\label{Hausd}
\dH(A,B)= \max\bigl\{\sup_{x\in A}\dist(x,B), \sup_{y\in B}\dist(y,A)\bigr\},
\end{equation}
where $\dist(x,A)$ is the Euclidean distance of~$x$ and~$A$.

Our second main theorem is then as follows:

\begin{theorem}
\label{main-result} Let~$\beta>\betac$ and suppose that the limit in \eqref{Delta-lim} exists 
with~$\Delta\in(0,\infty)$. Recall that~$W$ denotes the Wulff shape of a unit volume. 
Given $\varkappa,s,L\in(0,\infty)$,
let~$\AA_{\varkappa,s,L}$ be the event that 
any external contour~$\gamma$ for which~$\diam\gamma\ge s$ must also satisfy 
$\diam\gamma
>\varkappa\sqrt{v_L}$. 
Next, for each~$\epsilon>0$, let~$\BB_{\epsilon,s,L}$ be the event 
that there is at most one external contour~$\gamma_0$ in~$\Lambda_L$ with~$\diam\gamma_0\ge s$ and, 
whenever such a contour~$\gamma_0$ exists, it satisfies the conditions
\begin{equation}
\label{hausbd}
\inf_{z\in\R^2}\dH\bigl(V(\gamma_0),z+\sqrt{|V(\gamma_0)|}\,W\bigr)\le \sqrt{\epsilon v_L}
\end{equation}
and
\begin{equation}
\label{volbd}
\Phi_\Delta\bigl(v_L^{-1}|V(\gamma_0)|\bigr)\le
\inf_{0\le \lambda'\le1}\Phi_\Delta(\lambda')+\epsilon.
\end{equation}
In addition, the event~$\BB_{\epsilon,s,L}$ also requires that 
the magnetization inside~$\gamma_0$ obeys the constraint
\begin{equation}
\label{magbd}
\biggl|\sum_{x\in\V(\gamma_0)}(\sigma_x+\mstar)\biggr|\le\epsilon v_L.
\end{equation}
There exists a constant~$\varkappa_0>0$ such that for each $\zeta>0$ and each $\epsilon>0$ there exist numbers $K_0<\infty$ and $L_0<\infty$ such that
\begin{equation}
\label{1.16}
P_L^{+,\beta}\bigl(\AA_{\varkappa,s,L}\cap\BB_{\epsilon,s,L} \big|M_L=\mstar\,|\Lambda_L|-
2\mstar\, v_L\bigr)\ge1-L^{-\zeta}
\end{equation}
holds provided~$\varkappa\le\varkappa_0$ and~$s=K\log L$ with~$K\ge K_0$ and~$L\ge L_0$, . 
\end{theorem}

Thus, simply put, whenever there is a large droplet in the system, its shape rarely deviates 
from that of the Wulff shape and its volume (in units of~$v_L$) is almost always given 
by a value of~$\lambda$ nearly minimizing~$\Phi_\Delta$. 
Moreover, all other droplets in the system are at most of a logarithmic~size.

\smallskip
Most of the physically interesting behavior of this system is simply 
a consequence of where~$\Phi_\Delta$ achieves its minimum and how this minimum depends 
on~$\Delta$. The upshot, which is stated concisely in Proposition~\ref{prop-Phi} below, 
is that there is a critical value of~$\Delta$, given by
\begin{equation}
\label{Delta-c2} 
\Deltac=\frac12\Bigl(\frac32\Bigr)^{3/2},
\end{equation}
such that if~$\Delta<\Deltac$, then~$\Phi_\Delta$ has the unique minimizer at~$\lambda=0$, 
while for~$\Delta>\Deltac$, the unique minimizer of~$\Phi_\Delta$ is nontrivial. 
More explicitly, for~$\Delta\ne\Deltac$, the function~$\Phi_\Delta$ is minimized by
\begin{equation}
\label{xDelta2}
\lambda_\Delta=\begin{cases}0,\qquad&\text{if }\Delta<\Deltac,\\
\lambda_+(\Delta),\qquad&\text{if }\Delta>\Deltac,
\end{cases}
\end{equation}
where~$\lambda_+(\Delta)$ is the maximal positive solution to the equation
\begin{equation}
%\label{}
4\Delta\sqrt \lambda(1-\lambda)=1.
\end{equation}
The reason for the changeover is that, as~$\Delta$ increases through~$\Deltac$, 
a local minimum becomes a global minimum, see the proof of Proposition~\ref{prop-Phi}. 
As a consequence, the minimizing fraction~$\lambda$ does \emph{not} 
tend to zero as~$\Delta\downarrow\Deltac$; in particular, it tends to~$\lamc=2/3$.

\smallskip
Using the information about the unique minimizer of~$\Phi_\Delta$ for~$\Delta\ne\Deltac$, 
it is worthwhile to reformulate Theorem~\ref{main-result} as follows:

\begin{corollary}
\label{cor-main-result} 
Let~$\beta>\betac$ and suppose that the limit in \eqref{Delta-lim} exists with~$\Delta\in(0,\infty)$. 
Let~$\Deltac$ and~$\lambda_\Delta$ be as in \eqref{Delta-c2} and \eqref{xDelta2}, respectively. 
Let~$K$ be sufficiently large (i.e., $K\ge K_0$, where~$K_0$ is as in Theorem~\ref{main-result}). 
Considering the conditional distribution~$P_L^{+,\beta}(\,\cdot\,|M_L=\mstar\,|\Lambda_L|-2\mstar\, v_L)$,
the following holds with probability tending to one as~$L\to\infty$:
\begin{enumerate}
\item[(1)]
If~$\Delta<\Deltac$, then all contours~$\gamma$ in~$\Lambda_L$ satisfy~$\diam\gamma\le K\log L$.
\item[(2)]
If~$\Delta>\Deltac$, then there is exactly one external
contour~$\gamma_0$ with~$\diam\gamma_0>K\log L$
and all other external contours~$\gamma$ satisfy~$\diam\gamma\le K\log L$. 
Moreover, the unique ``large''  external contour~$\gamma_0$ asymptotically satisfies 
the bounds \twoeqref{hausbd}{magbd} for all~$\epsilon>0$.
In particular,~$|V(\gamma_0)|= v_L(\lambda_\Delta+o(1))$ with probability tending to one as~$L\to\infty$. 
\end{enumerate}
\end{corollary}

We remark that although the situation at~$\Delta=\Deltac$ is not fully resolved, 
we must have either a single large droplet or no droplet at all; i.e., 
the outcome must mimic the case~$\Delta>\Deltac$ or~$\Delta<\Deltac$. 
A better understanding of the case~$\Delta=\Deltac$ will certainly require a more refined analysis; 
e.g., the second-order large-deviation behavior of the measure~$P_L^{+,\beta}(\cdot)$.

\begin{remark}
We note that in the course of this work, the phrase ``$\beta > \betac$'' appears in three disparate meanings.  
First, for~$\beta > \betac$, the magnetization is positive, second, for~$\beta > \betac$, 
the surface tension is positive and third, for~$\beta > \betac$, truncated correlations decay exponentially.  
The facts that the transition temperatures associated with these properties all coincide 
\emph{and} that~$\betac$ is given by the self-dual condition plays no essential role in our arguments.  
Nor are any other particulars of the square lattice really used.  
Thus, we believe that  our results could be extended to other planar lattices without much modification.  
However, in the cases where the coincidence has not yet been (or cannot be) established,
we would need to define ``$\betac$'' so as to satisfy all three criteria.
\end{remark}

\subsection{Discussion and outline}
\label{heuristic}\noindent
The mechanism which drives the droplet formation/dissolution phenomenon described
in the above theorems is not difficult to understand on a heuristic level. 
This heuristic derivation (which applies to all dimensions~$d\ge2$)
has been discussed in detail elsewhere~\cite{BCK}, 
so we will be correspondingly brief. The main ideas are best explained in the context 
of the large-deviation theory for the ``grand canonical'' distribution and, as a matter of fact, 
the actual proof also follows this~path.

Consider the Ising model in the box~$\Lambda_L$ and suppose we wish 
to observe a magnetization deficiency~$\delta M=2\mstar v_L$ from the nominal value 
of~$\mstar|\Lambda_L|$. Of course, this can be achieved in one shot by the formation 
of a Wulff droplet at the cost of about~$\exp\{-w_1\sqrt{v_L}\}$. 
Alternatively, if we demand that this deficiency emerges out of the background fluctuations, 
we might guess on the basis of fluctuation-dissipation arguments that the cost 
would be of the order
\begin{equation}
%\label{}
\exp\Bigl\{-
\frac{(\delta M)^2}{2\text{Var}(M_L)}\Bigr\}\approx
\exp\Bigl\{-
2\frac{(\mstar\,v_L)^2}{\chi|\Lambda_L|}\Bigr\},
\end{equation}
where~$\chi$ is the susceptibility 
and $\text{Var}(M_L)=(\chi+o(1))|\Lambda_L|$ is the variance of~$M_L$ in distribution~$P_L^{+,\beta}$. 
Obviously, the former mechanism dominates when 
$\sqrt{v_L}\ll v_L^2/|\Lambda_L|$, i.e., when~$v_L\gg L^{4/3}$, while the latter dominates 
under the opposite extreme conditions, i.e., when~$v_L\ll L^{4/3}$. 
(These are exactly the regions previously treated in \cite{DS,Bob+Tim} 
where the corresponding statements have been established in full rigor.) 
In the case when~$v_L/L^{4/3}$ tends to a finite limit 
we now find that the two terms are comparable. 
This is the basis of our parameter~$\Delta$ defined in \eqref{Delta-lim}. 

Assuming~$v_L^{3/2}/|\Lambda_L|$ is essentially at its limit, let us instead try a droplet 
of volume~$\lambda v_L$, where~$0\le\lambda\le1$. 
The droplet cost is now reduced to
\begin{equation}
\label{1.22a}
\exp\bigl\{-w_1\sqrt\lambda\sqrt{v_L}\bigr\},
\end{equation} 
but we still need to account for the remaining fraction of the deficiency.
Assuming the fluctuation-dissipation reasoning can still be applied, this is now
\begin{equation}
%\label{}
\exp\Bigl\{-2\frac{(\mstar\,v_L)^2}{\chi|\Lambda_L|}(1-\lambda)^2\Bigr\}=
\exp\bigl\{-w_1\sqrt{v_L}(1-\lambda)^2\Delta\bigr\}.
\end{equation}
Putting these together we find that the total cost of achieving the deficiency 
$\delta M=2\mstar v_L$ using a droplet of volume~$\lambda v_L$ is given in the leading order by
\begin{equation}
\label{1.22}
\exp\bigl\{- w_1\Phi_\Delta(\lambda)\sqrt{v_L}\bigr\}.
\end{equation}
An optimal droplet size is then found by minimizing~$\Phi_\Delta(\lambda)$ over~$\lambda$. 
This is exactly the content of Theorem~\ref{LDP-thm}. 
We remark that even on the level of heuristic understanding, 
some justification is required for the decoupling of these two mechanisms. 
In~\cite{BCK}, we have argued this case on a heuristic level; 
in the present work, we simply provide a complete proof.

The pathway of the proof is as follows: The approximate equalities \twoeqref{1.22a}{1.22} must be proved 
in the form of upper and lower bounds which agree in the~$L\to\infty$ limit. 
(Of course, we never actually have to go through the trouble of establishing 
the formulas involving~$\Phi_\Delta(\lambda)$ for non-optimal values of~$\lambda$.) 
For the lower bound (see Theorem~\ref{lowerbound}) 
we simply shoot for the minimum of~$\Phi_\Delta(\lambda)$: 
We produce a near-Wulff droplet of the desired area and, on the complementary region, 
allow the background fluctuations to account for the rest. 
Here, as a bound, we are permitted to use a contour ensemble with restriction to 
contours of \emph{logarithmic} size which ensures the desired Gaussian behavior.

The upper bound requires considerably more effort. 
The key step is to show that, with probability close to one, 
there are no droplets at any scale larger than~$\log L$ 
or smaller than~$\sqrt{v_L}$. 
Notwithstanding the technical difficulties, 
the result (Theorem~\ref{upperbound}) 
is of independent interest because it applies for all~$\Delta\in(0,\infty)$, 
including the case~$\Delta=\Deltac$. 
Once the absence of these ``intermediate'' contour scales has been established, 
the proof of the main results directly follow.

\smallskip
We finish with a brief outline of the remainder of this paper. In the next section 
we collect the necessary technical statements needed for the proof of both the upper 
and lower bound. Specifically, in Section~\ref{sec2.1} we discuss in detail the minimizers 
of~$\Phi_\Delta$, in Section~\ref{sec2.2} we introduce the concept of skeletons 
and in Section~\ref{sec2.3} we list the needed properties of the logarithmic contour ensemble. 
Section~\ref{sec3} contains the proof of the lower bound, while Section~\ref{sec4} establishes 
the absence of contour on scales between~$\log L$ and the anticipated droplet size. 
Section~\ref{sec5} assembles these ingredients together into the proofs of the main results.

\section{Technical ingredients}
\smallskip\noindent
This section contains three subsections: Section~\ref{sec2.1} presents the 
solution of the variational problem for function~$\Phi_\Delta$ on the right-hand 
side of~\eqref{Phi}, while Sections~\ref{sec2.2} and~\ref{sec2.3} collect the 
necessary technical lemmas concerning the skeleton calculus and the small-contour 
ensemble.
We remark that a variety of closely related results have appeared in 
literature; in particular, in \cite{Bob+Tim} (and the earlier  \cite{DKS,DS,Pfister}).
For completeness, we will provide proofs, but keep them as brief as possible.
Readers familiar with these topics (or who are otherwise uninterested) are 
invited to skip the entire section on a preliminary run-through, 
referring back only for definitions when reading through Sections~\ref{sec3}--\ref{sec5}.

\subsection{Variational problem}
\label{sec2.1}\noindent
Here we investigate the global minima of the function
$\Phi_\Delta$ that was introduced in \eqref{Phi}. Since the general picture
is presumably applicable in 
higher dimensions as well 
(certainly at the level of heuristic arguments,
see~\cite{BCK}), we might as well carry out the analysis in the case of a general
dimension $d\ge2$. For the purpose of this subsection, let
\begin{equation}
\label{Phi-d} \Phi_\Delta(\lambda)=\lambda^{\frac{d-1}d}+\Delta
(1-\lambda)^2,
\qquad 0\le \lambda\le 1.
\end{equation}
We define
\begin{equation}
\label{Phi*-d} \Phi_\Delta^\star=\inf_{0\le \lambda\le1}\Phi_\Delta(\lambda)
\end{equation}
and note that~$\Phi_\Delta^\star>0$ once~$\Delta>0$. 
Let us introduce the~$d$-dimensional 
version of \eqref{Delta-c2},
\begin{equation}
\label{Delta-c-d}
\Deltac=\frac1d\Bigl(\frac{d+1}2\Bigr)^{\!\!\textstyle\frac{d+1}d}.
\end{equation}
The minimizers of~$\Phi_\Delta$ are then characterized as follows:

\begin{proposition}
\label{prop-Phi} 
Let~$d\ge2$ and, for any~$\Delta\ge0$, let~$\frakM_\Delta$ denote the set of 
all global minimizers of~$\Phi_\Delta$ on~$[0,1]$.
Then we have:

(1) If~$\Delta<\Deltac$, then~$\frakM_\Delta=\{0\}$.

(2) If~$\Delta=\Deltac$, then~$\frakM_\Delta=\{0,\lamc\}$, where
\begin{equation} 
\label{lambdac}
\lamc=\frac2{d+1}.
\end{equation}

(3) If~$\Delta>\Deltac$, then~$\frakM_\Delta=\{\lambda_0\}$, where 
$\lambda_0$ is the maximal
positive solution to the equation
\begin{equation}
\label{crit-eq} \frac{2d}{d-1}\Delta\, \lambda^{\frac1d}(1-\lambda)=1.
\end{equation}
\indent\phantom{(3)} In particular,~$\lambda_0>\lamc$.
\end{proposition}

\begin{proofsect}{Proof}
A simple calculation shows that~$\lambda=0$ is always a (one-sided)
local minimum of~$\lambda\mapsto\Phi_\Delta(\lambda)$, while~$\lambda=1$ is always 
a (one-sided) local maximum. Moreover, the stationary points of~$\Phi_\Delta$ in~$(0,1)$ have to satisfy \eqref{crit-eq}.
Consider the quantity
\begin{equation}
%\label{}
q(\lambda)=\frac1\Delta\bigl(1-\tfrac{d}{d-1}\lambda^{1/d}\Phi_\Delta'(\lambda)\bigr)
=\frac{2d}{d-1} \lambda^{1/d}(1-\lambda),
\end{equation}
i.e., $q(\lambda)$  is essentially 
the left-hand side of \eqref{crit-eq}. A simple calculation shows that~$q(\lambda)$ achieves 
its maximal value on~$[0,1]$ at
$\lambda=\lambda_d=\frac1{d+1}$, where it equals~$\Delta_d^{-1}=2d^2(d^2-1)^{-1}(d+1)^{-1/d}$, 
and is strictly increasing for~$\lambda<\lambda_d$ and strictly decreasing 
for~$\lambda>\lambda_d$. On the basis of these observations, 
it is easy to verify the following facts:
\begin{enumerate}
\item[(1)]
For~$\Delta\le\Delta_d$, 
we have~$\Delta q(\lambda)<1$ for all $\lambda\in[0,1]$ (except perhaps at $\lambda=\lambda_d$ when~$\Delta$ equals~$\Delta_d$). Consequently,~$\lambda\mapsto\Phi_\Delta(\lambda)$ is strictly increasing throughout~$[0,1]$. 
In particular,~$\lambda=0$ is the unique global minimum of~$\Phi_\Delta(\lambda)$ in~$[0,1]$.
\item[(2)]
For~$\Delta>\Delta_d$, \eqref{crit-eq}, resp., $\Delta q(\lambda)=1$ has two distinct solutions in $[0,1]$. Consequently,~$\lambda\mapsto\Phi_\Delta(\lambda)$ has two local 
extrema in~$(0,1)$: A local maximum at~$\lambda=\lambda_-(\Delta)$ and a local minimum at
$\lambda=\lambda_+(\Delta)$, where~$\lambda_-(\Delta)$ and
$\lambda_+(\Delta)$ are the minimal and maximal positive solutions to
\eqref{crit-eq}, respectively.
\end{enumerate}
As a simple calculation shows,
the function~$\Delta\mapsto\lambda_+(\Delta)$ is 
strictly increasing on its domain with 
$\lambda_+(\Delta)\sim 1-\frac{d-1}{2d}\frac1\Delta$ as~$\Delta\to\infty$.

In order to decide which of the two previously described local minima ($\lambda=0$
or~$\lambda=\lambda_+(\Delta)$) gives rise to the global minimum, we first
note that, while~$\Phi_\Delta(0)=\Delta$ tends to infinity as~$\Delta\to\infty$, 
the above asymptotics of~$\lambda_+(\Delta)$ shows 
that~$\Phi_\Delta(\lambda_+(\Delta))\to1$ 
as~$\Delta\to\infty$. Hence,~$\lambda_+(\Delta)$ is the unique global minimum of 
$\Phi_\Delta$ once~$\Delta$ is sufficiently large. 
Thus, it remains to show that 
the two local minima interchange their roles at~$\Delta=\Deltac$. To that end we 
compute
\begin{equation}
\frac{\textd}{\textd \Delta}\Phi_\Delta\bigl(\lambda_+(\Delta)\bigr)=
\frac{\partial}{\partial\Delta}\Phi_\Delta\bigl(\lambda_+(\Delta)\bigr)
=\bigl(1-\lambda_+(\Delta)\bigr)^2
<1,
\end{equation}
where we used that~$\lambda_+(\Delta)$ is a stationary point of
$\Phi_\Delta$ to derive the first equality. Comparing this with
$\frac{\textd}{\textd \Delta}\Phi_\Delta(0)=1$, we see that
$\Delta\mapsto\Phi_\Delta(\lambda_+(\Delta))$ increases with~$\Delta$ strictly 
slower than~$\Delta\mapsto\Phi_\Delta(0)$ 
on any finite interval of~$\Delta$'s. 
Hence, there must be a unique value of~$\Delta$  for which 
$\Phi_\Delta(0)$ and~$\Phi_\Delta(\lambda_+(\Delta))$
are exactly equal. An elementary computation
shows that this happens at~$\Delta=\Deltac$, where~$\Deltac$ is
given by~\eqref{Delta-c-d}. This finishes the proof of (1) and
(3); in order to show that also (2) holds, we just need to note
that~$\lambda_+(\Deltac)$ is exactly~$\lamc$ as given in \eqref{lambdac}. 
\end{proofsect}

Proposition~\ref{prop-Phi} allows us to define a quantity
$\lambda_\Delta$ by formula \eqref{xDelta2}, where now~$\lambda_+(\Delta)$ is
the maximal positive solution to \eqref{crit-eq}. Since
$\lim_{\Delta\downarrow\Deltac}\lambda_\Delta=\lamc>0$, 
the function
$\Delta\mapsto\lambda_\Delta$ undergoes a jump at~$\Deltac$.

\subsection{Skeleton estimates}
\label{sec2.2}\noindent
In this section we introduce coarse-grained versions of contours called 
\emph{skeletons}. These objects will be extremely useful whenever 
an upper bound on the probability of large
contours is needed. Indeed, the introduction of skeletons will permit us to
effectively integrate out small fluctuations of contour lines and thus
express the contour weights directly in terms of the surface tension.
Skeletons were first introduced in~\cite{ACC,DKS};
here we use a modified version of the definition from~\cite{Bob+Tim}.

\subsubsection{Definition and geometric properties}
Given a scale~$s>0$, an~$s$-\emph{skeleton} is an~$n$-tuple
$(x_1,\dots,x_n)$ of points on the dual lattice,~$x_i\in(\Z^2)^*$, such that~$n>1$ and
\begin{equation}
\label{2.7cnd}
s\le\Vert x_{i+1}-x_i\Vert\le 2s,\qquad i=1,\dots,n.
\end{equation}
Here~$\Vert\cdot\Vert$ denotes the~$\ell^2$-distance on~$\R^2$
and~$x_{n+1}$ is identified with~$x_1$.
Given a skeleton~$S$, let~$\ssP(S)$ be the closed \emph{polygonal curve} 
in~$\R^2$ induced by~$S$.
We will use~$|\ssP(S)|$ to denote the total 
length of~$\ssP(S)$, in accord with our general
notation for the length of curves.

A contour~$\gamma$ is called \emph{compatible} with an
$s$-skeleton~$S=(x_1,\dots,x_n)$, if
\begin{enumerate}
\item[(1)]~$\gamma$, viewed as a simple closed path on~$\R^2$,
passes through all sites 
$x_i$,~$i=1,\dots,n$ in the corresponding order.
\item[(2)] 
$\dH(\gamma,\ssP(S))\le s$, where~$\dH$ is the Hausdorff distance
\eqref{Hausd}.
\end{enumerate}
We write~$\gamma\sim S$ if~$\gamma$ and~$S$ are compatible.
For each 
configuration~$\sigma$, we let~$\Gamma_s(\sigma)$ be the
set of all~$s$-\emph{large} contours~$\gamma$ in~$\sigma$; 
namely all~$\gamma$ in~$\sigma$ for which there is 
an~$s$-skeleton~$S$ such that~$\gamma\sim  S$.
Given a set of~$s$-skeletons~$\frakS=(S_1,\dots,S_m)$, we say
that a configuration~$\sigma$ is \emph{compatible} with
$\frakS$, if~$\Gamma_s(\sigma)=(\gamma_1,\dots,\gamma_m)$ and
$\gamma_k\sim S_k$ for all~$k=1,\dots,m$. We will write~$\sigma\sim\frakS$ 
to denote that~$\sigma$ and~$\frakS$ are compatible.

It is easy to see that~$\Gamma_s(\sigma)$ actually consists 
of all contours~$\gamma$ of the configuration~$\sigma$
such that~$\diam\gamma\ge s$. 
Indeed,~$\diam\gamma\ge s$ for every~$\gamma\in\Gamma_s(\sigma)$
by the conditions (1) and \eqref{2.7cnd} above.  
On the other hand, for any~$\gamma$ with~$\diam\gamma\ge s$, 
we will construct an~$s$-skeleton by the following 
procedure: Regard~$\gamma$ as a closed non-self-intersecting curve, 
$\gamma=(\gamma_t)_{0\le t\le1}$, where~$\gamma_0$ is chosen so that 
$\sup_{x\in\gamma}\Vert x-\gamma_0\Vert\ge s$.  
Then we let~$x_1=\gamma_0$ and~$x_2=\gamma_{t_2}$, 
where~$t_2=\inf\{t>0\colon \Vert\gamma_t-\gamma_0\Vert\ge s\}$. Similarly, 
if~$t_j$ has been defined and~$x_j=\gamma_{t_j}$, we let 
$x_{j+1}=\gamma_{t_{j+1}}$, where~$t_{j+1}=\inf\{t\in(t_j,1]\colon \Vert\gamma_t-
\gamma_{t_j}\Vert\ge s\}$. Note that this definition ensures that 
\eqref{2.7cnd} as well as the conditions (1) and (2) hold.
The consequence of this construction is that,
via the equivalence relation $\sigma\sim\frakS$, 
the set of all skeletons 
induces a \emph{covering} of the set of all spin configurations.

\begin{remark}
The reader familiar with \cite{DKS,Bob+Tim} will notice that we explicitly keep 
the stronger condition~(1) from \cite{DKS}. Without the requirement that contours 
pass through the skeleton points in the given order, Lemma~\ref{L:skeleti} and, 
more importantly, Lemma~\ref{lemma-skeletonUB} below would fail to 
hold.
\end{remark}

\smallskip
Next we will discuss some subtleties of the geometry of the skeletons stemming from 
the fact that the corresponding polygons (unlike contours) may have 
self-intersections. We will stay rather brief; 
a detailed account of the topic can be found in~\cite{DKS}. 

We commence with a few geometric
definitions: Let~$\frakP=\{\ssP_1,\dots,\ssP_k\}$ denote a finite
collection of polygonal curves. Consider a smooth self-avoiding path~$\LL$ from a 
point~$x$ to~$\infty$  
that is generic  with respect to the polygons from~$\frakP$
(i.e., the path~$\LL$ has a finite number of intersections with
each~$\ssP_j$ and this number does not change under small perturbations of~$\LL$).
Let~$\#(\LL\cap\ssP_j)$ be the number of intersections of~$\LL$ with~$\ssP_j$. 
Then we \emph{define}~$V(\frakP)\subset\R^2$ to 
be the set of 
points~$x\in\R^2$ such that the total number of intersections,
$\sum_{j=1}^n\#(\LL\cap\ssP_j)$, is odd for any path~$\LL$ 
from~$x$ to~$\infty$ with the above properties. We will use 
$|V(\frakP)|$ to denote the area of~$V(\frakP)$.

If~$\frakP$ happens to be a collection of skeletons,~$\frakP=\frakS$, the relevant set will be~$V(\frakS)$.
If~$\frakP$ happens to be a collection of Ising contours,~$\frakP=\Gamma$, the associated~$V(\Gamma)$ can be thought of as a union of plaquettes centered 
at sites of~$\Z^2$; we will use~$\V(\Gamma)=V(\Gamma)\cap\Z^2$ 
to denote the relevant set of sites. It is clear that if~$\Gamma$ are the contours 
associated with a spin configuration~$\sigma$ in~$\Lambda$ and the plus boundary 
condition on~$\partial\Lambda$, then~$\V(\Gamma)$ are exactly the sites 
$x\in\Lambda$ where~$\sigma_x=-1$. We proceed by listing a few important estimates 
concerning compatible collections of contours and their associated skeletons:

\begin{lemma}
\label{L:triv}
There is a finite geometric constant~$g_1$ such that if~$\Gamma$ is a collection of contours and~$\frakS$ is a collection of 
$s$-skeletons with~$\Gamma\sim\frakS$, then
\begin{equation}
\label{triv-s}
\sum_{\gamma\in\Gamma}|\gamma|\le g_1s\sum_{S\in\frakS}\bigl|\ssP(S)\bigr|.
\end{equation}
In particular, if~$\diam\gamma\le\varkappa$ for all~$\gamma\in\Gamma$,
then we also have, for some finite constant~$g_2$,
\begin{equation}
\label{triv-kappa}
\bigl|V(\Gamma)\bigr|
\le g_2\varkappa\sum_{S\in\frakS}\bigl|\ssP(S)\bigr|.
\end{equation}
\end{lemma}

\begin{proofsect}{Proof}
Immediate from the definition of~$s$-skeletons.
\end{proofsect}

Lemma~\ref{L:triv} will be useful because of the following
observation:
Let~$\frakS$ be a collection of~$s$-skeletons and recall
that the minimal value of the surface tension, 
$\tau_\mini=\inf_{\bn\in\CalS_1}\tau_\beta(\bn)$ is strictly positive,~$\tau_\mini>0$.
Then
\begin{equation}
\label{triv-W}
\sum_{S\in\frakS}\mathscr{W}_\beta\bigl(\ssP(S)\bigr)
\ge 
\tau_\mini\sum_{S\in\frakS}\bigl|\ssP(S)\bigr|.
\end{equation}
Thus the bounds in~\twoeqref{triv-s}{triv-kappa} will allow us to
convert a lower bound on the overall contour surface area/volume into
a lower bound on the Wulff functional of the associated skeletons.

\smallskip
A little less trivial is the estimate on the difference 
between the volumes of 
$V(\Gamma)$ and~$V(\frakS)$:

\begin{lemma}
\label{L:skeleti}
There is a finite geometric constant~$g_3$ such that if~$\Gamma$ is a collection of contours and~$\frakS$ is a collection of 
$s$-skeletons with~$\Gamma\sim\frakS$, then
\begin{equation}
\label{g3}
\Bigl|\bigl|V(\Gamma)\bigr|-\bigl|V(\frakS)\bigr|\Bigr|\le
\bigl|V(\Gamma)\triangle V(\frakS)\bigr|
\le g_3 s\sum_{S\in\frakS}\bigl|\ssP(S)\bigr|.
\end{equation}
Here~$V(\Gamma)\triangle V(\frakS)$ denotes the symmetric difference 
of~$V(\Gamma)$ and~$V(\frakS)$.
\end{lemma}

\begin{proofsect}{Proof}
Follows by the same arguments as used in the proof of Theorem~5.13 in~\cite{DKS}.
%M COMMENTING THE PROOF
\begin{comment}
We will just rephrase the proof of Theorem~5.13 in~\cite{DKS}.
Let~$\Gamma= (\gamma_1,\dots,\gamma_m)$ and fix~$\gamma_k\in\Gamma$.
Let~$S_k\in\frakS$ with~$S_k=(x_0,\dots,x_n)$ 
be the skeleton compatible with~$\gamma_k$ 
and let~$S_{k,j}$ denote the segment of the straight line between~$x_j$ and 
$x_{j+1}$. Since~$\gamma$ passes through the skeleton points \emph{in the given order}, 
for each~$j$ there is a corresponding piece,~$\gamma_{k,j}$, of~$\gamma$ 
which  connects~$x_j$ and~$x_{j+1}$.

Let~$U_{k,j}$ be the subset of~$\R^2$ enclosed 
``between'' 
$S_{k,j}$ and~$\gamma_{k,j}$
(i.e.,~$U_{k,j}$ is the union of all  bounded connected components of 
$\R^2\setminus(S_{k,j}\cup\gamma_{k,j})$). 
We claim that
\begin{equation}
\label{UUgetin}
V(\Gamma)\triangle V(\frakS)\subseteq\bigcup_{k,j}U_{k,j}.
\end{equation}
Indeed, let~$x\in V(\Gamma)\triangle V(\frakS)$ and let~$\LL$ 
be a path connecting~$x$ to infinity which is 
generic with respect to both~$\frakS$ and~$\Gamma$.
Then~$\LL$ has the same parity of the number of 
intersections with~$\gamma_k$ and~$S_k$, \emph{unless}~$x\in U_{k,j}$ for 
some~$k$ and~$j$. By inspecting the definitions of~$V(\Gamma)$ and~$V(\frakS)$, 
\eqref{UUgetin} is proved.

Let~$U_s(\ssP(S_k))$ be the~$s$-neighborhood of the 
polygonal curve~$\ssP(S_k)$.
Since~$\partial U_{k,j}\subset U_s(\ssP(S_k))$, by 
\eqref{UUgetin} we have that
$V(\Gamma)\triangle V(\frakS)\subseteq \bigcup_k U_s(\ssP(S_k))$. From here 
\eqref{g3} directly follows.
\end{comment}
\end{proofsect}

\subsubsection{Probabilistic estimates}
The main reason why skeletons are useful is the availability of 
the so called 
\emph{skeleton upper bound}, originally due to Pfister~\cite{Pfister}. 
Recall that, for each~$A\subset\Z^2$, we use~$P_A^{+,\beta}$ to denote the 
probability distribution on spins in~$A$ with plus boundary condition on the 
boundary of~$A$. Given a set of skeletons, we let~$P_A^{+,\beta}(\frakS)=P_A^{+,\beta}(\{\sigma\colon \sigma\sim\frakS\})$ be the probability that~$\frakS$ is a skeleton of \emph{some} configuration in~$A$. Then we have:

\begin{lemma}[Skeleton upper bound]
\label{lemma-skeletonUB} 
For all~$\beta>\betac$, all finite~$A\subset\Z^2$, all scales~$s$ and all
collections~$\frakS$ of~$s$-skeletons in~$A$, we have
\begin{equation}
\label{skeleton-UB}
P_A^{+,\beta}(\frakS)\le \exp\bigl\{-\mathscr{W}_\beta(\frakS)\bigr\},
\end{equation}
where
\begin{equation}
%\label{}
\mathscr{W}_\beta(\frakS)=\sum_{S\in\frakS}\mathscr{W}_\beta\bigl(\ssP(S)\bigr).
\end{equation}
\end{lemma}

\begin{proofsect}{Proof}
This is exactly Eq.~(1.3.1) in~\cite{Bob+Tim}. The proof goes back
to~\cite{Pfister}, Lemma~6.7. For our purposes, the key ``splitting'' argument
is provided in Lemma~5.4 of~\cite{Pf-Velenik}. A special 
case of the key estimate appears in 
Eq.~(5.51) from Lemma~5.5 of~\cite{Pf-Velenik} with the correct
interpretation of the left-hand side.
\end{proofsect}

The bound~\eqref{skeleton-UB} will be used in several ways: First,
to show that the~$K\log L$-large contours in a box of side-length
$L$ are improbable, provided~$K$ is large enough; this is a consequence
of Lemma~\ref{lemma-Peierls} below. The absence of such contours will be
wielded to rule out the likelihood of other improbable scenarios. Finally,
after all atypical situations have been dispensed with, the skeleton
upper bound will deliver the contribution corresponding to
the term~$\sqrt\lambda$ in \eqref{LDP}.

\smallskip
An important consequence of the skeleton upper bound is the
following generalization of the Peierls estimate, which will be
useful at several steps of the proof of our main theorems.

\begin{lemma}
\label{lemma-Peierls} 
Let~$s=K\log L$ and let~$\mathscr{S}_{L,K}$ denote the set of all
$s$-skeletons that arise from contours in~$\Lambda_L$.
For each~$\beta>\betac$ and~$\alpha>0$,
there is a~$K_0=K_0(\alpha,\beta)<\infty$, 
such that 
\begin{equation}
\sum_{\frakS\subset\mathscr{S}_{L,K}}
\exp\bigl\{-\alpha\mathscr{W}_\beta(\frakS)\bigr\}\le1
\end{equation}
for (all~$L$ and) all~$K\ge K_0$.
\end{lemma}

\begin{proofsect}{Proof}
Let~$\mathscr{S}_{L,K}^0$ be the set of all~$K\log L$-skeletons~$S$
such that~$S=(x_1,\dots,x_k)$ with~$x_1=0$. 
By translation invariance,
\begin{equation}
\label{bd1} 
\sum_{\frakS\subset\mathscr{S}_{L,K}}
e^{-\alpha\mathscr{W}_\beta(\frakS)}\le
\sum_{n\ge1}\Bigl(L^2\!\!\!\sum_{S\in\mathscr{S}_{L,K}^0}\!\!\!
e^{-\alpha\mathscr{W}_\beta(\ssP(S))}\Bigr)^n,
\end{equation}
where the
prefactor~$L^2$ accounts for the translation entropy of each
skeleton within~$\Lambda_L$. The latter sum can be estimated by
mimicking the proof of Peierls' bound, where contour entropy was
bounded by that of the simple random walk on~$\Z^2$. Indeed,
each skeleton can be thought of as a sequence of steps with step-length
entropy at most~$32s^2$, where~$s=K\log L$,  and 
with each step weighted by a factor not exceeding~$e^{-\tau_\mini s}$. 
This and \eqref{triv-W} yield
\begin{equation}
\label{bd2}
\sum_{S\in\mathscr{S}_{L,K}^0}\!\!\!e^{-\alpha\mathscr{W}_\beta(\ssP(S))}
\le\sum_{m\ge 1}\bigl(32s^2
e^{-\alpha\tau_\mini s}\bigr)^{m}.
\end{equation}
By choosing~$K_0$ sufficiently large, the
right-hand side is less than~$\frac12 L^{-2}$ for all~$K\ge K_0$. 
Using this in \eqref{bd1}, the claim
follows. 
\end{proofsect}

Lemmas~\ref{lemma-skeletonUB} and~\ref{lemma-Peierls} will be used in the form of 
the following corollary:

\begin{corollary}
\label{cor-help}
Let~$\beta>\betac$,~$L\ge1$ and~$\kappa>0$ be fixed, and let~$\AA$ be the 
set of of configurations~$\sigma$ such that
$\mathscr{W}_\beta(\frakS)\ge\kappa$ for 
at least one collection of~$s$-skeletons~$\frakS$ 
satisfying~$\frakS\sim\sigma$. Let~$\alpha\in(0,1)$, and let~$K_0(\alpha,\beta)$ be as in 
Lemma~\ref{lemma-Peierls}. If~$s=K\log L$ with~$K\ge K_0(\alpha,\beta)$, then
\begin{equation}
P_L^{+,\beta}(\AA)\le e^{-(1-\alpha)\kappa}.
\end{equation}
\end{corollary}

\begin{proofsect}{Proof}
By the assumptions of the Lemma, we have
\begin{equation}
%\label{}
P_L^{+,\beta}(\AA)\le\sum_{\begin{subarray}{c}
\frakS\subset\mathscr{S}_{K,L}\\
\mathscr{W}_\beta(\frakS)\ge\kappa
\end{subarray}}P_L^{+,\beta}(\frakS),
\end{equation}
where we used the notation~$P_L^{+,\beta}(\frakS)
=P_L^{+,\beta}(\{\sigma\colon \sigma\sim\frakS\})$.
Lemma~\ref{lemma-skeletonUB} then implies
\begin{equation}
%\label{}
\qquad
P_L^{+,\beta}(\AA)\le\sum_{\begin{subarray}{c}
\frakS\subset\mathscr{S}_{K,L}\\
\mathscr{W}_\beta(\frakS)\ge\kappa
\end{subarray}}
e^{-\mathscr{W}_\beta(\frakS)}
\le e^{-(1-\alpha)\kappa}\sum_{\frakS\subset\mathscr{S}_{K,L}}
e^{-\alpha\mathscr{W}_\beta(\frakS)}.
\qquad
\end{equation}
Here we wrote~$e^{-\mathscr{W}_\beta(\frakS)}=e^{-\alpha\mathscr{W}_\beta(\frakS)}
e^{-(1-\alpha)\mathscr{W}_\beta(\frakS)}$ and then invoked 
to bound~$\mathscr{W}_\beta(\frakS)\ge
\kappa$ to estimate
$e^{-(1-\alpha)\mathscr{W}_\beta(\frakS)}$ by~$e^{-(1-\alpha)\kappa}$. 
Finally, we dropped the constraint
to~$\mathscr{W}_\beta(\frakS)\ge\kappa$ in the last sum.
Since~$s=K\log L$ with~$K\ge K_0(\alpha,\beta)$, the last sum is less than one 
by Lemma~\ref{lemma-Peierls}.
\end{proofsect}

Ideas similar to those used in the proof of Lemma~\ref{lemma-Peierls} can be
used  to estimate the probability of the occurrence of an~$s$-large contour:

\begin{lemma}
\label{lemma-nolog} For each~$\beta>\betac$, there exist  
a constant~$\alpha(\beta)>0$ such that 
\begin{equation}
\label{blabla}
P_A^{+,\beta}\bigl(\Gamma_s(\sigma)\ne\emptyset\bigr)\le |A|e^{-\alpha(\beta)s}
\end{equation}
for any finite~$A\subset \Z^2$ and any scale~$s$.
\end{lemma}

\begin{proofsect}{Proof}
Fix $\alpha>0$ and suppose without loss of generality that~$|A|>1$ and~$s\ge\alpha^{-1}\log|A|$
for some~$\alpha>0$. If~$\Gamma_s(\sigma)\ne\emptyset$, the associated 
$s$-skeleton must satisfy~$\mathscr{W}_\beta(\frakS)\ge\tau_\mini s$. 
Invoking \eqref{skeleton-UB} a variant of the estimate \twoeqref{bd1}{bd2} 
(here is where~$s\ge\alpha^{-1}\log|A|$ enters into the play), we show that
$P_A^{+,\beta}(\Gamma_s(\sigma)\ne\emptyset)\le C|A|s^2e^{-\frac12\tau_\mini s}$, 
where $C>0$ is a constant.
From here the bound \eqref{blabla} follows by absorbing the factor~$Cs^2$ into the exponential.
\end{proofsect}

\subsubsection{Quantitative estimates around Wulff minimum}
The existence of a minimum for the functional \eqref{Wbeta} and a
coarse-graining scheme supplemented with a bound of the type in 
\eqref{skeleton-UB} tell us the following: Consider a 
collection~$\Gamma$ of contours, all of which
are roughly of the same scale and which enclose a fixed total 
volume, and
suppose that the value of the Wulff functional on a~$\frakS$ with~$\frakS\sim\Gamma$ is 
close to the Wulff minimum. Then~(1) it must be the case that~$\Gamma$ consists of a single contour and~(2) the shape of this 
contour must be close to the Wulff shape. A quantitative (and mathematically
precise) version of this statement is given in the 
forthcoming lemma:

\begin{lemma}
\label{lemma-onecontour} 
For any~$\beta\ge\betac$, there exist
constants~$\epsilon_0=\epsilon_0(\beta)\in(0,1)$,~$c=c(\beta)>0$, and
$C=C(\beta)<\infty$ such that the following holds for all
$\epsilon\in(0,\epsilon_0)$: Let~$\Gamma$ be a collection of contours such that
$\diam\gamma > c \epsilon\sqrt{|V(\Gamma)|}$ for all
$\gamma\in\Gamma$ and let~$s$ be a scale function satisfying 
$s\le\epsilon\sqrt{|V(\Gamma)|}$.
Let~$\frakS$ be a collection of~$s$-skeletons
compatible with~$\Gamma$,~$\frakS\sim\Gamma$, such that
\begin{equation}
\label{2.16}
\mathscr{W}_\beta(\frakS)\le w_1\sqrt{|V(\Gamma)|}(1+\epsilon).
\end{equation}
Then~$\Gamma$ consists of a single contour,~$\Gamma=\{\gamma\}$,
and there is an~$x\in\R^2$ such that
\begin{equation}
\label{4.24} \dH\bigl(V(\gamma),\sqrt{|V(\gamma)|}W+x\bigr) \le
c\sqrt{\epsilon}\sqrt{|V(\gamma)|},
\end{equation}
where~$W$ is the Wulff shape of unit area centered at the origin. Moreover,
\begin{equation}
\label{4.24a}
\bigl| | V(\gamma)|-|V(\frakS)|\bigr|\le C 
\epsilon|V(\gamma)|.
\end{equation}
\end{lemma}

\begin{proofsect}{Proof}
We begin by noting that, by the assumptions of the present Lemma,~$|V(\Gamma)|$ and~$|V(\frakS)|$
have to be of the same order of magnitude. 
More precisely, we claim that
\begin{equation}
\label{G-S}
\bigl| |V(\Gamma)|-|V(\frakS)|\bigr|\le C\epsilon\bigl|V(\Gamma)\bigr|
\end{equation}
holds with some~$C=C(\beta)<\infty$ 
independent of~$\Gamma$,~$\frakS$ and~$\epsilon$. Indeed, 
from \eqref{triv-W} and \eqref{2.16} we~have
\begin{equation}
\label{sumaP(S)}
\sum_{S\in\frakS}\bigl|\ssP(S)\bigr|\le
\tau_\mini^{-1}\mathscr{W}_\beta(\frakS)\le
w_1(1+\epsilon)\tau_\mini^{-1}\sqrt{|V(\Gamma)|},
\end{equation}
which, using Lemma~\ref{L:skeleti} and the bounds
$s\le\epsilon\sqrt{|V(\Gamma)|}$ and~$\epsilon\le1$, gives \eqref{G-S} 
with~$C=2g_3 w_1\tau_\mini^{-1}$.

The bound \eqref{G-S} essentially allows us to replace~$V(\Gamma)$ 
by~$V(\frakS)$ in \eqref{2.16}.
Applying Theorem~2.10 from~\cite{DKS} to the set of skeletons
$\frakS$ rescaled by~$|V(\frakS)|^{-1/2}$, we can conclude that there is 
point~$x\in\R^2$ and a skeleton
$S_0\in\frakS$ such that
\begin{equation}
\label{4.25} 
\dH\bigl(\ssP(S_0),\sqrt{|V(\frakS)|}\partial W+x\bigr) \le 
\alpha\sqrt{\epsilon}\sqrt{|V(\frakS)|},
\end{equation}
and
\begin{equation}
\label{4.26} \sum_{S\in\frakS
\setminus\{S_0\}}\bigl|\ssP(S)\bigr|\le
\alpha\epsilon\sqrt{|V(\frakS)|},
\end{equation}
where~$\alpha$ is a constant proportional to 
the ratio of the maximum and the minimum of the surface tension.
Using \eqref{G-S} once more, we can modify \twoeqref{4.25}{4.26} by replacing 
$V(\frakS)$ on the right-hand sides
by~$V(\Gamma)$ at the cost of changing~$\alpha$ to~$\alpha(1+C)$. Moreover, since 
\eqref{G-S} also implies that
$|\sqrt{|V(\Gamma)|}-\sqrt{|V(\frakS)|}|\le C
\epsilon \sqrt{|V(\Gamma)|}$,
we have 
\begin{equation}
\label{dH(G-S)}
\dH\bigl(\sqrt{|V(\Gamma)|}\partial W,
\sqrt{|V(\frakS)|}\partial W\bigr)\le C\epsilon\diam W
\sqrt{|V(\Gamma)|}.
\end{equation}
Let~$\gamma\in\Gamma$ be the contour corresponding to~$S_0$. By
the definition of skeletons,~$\dH(\gamma,\ssP(S_0))\le
s\le \epsilon\sqrt{|V(\Gamma)|}$.
Combining this with \eqref{dH(G-S)}, the modified bound \eqref{4.25}, and 
$\epsilon\le1$, we get
\begin{equation}
\label{4.25a} 
\dH\bigl(\gamma,\sqrt{|V(\Gamma)|}\partial W+x\bigr) \le
c\sqrt{\epsilon}\sqrt{|V(\Gamma)|}
\end{equation}
for any~$c\ge 1+\alpha(1+C)+C\diam W$.
(From the properties of~$W$, it is easily shown that~$\diam W$ is of the order of unity.)

Let us proceed by proving that~$\Gamma=\{\gamma\}$.
For any~$\gamma'\in\Gamma\setminus\{\gamma\}$, let~$S_{\gamma'}$ 
be the unique skeleton in~$\frakS$ such that~$\gamma'\sim S_{\gamma'}$.
Since~$\diam\gamma' \le |\ssP(S_{\gamma'})|+s$ and, since also 
$|\ssP(S_{\gamma'})|\ge s$, we have~$\diam\gamma'\le
2|\ssP(S_{\gamma'})|$.
Using the modified bound \eqref{4.26}, we get
\begin{equation}
\label{2.24}
\diam\gamma' \le 2\bigl|\ssP(S_{\gamma'})\bigr| \le 2\alpha (1+C) 
\epsilon\sqrt{|V(\Gamma)|}.
\end{equation}
If~$c$ also satisfies the inequality~$c>2\alpha(1+C)$, then this 
estimate 
contradicts the assumption that
$\diam\gamma'\ge c \epsilon\sqrt{|V(\Gamma)|}$ for all~$\gamma'\in\Gamma$. Hence, 
$\Gamma=\{\gamma\}$ as claimed.

Thus,~$V(\Gamma)=V(\gamma)$ and the bound \eqref{4.24a} is directly implied by 
\eqref{G-S}.
Moreover, \eqref{4.25a} holds with~$V(\Gamma)$ replaced by
$V(\gamma)$ on both sides. To prove \eqref{4.24}, it remains to show that the 
naked~$\gamma$ on the left-hand side of \eqref{4.25a} can be replaced by 
$V(\gamma)$. But that is trivial because~$\gamma$ is the boundary of~$V(\gamma)$ 
and the Hausdorff distance of two closed sets in~$\R^2$ equals the Hausdorff 
distance of their boundaries.
\end{proofsect}

\subsection{Small-contour ensemble}
\label{sec2.3}\noindent
The goal of this section is to collect some estimates for
the probability in~$P_L^{+,\beta}$ conditioned on the fact that all
contours are~$s$-small in the sense that~$\Gamma_s(\sigma)=\emptyset$.
Most of what is to follow appears, in various guises, in the existing literature (cf Remark~\ref{rem6}).
For some of the estimates (Lemmas~\ref{L:M-fluctuations} and~\ref{lemma-Gauss-positive}) we will 
actually provide a proof, while for others (Lemma~\ref{lemma-Gauss}) we can quote directly.

\subsubsection{Estimates using the GHS inequality}
The principal resource for what follows are two basic properties of the correlation function 
of Ising spins. Specifically, let~$\langle \sigma_x;\sigma_y\rangle^{+,\beta}_{A,\bh}$ 
denote the truncated correlation function of the Ising model in a set~$A\subset\Z^2$ 
with plus boundary condition, in non-negative inhomogeneous external fields 
$\bh=(h_x)$ and inverse temperature~$\beta$. Then:
\begin{enumerate}
\item[(1)]
If~$\beta>\betac$, then the correlations in infinite volume decay exponentially, i.e., 
we have 
\begin{equation}
\label{2.31}
\langle \sigma_x;\sigma_y\rangle^{+,\beta}_{\Z^2,\bh}\le e^{-\Vert x-y\Vert/\xi}
\end{equation}
for some~$\xi=\xi(\beta)<\infty$ and all~$x$ and~$y$.
\item[(2)]
The GHS inequality implies that the finite-volume correlation function,~$\langle 
\sigma_x;\sigma_y\rangle^{+,\beta}_{A,\bh}$, is dominated by the infinite-volume 
correlation function at any pointwise-smaller 
field:
\begin{equation}
\label{2.32}
0\le\langle 
\sigma_x;\sigma_y\rangle^{+,\beta}_{A,\bh}\le
\langle 
\sigma_x;\sigma_y\rangle^{+,\beta}_{\Z^2,\bh'}
\end{equation}
for all 
$A\subset\Z^2$ and all~$\bh'=(h'_x)$ with~$h_x'\in[0,h_x]$ for all 
$x$.
\end{enumerate}
Note that, via \eqref{2.32}, the exponential decay \eqref{2.31} holds uniformly in~$A\subset\Z^2$.
Part~(1) is a consequence of the main result of~\cite{CCS}, see~\cite{Schonmann-Shlosman}; the
GHS inequality from part~(2) dates back~to~\cite{GHS}.

\smallskip
Now we are ready to state the desired estimates.
Let~$A\subset\Z^2$ be a finite set and let~$s$ be a scale
function. Let~$P_A^{+,\beta,s}$ be the Gibbs measure of the Ising
model in~$A\subset\Z^2$ conditioned on the event~$\{\Gamma_s(\sigma)=\emptyset\}$ 
and let us use~$\langle-\rangle_A^{+,\beta,s}$ to denote the expectation with 
respect to
$P_A^{+,\beta,s}$. Then we have the following bounds:

\begin{lemma}
\label{L:M-fluctuations}
For each~$\beta>\betac$, 
there exist constants~$\alpha_1(\beta)$ and~$\alpha_2(\beta)$ 
such that 
\begin{equation}
\label{uno}
\bigl|\langle M_A\rangle_A^{+,\beta,s}-\mstar|A|\bigr|
\le \alpha_1(\beta)\bigl(|\partial A|+|A|^2e^{-\alpha_2(\beta)\,s}\bigr)
\end{equation}
for each finite set~$A\subset\Z^2$ and any scaling function~$s$.
Moreover, if~$A'\subset A$, then
\begin{equation}
\label{due}
\bigl|\langle M_A\rangle_A^{+,\beta,s}-\langle M_{A\smallsetminus 
A'}\rangle_{A\smallsetminus A'}^{+,\beta,s}\bigr|
\le \alpha_1(\beta)\bigl(|A'|+|A|^2e^{-\alpha_2(\beta)\,s}\bigr).
\end{equation}
\end{lemma}

\begin{proofsect}{Proof}
By Lemma~\ref{lemma-nolog}, we have 
$P_A^{+,\beta}(\Gamma_s(\sigma)\ne\emptyset)\le  
|A| e^{-\alpha_2 s}$ for some 
$\alpha_2>0$, independent of~$A$. 
Note that we can suppose that $|A| e^{-\alpha_2 s}$ does not exceed, e.g., $1/2$, because otherwise \twoeqref{uno}{due} can be ensured by deterministic estimates.
An easy bound then shows that, for some $\alpha_1'=\alpha_1'(\beta)<\infty$,
\begin{equation}
\label{sbd}
\bigl|\langle M_A\rangle^{+,\beta,s}_A-\langle M_A\rangle^{+,\beta}_A\bigr|\le
\alpha_1'|A|^2e^{-\alpha_2 s}.
\end{equation}
Therefore, it suffices to prove the bounds \twoeqref{uno}{due} without the restriction 
to the ensemble of~$s$-small contours.
%M JUST STATING THE INEQUALITY
The proof will use that, for any~$B\subset\Z^2$ we have
\begin{equation}
\label{f1}
0\le
\langle \sigma_x\rangle^{+,\beta}_B-\langle \sigma_x\rangle^{+,\beta}_{B\cup\{y\}}\le 
e^{-\Vert x-y\Vert/\xi}.
\end{equation}
This inequality is a direct consequence of properties~(1-2) above.
The original derivation goes back to~\cite{Bricmont-Lebowitz-Pfister}.

%%% PREVIOUS VERSION
\begin{comment}
Next  we claim that, for any~$B\subset\Z^2$ we have
\begin{equation}
\label{f1}
0\le
\langle \sigma_x\rangle^{+,\beta}_B-\langle \sigma_x\rangle^{+,\beta}_{B\cup\{y\}}\le 
e^{-\Vert x-y\Vert/\xi}.
\end{equation}
Indeed, the difference of the two expectations can 
be written as an integral 
$\int_0^\infty \textd h \langle 
\sigma_x;\sigma_y\rangle^{+,\beta}_{B\cup\{y\},\bh}$,
where~$\bh=(h_z)$ 
is such that~$h_z=h\delta_{y,z}$. 
By property~(2) of the truncated correlation 
function, we have that 
$\langle \sigma_x;\sigma_y\rangle^{+,\beta}_{B\cup\{y\},
\bh}$ is non-negative, which proves the left inequality in \eqref{f1}, while it is 
bounded above by the same correlation function with~$B=\Z^2$. The integral 
representation can be used again for the correlation function in the infinite 
volume with the result
\begin{equation}
%\label{} 
\langle \sigma_x\rangle^{+,\beta}_B-\langle\sigma_x\rangle^{+,\beta}_{B\cup\{y\}}
\le 
\frac{\bigl\langle \sigma_x\tfrac{1+\sigma_y}2\bigr\rangle^{+,\beta}}
{\bigl\langle \tfrac{1+\sigma_y}2\bigr\rangle^{+,\beta}}
-\langle \sigma_x\rangle^{+,\beta}
\le
\frac{\bigl\langle \sigma_x;\tfrac{1+\sigma_y}2\bigr\rangle^{+,\beta}}
{{\bigl\langle \tfrac{1+\sigma_y}2\bigr\rangle^{+,\beta}}}\le
\langle 
\sigma_x;\sigma_y\rangle^{+,\beta},
\end{equation}
where we used that 
$\langle1+\sigma_y\rangle^{+,\beta}\ge1$ to derive the last inequality.
Using the 
property~(1) above,  the right-hand side is bounded by~$e^{-\Vert x-y\Vert/\xi}$.
\end{comment}
%%%%

The bound \eqref{f1} immediately implies both \eqref{uno} and \eqref{due}. Indeed, 
using \eqref{f1} for all~$x\in A$ and~$y\in B\setminus A$, we have for all 
$A\subseteq B\subseteq\Z^2$ that
\begin{equation}
%\label{}
0\le\langle 
M_A\rangle^{+,\beta}_A-\langle M_A\rangle^{+,\beta}_B\le 
\sum_{x\in A}\sum_{y\in B\setminus A}e^{- \Vert x-y\Vert/\xi} \le\alpha_1''|\partial A|,
\end{equation}
where~$\alpha_1''=\alpha_1''(\beta)<\infty$. This and \eqref{sbd} 
directly imply \eqref{uno}. To get \eqref{due}, we also need to note that~$|M_A-
M_{A\setminus A'}|\le|A'|$.
\end{proofsect}

Our next claim concerns an upper bound on the probability that the magnetization 
in the plus state deviates from its mean by a positive amount:

\begin{lemma}
\label{lemma-Gauss-positive} Let~$\beta>\betac$ and let~$\chi=\chi(\beta)$ be the 
susceptibility. Then there exists
a constant~$K=K(\beta)$ such that 
\begin{equation}
\label{ubp}
P_A^{+,\beta,s}\bigl(M_A\ge \langle M_A\rangle_A^{+,\beta}+\mstar\,
v\bigr)\le 2
e^{-\frac{(v \mstar)^2}{2\chi|A|}}
\end{equation}
for any finite
$A\subset\Z^2$, any~$v\ge0$, and any~$s\ge K\log|A|$.
\end{lemma}

\begin{proofsect}{Proof}
Let~$\MM$ denote the event $\MM=\{\sigma\colon M_A\ge \langle
M_A\rangle_A^{+,\beta}+\mstar\, v\}$. By Lemma~\ref{lemma-nolog} we
have that~$P_A^{+,\beta,s}(\MM)\le 2P_A^{+,\beta}(\MM)$, so we just
need to estimate~$P_A^{+,\beta}(\MM)$. Consider the cumulant generating function
$F_A^{+,\beta}(h)=\log \langle e^{hM_A}\rangle_A^{+,\beta}$. The
exponential Chebyshev inequality then gives
\begin{equation}
%\label{}
\log P_A^{+,\beta}(\MM)\le F_A^{+,\beta}(h)-h\langle
M_A\rangle_A^{+,\beta}-h\mstar\, v,\qquad h\ge0.
\end{equation}
By the property~(2) of the truncated correlation function, we get
\begin{equation}
%\label{}
\frac{\textd^2 F_A^{+,\beta}}{\textd h^2}(h)=\langle
M_A;M_A\rangle_{A,\bh}^{+,\beta}
\le \langle
M_A;M_A\rangle_{A,\bzero}^{+,\beta},
\end{equation}
where~$\bh=(h_x)$ with~$h_x=h$ for all~$x\in\Z^2$ and where~$\bzero$ is the zero field.
Since~$F_A^{+,\beta}(0)=0$ and~$\frac{\textd}{\textd h}F_A^{+,\beta}(0)=\langle
M_A\rangle_A^{+,\beta}$, we get the bound
\begin{equation}
%\label{}
F_A^{+,\beta}(h)\le h\langle M_A\rangle_A^{+,\beta}+
\frac{h^2}2\langle M_A;M_A\rangle_{A,\bzero}^{+,\beta}.
\end{equation}
Now, once more by the property~(2) above,
\begin{equation}
%\label{}
|A|^{-1}\langle M_A;M_A\rangle_{A,\bzero}^{+,\beta}\le
|A|^{-1}\langle 
M_A;M_A\rangle_{\Z^2,\bzero}^{+,\beta}
\le |A|^{-1} \sum_{x\in A}\sum_{y\in\Z^2} \langle \sigma_x;\sigma_y\rangle^{+,\beta}=\chi,
\end{equation}
where the sums converge by the property~(1).
The claim follows by optimizing over~$h$. 
\end{proofsect}

\begin{remark}
\label{rem6}
The bound in Lemma~\ref{lemma-Gauss-positive} 
corresponds to
Eq.~(9.33) of Proposition 9.1
in~\cite{Pf-Velenik} proved with the help of Lemma 5.1 from
\cite{Pfister}. Similarly, the estimates in Lemma~\ref{L:M-fluctuations}
are closely related to the bounds in Lemma~2.2.1 of~\cite{Bob+Tim}. We included 
the proofs of both statements to pinpoint the exact formulation needed for our 
analysis as well as to reduce the number of extraneous 
references.
\end{remark}

\subsubsection{Gaussian control of negative deviations}
Our last claim concerns the deviations of the plus magnetization in 
the \emph{negative} direction. Unlike in the previous Section,
here the restriction to the small contour is
crucial because, obviously, if the deviation 
is too large, there is a possibility
of forming a droplet which cannot be 
controlled by bulk estimates.

\smallskip
Let~$\beta>\betac$ and let~$v$ be such
that~$\langle M_A\rangle_A^{+,\beta,s}-2\mstar\, v$ is an allowed
value of~$M_A$. Define 
$\Omega_{A}^s(v)$ by the expression
\begin{equation}
\label{Gauss-asymp} P_A^{+,\beta,s}\bigl(M_A=\langle
M_A\rangle_A^{+,\beta,s}-2v\mstar\,\bigr)=\frac1{\sqrt{2\pi\chi|A|}}
\exp\Bigl\{-2\frac{(\mstar)^2}{\chi|A|}\,v^2+
\Omega_{A}^s(v)\Bigr\}.
\end{equation}
Then we have:

\begin{lemma}[Gaussian estimate]
\label{lemma-Gauss} For each~$\beta>\betac$ and
each set of positive constants
$a_1,a_2,a_3$, there are constants~$C<\infty$ and 
$K<\infty$ such that if~$s=K\log L$, then
\begin{equation}
\label{2.45}
\bigl|\Omega_{A}^s(v)\bigr|\le C\max\Bigr\{K\frac{v^2}{L^3}\log L,\, 
\frac{v^3}{L^4}\Bigl\}
\end{equation}
for all allowed values of~$v$ such that
\begin{equation}
\label{frst} 0\le v\le a_1\frac{L^2}{\log L}
\end{equation}
and all connected sets~$A\subset\Z^2$ such that
\begin{equation}
\label{scnd} a_2L^2\le|A|\le L^2\quad\text{and}\quad |\partial
A|\le a_3 L
\log L.
\end{equation}
\end{lemma}

\begin{proofsect}{Proof} This is a reformulation of
(a somewhat nontrivial)
Lemma~2.3.3 from~\cite{Bob+Tim}. \end{proofsect}

\section{Lower bound}
\label{sec3}
\smallskip\noindent
In this Section we establish a lower bound for the asymptotic stated in \eqref{LDP}. 
In addition to its contribution to the proof of Theorem~\ref{LDP-thm}, this lower bound will
play an essential role in the proofs of Theorem~\ref{main-result} and 
Corollary~\ref{cor-main-result}.
A considerable part of the proof hinges on the Fortuin-Kasteleyn representation of the Ising
(and Potts) models, which makes the technical demands of this section rather different 
from those of the following sections.

\subsection{Large-deviation lower bound}
\label{sec3.1}\noindent
This section is devoted to the proof of the following theorem:

\begin{theorem}[Lower bound]
\label{lowerbound} 
Let~$\beta>\betac$ and let~$(v_L)$ be a
sequence of positive numbers such that~$\mstar\,\vert\Lambda_L\vert 
- 2\mstar\, v_L$ is an allowed value of~$M_L$ for all~$L$. 
Suppose that the limit \eqref{Delta-lim} exists with
$\Delta\in(0,\infty)$. Then there exists a sequence~$(\epsilon_L)$
with~$\epsilon_L\to 0$ such that
\begin{equation}
\label{LBeq} P_L^{+,\beta}
\bigl(M_L=\mstar\,\vert\Lambda_L\vert - 2\mstar\, v_L\bigr)
\ge \exp\bigl\{-w_1 \sqrt{v_L} \bigl(\,\inf_{0\le\lambda\le1}\Phi_\Delta(\lambda)
+\epsilon_L\bigr)\bigr\}
\end{equation}
holds for all~$L$.
\end{theorem}

\begin{remark}
It is worth noting that, unlike in the corresponding 
statements of the lower bounds in \cite{DKS,Bob+Tim}, we do not require any 
control over how fast the error~$\epsilon_L$ tends to zero as~$L\to\infty$. 
Indeed, it turns out that in the regime of finite~$\Delta$, the simple convergence~$\epsilon_L\to0$ 
will be enough to prove our main results. 
However, in the cases when~$v_L$ tends to infinity so fast that~$\Delta$ is infinite, a proof 
would probably need also \emph{some} information about the rate of the convergence~$\epsilon_L\to0$.
\end{remark}

The strategy of the proof will simply be to produce a near-Wulff
droplet that comprises a particular fraction of the volume~$v_L$. 
The droplet will account for its requisite share of the deficit magnetization 
and we then force the exterior to absorb
the rest. The probability of the latter event is estimated 
by using the truncated contour ensemble.

Let us first attend to the production of the droplet. Consider the
Wulff shape~$W$ of unit area centered at the origin and a closed, self-avoiding
polygonal curve~$\ssP\subset W$. We will assume that the vertices of~$\ssP$ have 
rational coordinates and, if~$N$ denotes the number of vertices of~$\ssP$, that 
each vertex is at most~$1/N$ away from the boundary of~$W$.
Let~$\Int\,\ssP$ denote the set of points~$x\in\R^2$ surrounded by~$\ssP$.
For any~$t,r>1$, let~$\ssP_0, \ssP_1, \ssP_2,
\ssP_3$ be four magnified copies of~$\ssP$ obtained by 
rescaling~$\ssP$ by factors~$t$,~$t+r$,
$t+2r$, and
$t+3r$, respectively. (Thus, for instance,~$\ssP_0=\{x\in\R^2\colon 
x/t\in\ssP\}$.) This yields three ``coronas''~$K^{\text{\rm I}}_{t,r}
=\Int\,\ssP_1\setminus\Int\,\ssP_0$,~$K^{\text{\rm II}}_{t,r}
=\Int\,\ssP_2\setminus\Int\,\ssP_1$, and~$K^{\text{\rm III}}_{t,r}
=\Int\,\ssP_3\setminus\Int\,\ssP_2$ surrounding~$\ssP_0$. 
Let~$\K_{t,r}^{\text{\rm I}}=K^{\text{\rm I}}_{t,r}\cap\Z^2$, and similarly for
$\K_{t,r}^{\text{\rm II}}$ and~$\K_{t,r}^{\text{\rm III}}$.

Recall that a~$*$-connected circuit in~$\Z^2$ is a closed
path on vertices of~$\Z^2$ whose elementary steps connect either nearest
or next-nearest neighbors. Let~$\EE_{t,r}$ be the set of
configurations~$\sigma$ such that~$\K_{t,r}^{\text{\rm
I}}$ contains a~$*$-connected circuit of sites~$x\in\Z^2$ with
$\sigma_x=-1$ and~$\K_{t,r}^{\text{\rm III}}$ contains a
$*$-connected circuit of sites~$x\in\Z^2$ with~$\sigma_x=+1$. The
essential part of our lower bound comes from the following
estimate:

\begin{lemma}
\label{circuit-lemma} 
Let~$\beta>\betac$ and let~$\ssP$ be a polygonal curve
as specified above. For any pair of sequences~$(t_L)$ and~$(r_L)$ tending 
to infinity as~$L\to\infty$ in such a way that
\begin{equation}
\label{t-r-cond} 
t_LL^{-1}\to0,\quad
t_L r_L e^{-r_L\tau_\mini/3}\to0 
\quad\text{\rm and}\quad r_Lt_L^{-1}\to0,
\end{equation}
there is a sequence~$(\epsilon_L')$ with~$\epsilon'_L\to0$
such that
\begin{equation}
P_L^{+,\beta}(\EE_{t_L,r_L})\ge\exp\bigl\{-
t_L\mathscr{W}_\beta(\ssP)(1+\epsilon_L')\bigr\},
\end{equation}
for all~$L\ge1$.
\end{lemma}

The proof of this lemma requires some substantial preparations and
is therefore deferred to Section~\ref{sec3.2}. Using
Lemma~\ref{circuit-lemma}, we can prove Theorem~\ref{lowerbound}.

\begin{figure}[t]
\refstepcounter{obrazek}
\vglue0.3cm
\ifpdf \centerline{\includegraphics[width=3.0truein]{polygons.pdf}}
\else
\centerline{\epsfxsize=3.0truein\epsffile{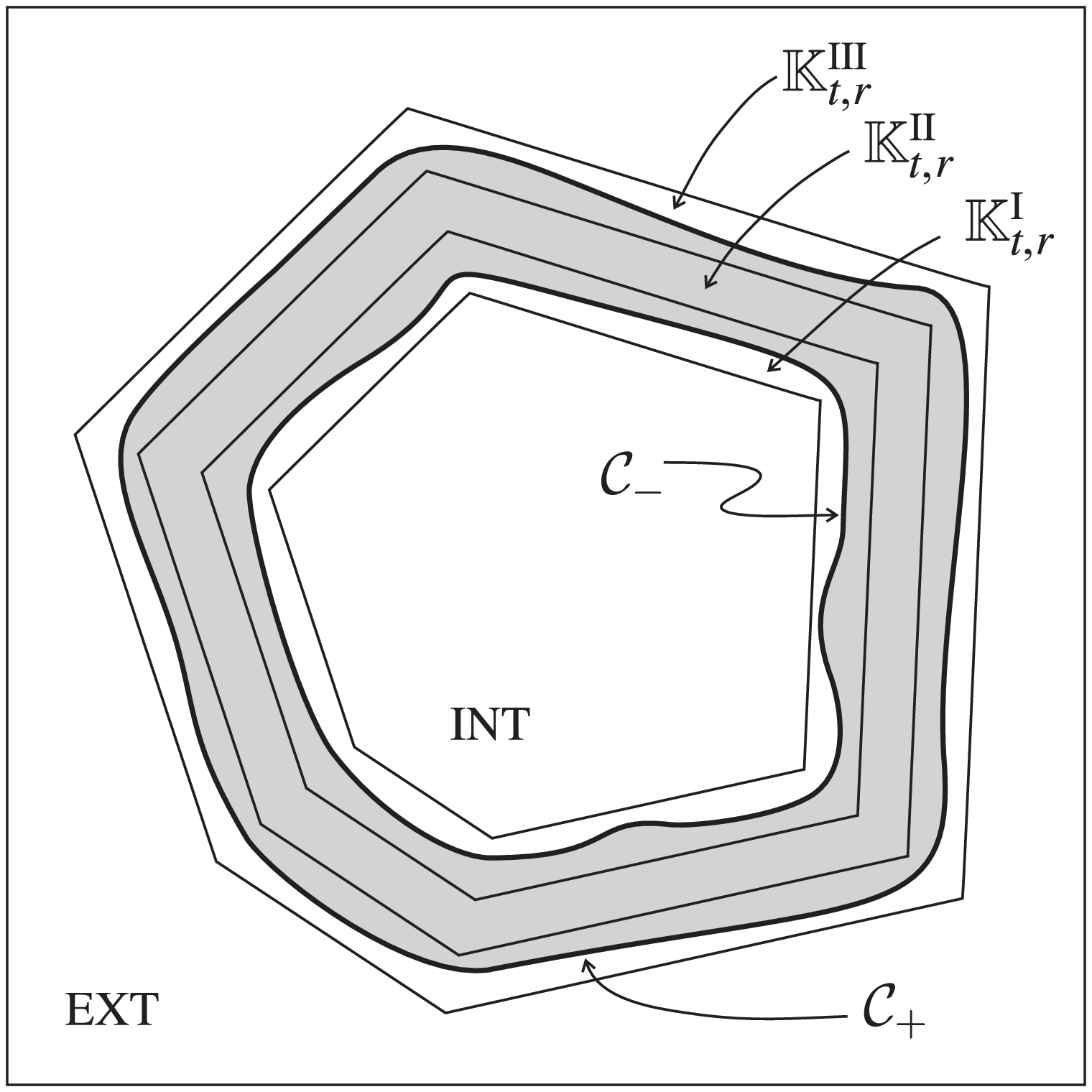}} 
\fi
\bigskip
\begin{quote}
%\small
\fontsize{9.5}{7}\selectfont
{\sc Figure~\theobrazek.\ }
\label{fig2}
An illustration of the 
``coronas''~$\mathbb K_{t,r}^{\text{\rm I}}$,~$\mathbb K_{t,r}^{\text{\rm II}}$,
$\mathbb K_{t,r}^{\text{\rm III}}$, 
the sets~$\IN$ and~$\OUT$, and the~$*$-connected circuits~$\mathcal C_+$ and~$\mathcal C_-$ 
of plus and minus sites, respectively, which are used in Lemma~\ref{circuit-lemma} 
and the proof of Theorem~\ref{lowerbound}. 
Going from inside out, the four polygons correspond to~$\ssP_0$,~$\ssP_1$,~$\ssP_2$ and~$\ssP_3$; 
the shaded region denotes the set~$A_\pm$.
\normalsize
\end{quote}
\end{figure}

\begin{proofsect}{Proof of Theorem~\ref{lowerbound}}
Let us introduce the abbreviation
\begin{equation}
%\label{}
\MM_L=\bigl\{\sigma\colon M_L=\mstar\,|\Lambda_L|-2\mstar\, v_L\bigr\}
\end{equation}
for the central event in question. 
Suppose first that~$\Delta\le\Deltac$, where~$\Deltac$ is as in~\eqref{Delta-c2}.
Proposition~\ref{prop-Phi} then guarantees that 
$\inf_{0\le\lambda\le1}\Phi_\Delta(\lambda)=\Phi_\Delta(0)=\Delta$. In particular, 
there is no need to produce a droplet in the system.
Let~$s=K\log L$. By restricting to the set of configurations 
$\{\sigma\colon\Gamma_s(\sigma)=\emptyset\}$ we get 
\begin{equation}
\label{MLlbd}
P_L^{+,\beta}(\MM_L)\ge P_L^{+,\beta,s}(\MM_L)P_L^{+,\beta}
\bigl(\Gamma_s(\sigma)=\emptyset\bigr).
\end{equation}
The resulting lower bound is then a consequence of \eqref{Gauss-asymp}, 
Lemma~\ref{lemma-Gauss} and Lemma~\ref{lemma-nolog}, provided~$K$ is sufficiently 
large. 

To handle the remaining cases,~$\Delta>\Deltac$, we \emph{will} have to produce 
a droplet.
Fix a polygon~$\ssP$ with the above properties, let~$\Vol(\ssP)$ denote the 
two-dimensional Lebesgue volume of its interior, and let~$|\ssP|$ denote the size 
(i.e., length) of its boundary. 
Let~$\lambda=\lambda_\Delta$, where~$\lambda_\Delta$ is as defined in 
\eqref{xDelta2}, and recall that, for this choice of~$\lambda$, we have 
$\Phi_\Delta(\lambda)=\inf_{0\le\lambda'\le1}\Phi_\Delta(\lambda')$ and 
$\lambda\ge\lamc>0$. 
Since the goal is to produce a droplet of volume~$\lambda v_L$, we let~$t_L=\sqrt{\lambda v_L}$ and pick 
$r_L$ be such that \eqref{t-r-cond} holds as~$L\to\infty$. 
Abbreviating~$\EE_L=\EE_{t_L,r_L}$, we let~$(\epsilon_L')$ denote 
the corresponding sequence from
Lemma~\ref{circuit-lemma}. (Note that~$\epsilon_L'$ may depend on~$\ssP$.)

For configurations in~$\EE_L$, 
let~$\CC_+$ be the innermost~$*$-connected circuit of 
plus spins in~$\K_{t,r}^{\text{\rm III}}$ and let~$\CC_-$ denote the
outermost~$*$-connected circuit of minus spins in
$\K_{t,r}^{\text{\rm I}}$. 
Let~$\IN$ be the set of sites
in the interior of~$\CC_-$ and let~$\OUT$ be the set of sites 
in~$\Lambda_L$ that are in the exterior of~$\CC_+$.
(Thus, we have~$\IN\,\cap\,\CC_-=\OUT\,\cap\,\CC_+=\emptyset$.)
Further, let~$A_\pm=\Lambda_L\setminus (\IN\cup \OUT)$ and use
$\sigma_\pm$ to denote the spin configuration on~$A_\pm$. 
Let~$M_\IN$,~$M_\OUT$ and~$M_\pm$ denote the overall 
magnetization in~$\IN$,~$\OUT$ and~$A_\pm$, respectively. Finally, 
let us abbreviate~$\mu_\IN=\lfloor\langle M_\IN\rangle_\IN^{+,\beta,s}\rfloor$ 
and introduce the event~$\EE_L'=\{\sigma\in\EE_L\colon M_\IN=-\mu_\IN\}$.

The lower bound on~$P_L^{+,\beta}(\MM_L)$ will be derived by
restricting to the event~$\EE_L'$, conditioning on~$\sigma_\pm$,
extracting the probability of having the correct magnetization 
in~$\Lambda_L\setminus A_\pm$, and applying
Lemma~\ref{lemma-Gauss} to retrieve the contribution from droplet
surface tension. The first  two steps of this program give
\begin{equation}
\label{summa} P_L^{+,\beta}(\MM_L)\ge
P_L^{+,\beta}(\MM_L\cap\EE_L')\ge\sum_{\sigma_\pm}
P_L^{+,\beta}(\MM_L\cap\EE_L'|\sigma_\pm)
P_L^{+,\beta}(\sigma_\pm).
\end{equation}
Our next goal is to produce a lower bound of the type \eqref{LBeq}
on~$P_L^{+,\beta}(\MM_L\cap\EE_L'|\sigma_\pm)$, uniformly
in~$\sigma_\pm$. The advantage of conditioning on a fixed
configuration is that, if~$\MM_L\cap\EE_L'\cap\{\sigma_\pm\}$ occurs, the overall
magnetizations in~$\IN$ and~$\OUT$ are fixed. Thus, on
$\MM_L\cap\EE_L'\cap\{\sigma_\pm\}$ we get
\begin{equation}
\label{MExt}
M_\OUT= M_L - M_{\pm} -M_{\IN}=\langle
M_\OUT\rangle_\OUT^{+,\beta,s}-2\mstar\, v_L\bigl(1-\lambda\Vol(\ssP) -
\delta_L\bigr),
\end{equation}
where~$\delta_L=\delta_L(\sigma_\pm)$ is given by the equation
$2\mstar\,v_L\delta_L=\text{\rm I}+\text{\rm II}+\text{\rm III}+\text{\rm IV}$ 
with~$\text{\rm I--IV}$ defined by
\begin{alignat}{2}
\text{\rm I}&=\mu_\IN-\mstar\, |\IN|,&\qquad\qquad 
\text{\rm II}&=-\langle
M_\OUT\rangle_\OUT^{+,\beta,s}+\mstar\, |\OUT|,\\
\text{\rm III}&=-M_{\pm}+\mstar\,|A_\pm|,
&\qquad 
\text{\rm IV}&=2\mstar\,
\bigl(|\IN|-\lambda\Vol(\ssP) v_L\bigr).
\end{alignat}
To estimate~$\text{\rm I--IV}$, we first notice the geometric bounds
\begin{equation}
\begin{array}{c}
t_L^2\Vol(\ssP)-t_L|\ssP|\le 
|\IN|\le (t_L+r_L)^2\Vol(\ssP)+(t_L+r_L)|\ssP|,
\\*[2mm]
%\intertext{}
|A_\pm|\le (t_L+3r_L)^2-t_L^2+(t_L+3r_L)|\ssP|,
\end{array}
\end{equation}
and  recall that, since both~$\CC_+$ and~$\CC_-$ are contained in~$A_\pm$, we have 
$|\CC_-|,|\CC_+|\le |A_\pm|$.
Lemma~\ref{L:M-fluctuations} for~$s=K\log L$ then allows us to estimate
$|\text{\rm I}|\le \alpha_1(\beta)
(|A_\pm| +|\IN|^2 L^{-\alpha_2(\beta)K})$
and, similarly,~$|\text{\rm II}|\le \alpha_1(\beta)
(|A_\pm| +4L+ L^{4-\alpha_2(\beta)K})$, while the remaining two quantities
are bounded  by invoking
$|\text{\rm III}|\le 2|A_\pm|$ and 
$|\text{\rm IV}|\le 4r_L t_L+2r_L^2+2(t_L+r_L)|\ssP|$.
Using that~$r_L=o(\sqrt{v_L})$ and~$t_L=O(\sqrt{v_L})$, we have~$|A_\pm|=o(v_L)$ 
as~$L\to\infty$. Moreover, if~$K$ is so large that 
$4-\alpha_2(\beta)K<4/3$, 
we also have~$|\IN|^2L^{-\alpha_2(\beta)K}\le L^{4-\alpha_2(\beta)K}
=o(v_L)$ as~$L\to\infty$. 
Combining these bounds, it is easy to show that~$|\delta_L(\sigma_\pm)|\le \bar 
\delta_L$ for all~$\sigma_\pm$, where $\bar\delta_L$ is a sequence 
such that~$\lim_{L\to\infty}\bar\delta_L=0$.

Now we are ready to estimate the probability that both~$\IN$ and
$\OUT$ produce their share of magnetization deficit. Note first that
\begin{equation}
\label{P>Ps}
P_\IN^{-,\beta}(M_\IN=-\mu_\IN)\ge
P_\IN^{-,\beta,s}(M_\IN=-\mu_\IN)
P_\IN^{-,\beta}\bigl(\Gamma_s(\sigma)=\emptyset\bigr).
\end{equation}
Using Lemmas~\ref{lemma-Gauss} and~\ref{lemma-nolog}, we get
$P_\IN^{-,\beta}(M_\IN=-\mu_\IN)\ge C{L^{-2/3}}$
for some~$C=C(\beta)>0$.
On the other hand, letting~$\MM_\OUT=\{\sigma\colon
M_\OUT=\langle
M_\OUT\rangle_\OUT^{+,\beta,s}-2\mstar\, v_L(1-\lambda \Vol(\ssP)-\delta_L)\}$,
a bound similar to \eqref{P>Ps} for~$P^{+,\beta}_{\OUT}$ 
combined with Lemmas~\ref{lemma-Gauss} and~\ref{lemma-nolog} yields
\begin{equation}
\label{Extbd}
P_\OUT^{+,\beta}(\MM_\OUT)\ge 
\frac{C'}{\sqrt{|\OUT|}}
\exp\Bigl\{-2\frac{(\mstar\,v_L)^2}{\chi|\OUT|}\bigl(1-\lambda\Vol(\ssP)
-\delta_L\bigr)^2\Bigr\},
\end{equation}
where~$C'=C'(\beta)>0$ is independent of~$\sigma_\pm$
contributing to \eqref{summa}. Combining the previous estimates,
we can use Lemma~\ref{circuit-lemma} to extract the surface energy term. 
The result is
\begin{equation}
\label{almost} P_L^{+,\beta}(\MM_L)\ge
C''L^{-5/3}\exp\bigl\{-w_1\sqrt{v_L}\,\Phi_L
-\epsilon_L'\sqrt{v_L}\bigr\},
\end{equation}
where~$C''=C''(\beta)>0$ and where~$\Phi_L$ stands for the quantity
\begin{equation}
%\label{}
\Phi_L=\frac{\mathscr{W}_\beta(\ssP)}{w_1}\sqrt\lambda+\frac{2(\mstar)^2\chi^{-1}w_1^{-
1}v_L^{3/2}}{L^2-(t_L+r_L)^2}
\bigl(1-\lambda\Vol(\ssP)+\bar\delta_L\bigr)^2.
\end{equation}
As is clear from our previous reasoning, the quantity~$\Phi_L$ can be made 
arbitrary close to~$\Phi_\Delta(\lambda)$ by letting~$L\to\infty$ and optimizing 
over
$\ssP$ with the above properties. The existence of the desired 
sequence~$(\epsilon_L)$ then follows by the definition of the limit. 
\end{proofsect}

\subsection{Results using random-cluster representation}
\label{sec3.2}\noindent
In this section we establish some technical results necessary for
the completion of the proof of our lower bound. These results are
stated mostly in terms of the random cluster counterpart of the
Ising model; the crowning achievement, which is
Lemma~\ref{FKcircuit-lemma}, gives immediately in the proof of
Lemma~\ref{circuit-lemma}. We remark that the latter is the sum
total of what this section contributes to the proof of
Theorem~\ref{lowerbound}. The uninterested, or well-informed,
readers are invited to skip the entire section, provided they are
prepared to accept Lemma~\ref{circuit-lemma} without a proof.

\subsubsection{Preliminaries}
The \emph{random cluster} representation for the Ising (and Potts) ferromagnets 
is by
now a well established tool. The purpose of the following remarks
is to define our notation; for more background and details we
refer the reader to, e.g.,~\cite{Grimmett,BBCK} or 
the excellent review~\cite{GHM}. 

Let~$\T\subset\Z^2$ denote a finite graph. A \emph{bond configuration},
generically denoted by~$\omega$, is the assignment of a zero
(vacant) or a one (occupied) to each bond in~$\T$. The weight of a
configuration~$\omega$ is given, informally, by
$R^{|\omega|}q^{C(\omega)}$, where~$|\omega|$ denotes the number
of occupied bonds and~$C(\omega)$ denotes the number of connected
components. For the Ising system at hand we have~$q=2$ and~$R=e^{2\beta}-1$. The
precise meaning of~$C(\omega)$ depends on the boundary conditions;
of concern here are the so called \emph{free} and \emph{wired} boundary
conditions. In the former,~$C(\omega)$ is the usual number of
connected components including the isolated sites, while in the
latter all clusters touching the bond-complement of~$\T$ are identified as a single component.

The free and wired random-cluster measures in~$\Lambda_L$, denoted by
$P_{L,\text{FK}}^{\text{free},\beta}$ and~$P_{L,\text{FK}}^{\text{w},\beta}$, respectively, 
correspond to the
free and plus (or minus) boundary conditions in the Ising spin
system. 
Both random-cluster measures enjoy the FKG property and
the wired measure  stochastically dominates the free measure. 
The infinite volume limits of these measures also exist; we denote these
limiting objects by~$P_{\text{FK}}^{\text{free},\beta}$ and~$P_{\text{FK}}^{\text{w},\beta}$.
The most important type of event we shall consider
is the event that sites are connected
by paths of occupied bonds. Our notation is as follows: If
$x,y\in\T$, we define~$\{x\longleftrightarrow y\}$ to be the event
that there is such a connection. If we demand the existence of a
path using only bonds with both ends in some subgraph~$\A\subset\T$, we
write~$\{x\underset{\A}\longleftrightarrow y\}$.

The next concept we need to discuss is \emph{duality}.
For any~$\T\subset\Z^2$, the \emph{dual} graph~$\T^*$ 
is defined as follows: Each
bond of~$\T$ is transversal to a bond on
$(\Z+\frac12)\times(\Z+\frac12)=(\Z^2)^*$. These bonds are the
bonds of~$\T^*$; the sites of~$\T^*$ are the endpoints of these
bonds. Each configuration~$\omega$ induces a configuration
on the dual graph via the correspondence 
``direct occupied''~with~``dual vacant'' 
and \emph{vice versa}. 
It turns out that, if we start with either free
or wired boundary conditions on~$\T$, the weights for the dual
configurations are also random-cluster weights with parameters
$(q^*,R^*)=(q,q/R)$, provided we also interchange the designation of
``free'' and ``wired.'' 
Of course, the graph and its dual are not precisely
the same. For example, if we examine the relevant graph for the problem dual to 
the wired system in~$\Lambda_L$, this consists of an~$(L+1)\times(L+1)$ rectangle
with the corners missing. Moreover, because the boundary conditions on the dual
graph are free, all dual edges touching the boundary sites are occupied 
independently
of the rest of the configuration. Thus, ignoring these decoupled degrees of
freedom, the restricted measure is equivalent to a free 
measure~on~$\Lambda_{L-1}$.

In general, we will use~$\beta^*$ to denote the inverse
temperature dual to~$\beta$, which, for~$q=2$ and the normalization of the Hamiltonian \eqref{Ham}, 
is related to~$\beta$ via
$\beta^*=\frac12\log\coth\beta$.
The critical temperature is self dual, i.e.,
$\betac=\frac12\log\coth\betac$. 
For~$\beta>\betac$, the dual model is in the high-temperature phase. Hence,
the limiting free and wired measures at~$\beta^*$ coincide and,
using the well-known relation between the spin-correlations 
and the connectivity functions in the FK~representation,
we have
\begin{equation}
\label{3.12}
P_{\text{FK}}^{\text{free},\beta^*}(x\longleftrightarrow y)=
P_{\text{FK}}^{\text{w},\beta^*}(x\longleftrightarrow y)
=\langle\sigma_0\sigma_x\rangle^{+,\beta^*},
\end{equation}
for all~$x,y\in\Z^2$. 
Thus, the exponential decay of correlations in the spin system at high 
temperatures,~$\langle\sigma_0\sigma_x\rangle^{+,\beta^*}\le e^{-\Vert x-
y\Vert/\xi}$
where~$\xi=\xi(\beta^*)$ is the correlation length, corresponds
to an exponential decay of the connectivity probabilities.
In particular, the \emph{surface tension} at~$\beta>\betac$, as defined in 
\eqref{tau-beta-n} for unit vectors~$\bn$ with rationally related 
components, is the inverse of the correlation length for two point connectivity 
functions in the direction~$\bn$ at inverse temperature~$\beta^*$.

\subsubsection{Decay estimates}
Here we assemble two important ingredients for the proof of Lemma~\ref{circuit-lemma}. 
We begin by quantifying the decay of
the point-to-boundary connectivity function: 

\begin{lemma}
\label{lemma-connection} Consider the~$q=2$ random cluster model
at~$\beta<\betac$ (which corresponds to the
high-temperature phase of the Ising system). Then, 
\begin{equation}
\label{eq3.10} 
P_{\ell,\text{\rm FK}}^{\text{\rm w},\beta}
\bigl(\{0\longleftrightarrow\partial\Lambda_\ell\}\bigr)
\le 4\ell e^{-\ell/\xi}
\end{equation}
for all~$\ell\ge 1$.
\end{lemma}

\begin{proofsect}{Proof}
This is one portion of the proof of Proposition~4.1 in \cite{CCFS}.
\end{proofsect}

For the purposes of the next lemma, let~$\bn$ be a unit vector
with rationally related components and let~$\CC(\bn)$ be the set
of all pairs~$(a,b)$ of positive real numbers such that the~$a\times b$
rectangle with side~$b$ perpendicular to~$\bn$ can be
positioned in~$\R^2$ in such a way that all its four corners are
in~$\Z^2$. We will use~$R_{a,b}^\bn\subset\Z^2$ to denote a
generic~$a\times b$ rectangle with the latter property. If~$x$ and~$y$ are the two 
corners along the same~$b$-side of~$R_{a,b}^\bn$, 
we let~$\BB_{a,b}^\bn$ denote the event 
$\{x\underset{R_{a,b}^\bn}{\longleftrightarrow} y\}$.

\begin{lemma}
\label{lemma-dualST} Let~$\beta\in(0,\betac)$ and let 
$\beta^*=\frac12\log\coth\beta$. 
Let~$\bn$ be a unit vector with
rationally related 
components and suppose that~$L$,~$a_L$ and~$b_L$, 
with~$(a_L,b_L)\in\CC(\bn)$, tend to
infinity in such a way that~$a_L/L\to0$,~$b_L/L\to0$ and
$\dist(R_{a,b}^\bn,\Z^2\setminus\Lambda_L)/(
b_L+\log L)\to\infty$ as~$L\to\infty$. Then
\begin{equation}
\label{surface-tension} \lim_{L\to\infty}P_{L,\text{\rm
FK}}^{\text{\rm free},\beta}\bigl(\BB_{a_L,b_L}^\bn\bigr)^{1/b_L}\ge
e^{-\tau_{\beta^*}(\bn)}.
\end{equation}
\end{lemma}

\begin{proofsect}{Proof}
We will first establish the limit \eqref{surface-tension} for the
measure in infinite volume and then show that provided
$R_L^\bn$ are well separated from
$\Z^2\setminus\Lambda_L$ as specified, the finite volume effects
are not important. Throughout the proof, we will omit the subscript~$\beta^*$ for
the surface tension.

Fix~$\bn\in\CalS_1$ with rationally related components and let
$\beta<\betac$. Let
\begin{equation}
\theta^{\,\bn}_{a,b}=P_{\text{FK}}^{\text{w},\beta}
\bigl(\BB_{a,b}^\bn\bigr), 
\qquad (a,b)\in\CC(\bn),
\end{equation}
and note that if~$(a,b_1)\in\CC(\bn)$ and
$(a,b_2)\in\CC(\bn)$ with~$b_2\ge b_1$, then also
$(a,b_1+b_2)\in\CC(\bn)$ and
$(a,b_2-b_1)\in\CC(\bn)$. We begin by the 
claim that the events in
question enjoy a subadditive property:
\begin{equation}
\label{subadd} \theta^{\,\bn} _{a,b_1+b_2}
\ge\theta^{\,\bn}_{a,b_1}\theta^{\,\bn}_{a,b_2}, \qquad (a,b_1),\,(a,b_2)\in\CC(\bn).
\end{equation}
Indeed, we let~$R_{a,b_2}^\bn$ be
translated relative to~$R_{a,b_1}^\bn$ so that the
``left''~$a$-side of~$R_{a,b_2}^\bn$ coincides with
the ``right''~$a$-side of~$R_{a,b_1}^\bn$. Let~$x_1$
and~$y_1$ be the ``left'' and ``right'' bottom corners of
$R_{a,b_1}^\bn$ and let~$x_2$ and~$y_2$ be similar
corners of~$R_{a,b_2}^\bn$. By our construction,~$y_1$
and~$x_2$ coincide. Let~$R_{a,b_1+b_2}^\bn$ denote the
union~$R_{a,b_1}^\bn\cup R_{a,b_2}^\bn$.
Then
\begin{equation}
\bigl\{x_1 \underset{R_{a,b_1+b_2}^\bn}{\longleftrightarrow} y_2\bigr\}\supset\bigl\{x_1
\underset{R_{a,b_1}^\bn}{\longleftrightarrow}
y_1\bigr\}\cap \bigl\{x_2 \underset{R_{a,b_2}^\bn}{\longleftrightarrow} y_2\bigr\}.
\end{equation}
The inequality \eqref{subadd} then follows immediately from the
FKG property of the measure~$P_{\text{FK}}^{\text{w},\beta}$.

Let~$\AA(\bn)=\{a>0\colon \exists b>0,\,
(a,b)\in\CC(\bn)\}$
be the set of allowed values of~$a$.
As a consequence of subadditivity, for any~$a\in\AA(\bn)$
we have the
existence of the limit~$e^{-\varpi_a(\bn)}=\lim_{b\to\infty}(\theta^{\,\bn}_{a,b})^{1/b}$. 
(Here~$b$ only takes 
values such that~$(a,b)\in\CC(\bn)$.)
Further, if~$a_1,a_2\in\AA(\bn)$ with~$a_1\ge a_2$,
then there \emph{is} a~$b$ such that both~$(a_1,b)\in\CC(\bn)$ and~$(a_2,b)\in\CC(\bn)$, 
and, for any such~$b$, we have~$\theta^{\,\bn}_{a_1,b}\ge\theta^{\,\bn}_{a_2,b}$. Thence
$\varpi_{a_1}(\bn)\le \varpi_{a_2}(\bn)$ whenever 
$a_1,a_2\in\AA(\bn)$ satisfy~$a_1\ge a_2$. Let
$\varpi(\bn)=
\lim_{a\to\infty}\varpi_a(\bn)$, where~$a$'s are restricted to 
$\AA(\bn)$. 
Now the quantity~$\theta^{\,\bn}_{\infty,b}=\lim_{a\to\infty}\theta^{\,\bn}_{a,b}$, where 
$(a,b)\in\CC(\bn)$, still obeys the subadditivity relation 
\eqref{subadd} and, in particular, the \emph{half-space} surface tension 
$\tau_{\text{h}}(\bn)$ is well defined by the limit
\begin{equation}
e^{-\tau_{\text{h}}(\bn)}=\lim_{b\to\infty}\lim_{\begin{subarray}{c}
(a,b)\in\CC(\bn)\\ a\to\infty
\end{subarray}}(\theta^{\,\bn}_{a,b})^{1/b}.
\end{equation}
Moreover,~$\theta^{\,\bn}_{\infty,b}\ge\theta^{\,\bn}_{a,b}$ for all~$a$ and~$b$ such that 
$(a,b)\in\CC(\bn)$ and, 
therefore,~$\tau_{\text{h}}(\bn)\le
\varpi(\bn)$. Our goal is to demonstrate that
$\tau_{\text{h}}(\bn)=\varpi(\bn)$ and that the half-space
surface tension~$\tau_{\text{h}}(\bn)$ equals the full space
surface tension~$\tau(\bn)$.

Let~$\epsilon>0$. Then there is a~$b^\star$ such that
$\theta^{\,\bn}_{\infty,b^\star}\ge
e^{-b^\star(\tau_{\text{h}}(\bn)+\epsilon)}$. However, since
$\theta^{\,\bn}_{\infty,b^\star}$ simply \emph{equals}
the limit of~$\theta^{\,\bn}_{a,b^\star}$ as
$a\to\infty$, there is an~$a^\star$ such that
$\theta^{\,\bn}_{a^\star,b^\star}\ge
e^{-b^\star(\tau_{\text{h}}(\bn)+2\epsilon)}$. Thence
$\varpi(\bn)\le\tau_{\text{h}}(\bn)$ and the equality of 
$\tau_{\text{h}}(\bn)$ and~$\varpi(\bn)$ follows. To remove 
the half-space constraint, consider the analogue of
the previously defined events. Let~$x$ and~$y$ be related to~$R_{a,b}^{\bn}$ 
as in the definition of event~$\BB_{a,b}^\bn$ 
and let~$D_{a,b}^\bn$
denote the union of~$R_{a,b}^\bn$ and its reflection
through the line joining~$x$ and~$y$. Let
\begin{equation}
\rho^{\,\bn}_{a,b}=P_{\text{FK}}^{\text{w},\beta}\bigl(\{x
\underset{D_{a,b}^\bn}{\longleftrightarrow} y\}\bigr).
\end{equation}
Reasoning identical to that employed thus far yields
\begin{equation}
e^{-\tau(\bn)}=\lim_{b\to\infty}\lim_{a\to\infty}
(\rho^{\,\bn}_{a,b})^{1/b}
=\lim_{a\to\infty}\lim_{b\to\infty}(\rho^{\,\bn}_{a,b})^{1/b},
\end{equation}
where we tacitly assume~$(a,b)\in\CC(\bn)$ for the production of both limits.
Now, obviously,~$\rho^{\,\bn}_{a,b}\ge\theta^{\,\bn}_{a,b}$ and
hence~$\tau(\bn)\le\tau_{\text{h}}(\bn)$. To derive the
opposite inequality, we note that for each~$a\in\AA(\bn)$, there is a 
$g(a)>0$
such that
\begin{equation}
\label{lastbd} \theta^{\,\bn}_{2a,b}\ge g(a)\rho^{\,\bn}_{a,b},\qquad
(a,b)\in\CC(\bn).
\end{equation}
Indeed, the event giving
rise to~$\theta^{\,\bn}_{2a,b}$ can certainly be
achieved by connecting the bottom corners 
of~$R_{2a,b}^\bn$ directly to the middle points and then
connecting the middle points on the opposite~$a$-sides of
$R_{2a,b}^\bn$. 
Then \eqref{lastbd} follows by FKG. (To get that~$g(a)>0$, we
also used that~$\beta>0$.)
Taking the~$1/b$-th power of both sides of \eqref{lastbd} and
letting~$b\to\infty$ followed by~$a\to\infty$ we arrive at~$\varpi(\bn)=\tau_{\text{h}}(\bn)=
\tau(\bn)$ as promised.

To finish the proof, we must account for the effects of finite
volume. Consider the event~$\FF_{a,b}^\bn=\{\partial
R_{a,b}^\bn\leftrightarrow\partial\Lambda_L\}$. Should~$\FF_{a,b}^\bn$
\emph{not} occur, a vacant ring separates~$R_{a,b}^\bn$ from~$\partial\Lambda_L$ and, 
using fairly standard arguments, we have
\begin{equation}
\label{circuittrick}
P_{L,\text{FK}}^{\text{free},\beta}
(\BB_{a,b}^\bn)
\ge
P_{\text{FK}}^{\text{w},\beta}
\bigl(\BB_{a,b}^\bn\big|(\FF_{a,b}^{\bn})^{\text{c}}\bigr).
\end{equation}
On the other hand, by Lemma~\ref{lemma-connection}, we have
\begin{equation}
P_{\text{FK}}^{\text{w},\beta}(\FF_{a,b}^\bn)\le 
P_{L,\text{FK}}^{\text{w},\beta}(\FF_{a,b}^\bn)\le
8L(a+b)\,e^{-\dist(\partial R_{a,b}^\bn,\partial\Lambda_L)/\xi}.
\end{equation}
Thus if the distance between~$\partial R_{a,b}^\bn$ and~$\partial\Lambda_L$ exceeds a large multiple of 
$b_L+\log L$, the dominant contribution to~$P_{\text{FK}}^{\text{w},\beta}
(\BB_{a,b}^\bn)$ comes from~$P_{\text{FK}}^{\text{w},\beta}
(\BB_{a,b}^\bn\big|(\FF_{a,b}^\bn)^{\text{c}})$. Using 
\eqref{circuittrick}, the claim follows. 
\end{proofsect}

\subsubsection{Corona estimates}
We recall the ``corona'' regions 
$\K_{t,r}^{\text{\rm I}}$--$\K_{t,r}^{\text{\rm III}}$
associated with some given polygon~$\ssP$. In addition, we will
also need to consider the collection of dual sites
$\K_{t,r}^{*\text{\rm II}}=K_{t,r}^{\text{\rm II}}\cap(\Z^2)^*$,
where~$(\Z^2)^*$ is the lattice dual to~$\Z^2$. (This differs slightly from
the graph dual to~$\K_{t,r}^{\text{\rm II}}$ by some boundary sites.)
In the context of the random cluster model (and its dual) we will consider three
events: The first event, to be denoted~$\EE^{\text{\rm I}}_{t,r}$, 
takes place in~$\K^{\text{\rm I}}_{t,r}$ and is defined by
\begin{equation}
\EE^{\text{\rm I}}_{t,r}=\bigl\{\omega\in\Omega\colon
\emph{ there is a circuit of occupied bonds in } \K^{\text{\rm I}}_{t,r} 
\emph{ surrounding the origin}\bigr\}.
\end{equation}
The event~$\EE^{\text{\rm III}}_{t,r}$ is defined similarly
except that the circuit takes place in the region~$\K^{\text{\rm III}}_{t,r}$. 
Finally, one more circuit, this time a dual circuit in the region 
$\K^{\text{\rm II}*}_{t,r}$. We define
\begin{equation}
\EE^{\text{\rm II}*}_{t,r}=\bigl\{\omega\in\Omega\colon
\emph{ there is a dual circuit of vacant bonds in } \K^{*\text{\rm II}}_{t,r} 
\emph{ surrounding the origin}\bigr\}.
\end{equation}
As we will see in the proof of Lemma~\ref{circuit-lemma}, the event
$\EE^{\text{\rm I}}_{t,r}\cap\EE^{\text{\rm II}*}_{t,r}\cap\EE^{\text{\rm III}}_{t,r}$ 
more or less implies the desired event~$\EE_{t,r}$. The desired lower
bound will then be an immediate consequence of the following lemma:

\begin{lemma}
\label{FKcircuit-lemma} Let~$\beta>\betac$ and let~$\ssP$ be as in 
Lemma~\ref{circuit-lemma}. 
For any sequences~$(t_L)$ and~$(r_L)$ satisfying 
\eqref{t-r-cond}, there is a sequence~$(\epsilon_L'')$ such that
$\epsilon_L''\to0$ and, for all~$L$,
\begin{equation}
%\label{}
P_{L,\text{\rm FK}}^{\text{\rm w},\beta}
\bigl(\EE^{\text{\rm I}}_{t_L,r_L}\cap
\EE^{\text{\rm II}*}_{t_L,r_L}\cap
\EE^{\text{\rm III}}_{t_L,r_L}\bigr)\ge
\exp\bigl\{-t_L\mathscr{W}_\beta(\ssP)(1+\epsilon_L'')\bigr\}.
\end{equation}
\end{lemma}

\begin{proofsect}{Proof}
In the course of this proof, let us abbreviate~$\EE^{\text{\rm
I}}_L=\EE^{\text{\rm I}}_{t_L,r_L}$, and similarly for
$\EE^{\text{\rm II}*}_L$ and~$\EE^{\text{\rm III}}_L$, as well as
$\K_L^{\text{\rm I}}$,~$\K_L^{*\text{\rm II}}$, and
$\K_L^{\text{\rm III}}$. We will start with an estimate for
$P_{L,\text{\rm FK}}^{\text{\rm w},\beta}(\EE^{\text{\rm
II}*}_L)$, which is in any case the central ingredient of this
lemma. Let~$T$ be the smallest integer~$T\ge2$ such that the polygon~$\ssP$
magnified by~$T$ has all vertices on~$\Z^2$. Let~$u_L=T\lfloor 
(t_L+r_L)/T\rfloor+T$ and let~$x_1,\dots,x_N$ be 
the vertices of the polygon~$\ssP$ magnified by~$u_L$. 
Let~$x_1^*,\dots,x_N^*$ be the corresponding vertices of the 
polygon~$\ssP$ magnified by~$u_L$ and translated by~$(-\tfrac12,-\tfrac12)$.
Notice that (once~$t_L$ and~$r_L$ are large enough) 
the sites~$x_1^*,\dots,x_N^*$
lie inside the ``corona''~$\K_L^{*\text{\rm II}}$.
We use~$\bn_i$ to denote the unit
vector constituting the outer normal to the side between~$x_{i+1}^*$ 
and~$x_i^*$ (where~$x_{N+1}^*$ is identified with~$x_1^*$). By our
construction,~$x_1,\dots,x_N\in\Z^2$,~$x_1^*,\dots,x_N^*\in(\Z^2)^*$ and~$\bn_i$ have rationally
related~components.

For~$i=1,\dots, N$, let us consider the rectangles
$R_{a_i,b_i}^{\bn_i}$ with the base coinciding with the
line between~$x_i^*$ and~$x_{i+1}^*$. Here~$a_i$ is the largest
possible number such that~$(a_i,b_i)\in\CC(\bn_i)$ and
$R_{a_i,b_i}^{\bn_i}\subset\K_L^{*\text{\rm II}}$. 
We remark that all~$(a_i)$ and~$(b_i)$ have
$L$-dependence which is notationally suppressed and that these
tend to infinity as~$L\to\infty$. In particular, the~$b_i$'s scale
with~$u_L$. Let us denote
\begin{equation}
\label{scal-rel} \frakb_i=\lim_{L\to\infty}\frac{b_i}{t_L},\qquad i=1,\dots,N,
\end{equation}
where the limit exists by the construction of~$b_i$'s and where we noted that 
$t_L/u_L\to1$ as~$L\to\infty$.

Let~$\BB_i^*$ be the event that there is a dual vacant connection 
$x_i^*\longleftrightarrow x_{i+1}^*$ in the box 
$R_{a_i,b_i}^{\bn_i}$
and let~$\BB_i$ be the corresponding ``direct'' event that there is
a direct occupied path
$x_i\longleftrightarrow x_{i+1}$ contained in~$(\tfrac12,\tfrac12)$-translate of 
$R_{a_i,b_i}^{\bn_i}$.
It is clear that the intersection~$\bigcap_{i=1}^N\BB_i^*$
produces the event~$\EE_L^{\text{\rm II}*}$ and that these events are
FKG-correlated. Moreover, by duality, we have
\begin{equation}
P_{L,\text{FK}}^{\text{w},\beta}(\BB_i^*)= P_{L-1,\text{FK}}^{\text{free},\beta^*}
(\BB_i)
\end{equation}
(c.f., the paragraph before \eqref{3.12}).
Now we are perfectly positioned to apply
Lemma~\ref{lemma-dualST}: Using FKG, the scaling relation
\eqref{scal-rel}, and the fact that also the~$a_j$'s tend to
infinity by our construction, we have as a consequence of the
above-mentioned lemma that
\begin{equation}
\label{lim-eq} \lim_{L\to\infty}
P_{L,\text{FK}}^{\text{w},\beta}\bigl(\EE_L^{\text{\rm II}*}\bigr)
^{1/t_L}= \exp\Bigl\{-\sum_{j=1}^N \frakb_j\tau_\beta(\bn_j)\Bigr\}.
\end{equation}

The remainder of the proof concerns the estimate of the probability
$P_{L,\text{FK}}^{\text{w},\beta}(\EE_L^{\text{\rm I}}
\cap\EE_L^{\text{\rm III}}|\EE_L^{\text{\rm II}*})$. We claim that this
conditional probability tends to one as~$L\to\infty$. 
First, as a worst-case
scenario, consider the event~$V_L^{\text{\rm II}*}$ that all bonds in
$\K_L^{*\text{\rm II}}$ are vacant. By monotonicity in boundary
conditions and the strong FKG property of~$P_{L,\text{FK}}^{\text{w},\beta}$ it is 
seen that
\begin{equation}
P_{L,\text{FK}}^{\text{w},\beta}\bigl(\EE_L^{\text{\rm I}}
\cap\EE_L^{\text{\rm III}}\big|\EE_L^{\text{\rm II}*}\bigr)\ge
P_{L,\text{FK}}^{\text{w},\beta}\bigl(\EE_L^{\text{\rm I}}
\cap\EE_L^{\text{\rm III}}\big|V_L^{\text{\rm II}*}\bigr).
\end{equation}
Under the condition that~$V_L^{\text{\rm II}*}$ occurs, 
$\EE_L^{\text{\rm I}}$ and
$\EE_L^{\text{\rm III}}$ are independent and we may treat them
separately. The arguments are virtually identical
for both events, so we need only be explicit 
about $P_{L,\text{FK}}^{\text{w},\beta}(\EE_L^{\text{\rm I}}
|V_L^{\text{\rm II}*})$.

Let~$\ell_L$ be a maximal integer such that there is a circuit of dual cites, 
$z_1^*,\dots,z_m^*$, separating the boundaries of 
$\K_L^{\text{\rm I}}$ with the property that,
if~$\Lambda_{\ell_L}^*(z_j^*)$ is the translate of
$\Lambda_{\ell_L}^*$ by (the vector)~$z_j^*$, then
$\Lambda_{\ell_L}^*(z_j^*)\subset\K_L^{\text{\rm I}}$.
Note that~$\liminf_{L\to\infty}\ell_L/r_L>1/3$.
Now, for the event~$\EE_L^{\text{\rm I}}$
\emph{not} to occur, there must be a dual occupied path connecting some
dual site on the outer boundary of~$\K_L^{\text{\rm I}}$ to another on
the inner boundary and hence at least one~$z_j^*$ has
to be connected to the boundary of its~$\Lambda_{\ell_L}^*(z_j^*)$
by a path of dual occupied bonds. Using subadditivity of the probability measure, 
we find
\begin{equation}
1-P_{L,\text{FK}}^{\text{w},\beta}\bigl(\EE_L^{\text{\rm I}}
\big|V_L^{\text{\rm II}*}\bigr)\le\sum_{j=1}^m
P_{L,\text{FK}}^{\text{w},\beta}\bigl(
z_j^*\longleftrightarrow\partial\Lambda_{\ell_L}^*(z_j^*)
\big|V_L^{\text{\rm II}*}\bigr).
\end{equation}
Now, again invoking monotonicity in the boundary conditions, the
probability of the above connection events may be estimated from
above by placing dual wired (i.e., direct free) boundary
conditions on~$\Lambda_{\ell_L}^*(z_j^*)$. But then, by duality,
we have exactly the event which is the subject of
Lemma~\ref{lemma-connection}. Explicitly,
\begin{equation}
P_{L,\text{FK}}^{\text{w},\beta}\bigl(
z_j^*\longleftrightarrow\partial\Lambda_{\ell_L}^*(z_j^*)
\big|V_L^{\text{\rm II}*}\bigr)\le
P_{\ell_L,\text{FK}}^{\text{w},\beta^*}\bigl(
0\longleftrightarrow\partial\Lambda_{\ell_L}\bigr)
\end{equation}
holds for all~$j=1,\dots,m$, and the bound in \eqref{eq3.10} can
be applied. Now the number of sites~$z_j^*$ which comprise the
circuit does not exceed a multiple of~$t_L$. Thus, for some
constant~$C$ independent of~$L$ we have
\begin{equation}
P_{L,\text{FK}}^{\text{w},\beta}\bigl(\EE_L^{\text{\rm I}}
\big|V_L^{\text{\rm II}*}\bigr)\ge1-C \ell_Lt_Le^{-\ell_L/\xi}.
\end{equation}
By the condition stated in \eqref{t-r-cond}, the fact that~$r_L\ge\ell_L\ge r_L/3$ 
for sufficiently large~$L$, and the observation that~$\xi^{-1}=\tau_\mini$, the 
desired
result for~$\EE_L^{\text{\rm I}}$ follows. Similarly for
$\EE_L^{\text{\rm III}}$. 
\end{proofsect}

\begin{proofsect}{Proof of Lemma~\ref{circuit-lemma}}
We make liberal use of the correspondence between the graphical
configurations~$\omega$ and (sets of) spin configurations as
described, e.g., in~\cite{ACCN,ES,BBCK}. Each connected cluster in
$\omega$ represents the spin configurations in which all sites of
the cluster have spins of the same type. Thus, if
$\EE_L^{\text{\rm I}}\cap\EE_L^{\text{\rm II}*}\cap\EE_L^{\text{\rm III}}$
occurs, then the inner circuit of occupied bonds in~$\K_L^{\text{\rm I}}$
forces the spins on these sites to be of the same type. Since
these are disconnected from the boundary of~$\Lambda_L$ by the
dual vacant circuit in~$\K_L^{*\text{\rm II}}$, with probability
one-half, all spins on the circuit are minus. Similarly, the outer
circuit of bonds in~$\K_L^{\text{\rm III}}$ is plus-type with
probability one if it is connected to~$\partial\Lambda_L$ and
with probability~$1/2$ otherwise. Thus,
$P_L^{+,\beta}(\EE_{t_L,r_L}|\EE_L^{\text{\rm I}}
\cap\EE_L^{\text{\rm II}*}\cap\EE_L^{\text{\rm III}})$ is certainly
bigger than~$1/4$, and the claim follows using Lemma~\ref{FKcircuit-lemma}. 
\end{proofsect}

\section{Absence of intermediate contour sizes}
\label{sec4}
\subsection{Statement and outline}
The goal of this section is to 
prove that, with probability tending to one as~$L\to\infty$, there will be no 
contours with a diameter between the scales of~$\log L$ and~$\sqrt{v_L}$ in the 
``canonical'' ensemble of the Ising model in volume~$\Lambda_L$. This result is by far 
the most difficult part of the proof of our main results stated in 
Section~\ref{sec1.3}. 

We start with a standard notion from contour theory. 
Let~$\Gamma(\sigma)$ denote the set of all contours of a configuration~$\sigma$ in~$\Lambda_L$ 
with plus boundary condition.
Applying the rounding rule, contours are self-avoiding simple curves in~$\R^2$. 
Recall that~$\Gamma_s(\sigma)$ is the set of contours of~$\sigma$ that have a non-trivial~$s$-skeleton.
We say that~$\gamma\in\Gamma(\sigma)$ is an \emph{external} contour, if  
it is not surrounded by any other contour from~$\Gamma$.
We will use~$\Gamma_s^\ext(\sigma)$ to denote the set of external contours of~$\Gamma_s(\sigma)$.
(We remark that~$\Gamma_s^\ext(\sigma)$, namely the external contours of~$\Gamma(\sigma)$
which are big enough to have an~$s$-skeleton, coincides exactly with the set of external contours of the
collection~$\Gamma_s(\sigma)$.)

Using this notation, the event~$\AA_{\varkappa,s,L}$ from Theorem~\ref{main-result} 
is best described via its complement:
\begin{equation}
%\label{}
\AA^{\text{\rm c}}_{\varkappa,s,L}=\bigl\{\sigma\colon 
\exists\gamma\in\Gamma_s^\ext(\sigma),\,
\diam\gamma\le\varkappa\sqrt{v_L}\bigr\}.
\end{equation}
The relevant claim is then restated as follows:

\begin{theorem}
\label{upperbound}
Let~$\beta>\betac$ and let~$(v_L)$ be a sequence of positive 
numbers that make~$\mstar|\Lambda_L|-2\mstar\,v_L$ an allowed value of~$M_L$ for 
all~$L$. Suppose the limit~$\Delta$ in \eqref{Delta-lim} obeys 
$\Delta\in(0,\infty)$. 
For each~$c_0>0$ there exist $\varkappa>0$, $K_0<\infty$ and $L_0<\infty$ such that
if $K\ge K_0$, $L\ge L_0$ and $s=K\log L$, then
\begin{equation}
\label{intermediate}
P_L^{+,\beta}\bigl(\AA_{\varkappa,s,L}^{\text{\rm c}}\big|
M_L=\mstar|\Lambda_L|-2\mstar\,v_L\bigr)\le L^{-c_0}
\end{equation}
\end{theorem}

\smallskip
Let~$s=K\log L$ be a scale function and recall that a 
contour~$\gamma$ is \textit{$s$-large} if~$\gamma\in\Gamma_s(\sigma)$. 
For $\varkappa>0$, a contour~$\gamma$ large enough to be an~$s$-large contour but satisfying
$\diam\gamma\le\varkappa\sqrt{v_L}$ will be called a \textit{$\varkappa$-intermediate} contour.
Thus, Theorem~\ref{upperbound} shows that, 
in the canonical ensemble with the magnetization fixed to~$\mstar|\Lambda_L|-
2\mstar\,v_L$, there are no~$\varkappa$-intermediate contours with probability tending to one 
as~$L$ tends to infinity. 
This statement, which is of interest in its own right, reduces the 
proof of our main result to a straightforward application of isoperimetric 
inequalities for the Wulff functional as formulated in Lemma~\ref{lemma-onecontour}.

\begin{remark}
The reason why a \textit{power} of~$L$ appears on the 
right-hand side is because we only demand the absence of contours with sizes 
over~$K\log L$.
Indeed, for a general~$s$,  the right-hand side of \eqref{intermediate} could be replaced by
$e^{-\alpha  s}$ for some constant~$\alpha>0$. 
In particular, the decay can be made substantially faster 
by easing the lower limit of what we chose to call an intermediate size contour.
Finally, we note that~$L_0$ in Theorem~\ref{upperbound} depends not only on  
$\beta$, $\Delta$, and  $c_0$, 
but also  on how fast the  limit~$v_L^{3/2}/|\Lambda_L|$ is achieved.
\end{remark}

The proof of Theorem~\ref{upperbound} will require some preparations. In particular, 
we will need to estimate the (conditional) probability of five highly unprobable events  
that we would like to exclude explicitly from the further considerations.
All five events are defined with reference to a positive number~$\varkappa$ which, more or less, 
is the same~$\varkappa$
that appears in Theorem~\ref{upperbound}. 

The first event,~$\RR^1_{\varkappa,s,L}$, 
collects the configurations for which the combined length of all~$s$-large contours 
in~$\Lambda_L$ exceeds~$\varkappa^{-1}s\sqrt{v_L}$.
These configurations need to be \textit{a priori} excluded 
because all of the crucial Gaussian estimates from Section~\ref{sec2.3} can only be applied to 
regions  with a moderate  surface-to-volume ratio.
Next, we show that one can ignore configurations whose large contours
occupy too big volume.
This is the basis of the event~$\RR^2_{\varkappa,s,L}$.
The remaining  three events concern the magnetization deficit in two random subsets of 
$\Lambda_L$: A set~$\oInnt\subset\V(\Gamma_s^\ext(\sigma))$ of sites enclosed by an~$s$-large contour 
and a  set~$\oExxt$ of sites outside all~$s$-large contours. 
The precise 
definitions of these sets is given
in Section~\ref{sec4.2}. 
The respective events~are:
\begin{enumerate}
\item[(3)]
The event~$\RR^3_{\varkappa,s,L}$ that 
$M_{\oInnt}\le -\mstar|\oInnt|-\varkappa^{-1}s v_L^{3/4}$.
\item[(4)]
The event 
$\RR^4_{\varkappa,s,L}$ that~$M_{\oExxt}\ge\mstar|\oExxt|-2\varkappa\mstar 
v_L$.
\item[(5)]
The event~$\RR^5_{\varkappa,s,L}$ that 
$M_{\oExxt}\le\mstar|\oExxt|-2(1+\varkappa^{-1})\mstar v_L$.
\end{enumerate}
By choosing~$\varkappa$ sufficiently small, the 
events~$\RR^1,\dots,\RR^5$ will be shown to have a probability vanishing exponentially fast 
with~$\sqrt{v_L}$. These estimates are the content of Lemma~\ref{lemma-length} and 
Lemmas~\ref{lemma-nomag}-\ref{lemma-muchmag}.

Once the preparatory statements have been proven, we consider a rather extreme version 
of the restricted contour ensemble, namely, one in which no contour that is larger
than~$\varkappa$-intermediate is allowed to appear. We show, in a rather difficult
Lemma~\ref{lemma-intermediate-help}, that despite this restriction, bounds similar to those of
\eqref{intermediate} still hold.  
The final step---the proof of Theorem~\ref{upperbound}---is now achieved by conditioning
on the location(s) of the large contour(s), which by the ``$\RR$-lemmas'' are typically not \emph{too} big and
not \emph{too} rough. By definition, the exterior region is now in the restricted ensemble featured in
Lemma~\ref{lemma-intermediate-help} and the result derived therein allows a relatively easy endgame.

\smallskip
Throughout Sections~\ref{sec4.2}-\ref{sec4.4} we will let~$\beta>\betac$ 
be fixed and let~$(v_L)$ be a sequence of positive numbers such that
$\mstar|\Lambda_L|-2\mstar\,v_L$ is an allowed value of~$M_L$ for all~$L$.
Moreover, we will assume that~$(v_L)$ is such that the limit~$\Delta$ in \eqref{Delta-lim} exists with
$\Delta\in(0,\infty)$. 

\subsection{Contour length and volume}
\label{sec4.2}\noindent
In this section we will prepare the grounds for the proof
of Theorem~\ref{upperbound}. In particular, we derive rather
crude estimates on the total length of large contours
and the volume inside and outside large external
contours. These results come as 
Lemmas~\ref{lemma-length} and~\ref{lemma-volume} 
below. 

\subsubsection{Total contour length}
We begin by estimating the combined length of large contours. Let~$s$ 
be a scale function and, for any~$\varkappa>0$, let
$\RR^1_{\varkappa,s,L}$ be the event
\begin{equation}
\RR^1_{\varkappa,s,L}=\Bigl\{\,\sigma\colon\!\!
\sum_{\gamma\in\Gamma_s(\sigma)}|\gamma|\ge \varkappa^{-1}
s\sqrt{v_L}\Bigr\}.
\end{equation}
The probability of event~$\RR^1_{\varkappa,s,L}$ is then estimated as
follows:

\begin{lemma}
\label{lemma-length} For each~$c_1>0$ 
there exist~$\varkappa_0>0$,~$K_0<\infty$ and~$L_0<\infty$
such that 
\begin{equation}
\label{A1bd} P_L^{+,\beta}
\bigl(\RR^1_{\varkappa,s,L}\big|M_L=\mstar\,|\Lambda_L|-2\mstar\,
v_L\bigr)\le e^{-c_1\sqrt{v_L}}
\end{equation}
holds for all~$\varkappa\le\varkappa_0$,~$K\ge K_0$,~$L\ge
L_0$, and~$s=K\log L$.
\end{lemma}

\begin{proofsect}{Proof}
Let~$K_0$ be the quantity~$K_0(\frac12,\beta)$ from
Lemma~\ref{lemma-Peierls} and let us recall that~$\tau_\mini$ denotes the
minimal value of the surface tension. We claim that it suffices to
show that, for all~$c_1'>0$ and an appropriate choice of~$\varkappa$, the bound
\begin{equation}
\label{jinej} P_L^{+,\beta}(\RR^1_{\varkappa,s,L})\le
e^{-c_1'\sqrt{v_L}}
\end{equation}
holds true once~$L$ is sufficiently large.
Indeed, if \eqref{jinej} is established, we just choose~$c_1'$ so large that the 
difference~$c_1'-c_1$ exceeds the rate constant from the lower bound in 
Theorem~\ref{lowerbound} and the estimate \eqref{A1bd} immediately follows.

In order to prove \eqref{jinej}, fix~$c_1'>0$ and let~$\varkappa_0^{-1}
= 2g_1c_1'/\tau_\mini$, where~$g_1$ is as in \eqref{triv-s}. Let~$K\ge K_0$, 
$\varkappa\le\varkappa_0$ and
$s=K\log L$. 
We claim that if~$\sigma\in\RR^1_{\varkappa,s,L}$
and~$\frakS$ is a collection of~$s$-skeletons such that 
$\frakS\sim\sigma$, then \eqref{triv-s} and 
\eqref{triv-W} force
\begin{equation}
\varkappa^{-1}
s\sqrt{v_L}\le\sum_{\gamma\in\Gamma_s(\sigma)}|\gamma|\le g_1s
\sum_{S\in\frakS}\bigl|\ssP(S)\bigr|\le g_1s\tau_\mini^{-1}
\mathscr{W}_\beta(\frakS).
\end{equation}
Hence, for each~$\sigma\in\RR^1_{\varkappa,s,L}$ there is at least one 
$\frakS$ such that~$\frakS\sim\sigma$ and~$\mathscr{W}_\beta(
\frakS)\ge2c_1'\sqrt{v_L}$. By Corollary~\ref{cor-help} with 
$\kappa =2c_1'\sqrt{v_L}$ 
and~$\alpha=\frac12$, and our choice of~$K_0$, 
\eqref{jinej} follows. 
\end{proofsect}

\subsubsection{Interiors and exteriors}
Given a scale function 
$s$ and a configuration~$\sigma$, 
let~$\Gamma_s^\ext(\sigma)$ be the set of 
external contours in 
$\Gamma_s(\sigma)$. (Note that these contours will also be 
external
in the set of all contours of~$\sigma$.) Define~$\Innt=\Innt_{s,L}(\sigma)$ 
to 
be the set of all
sites in~$\Lambda_L$ enclosed by some 
$\gamma\in\Gamma_s^\ext(\sigma)$
and let~$\Exxt=\Exxt_{s,L}(\sigma)$ be the 
complement of~$\Innt$, i.e.,~$\Exxt=\Lambda_L\setminus\Innt$. 

Given a set of external contours~$\Gamma$, 
we claim that under the condition that~$\Gamma_s^\ext(\sigma)=\Gamma$, the measure~$P_L^{+,\beta}$
is a product of independent measures on~$\Exxt$ and~$\Innt$. 
A coarse look might suggest a product of
plus-boundary condition measure on~$\Exxt$ and the minus measure on~$\Innt$. 
Indeed, all spins in~$\Exxt$ up against a piece of~$\Gamma$ are 
necessarily pluses and similarly all spins on the~$\Innt$ sides of these contours are minuses.
But this is not quite the end of the story, two small points are in order: 
First, we have invoked a rounding rule.
Thus, for example, certain spins in~$\Exxt$ (at some corners but not up against the contours) are
\emph{forced} to be plus otherwise the rounding rule would have drawn the contour differently. On the other
hand, some corner spins \emph{are} permitted either sign because the rounding rule would separate any such
resulting contour. Fortunately, the upshot of these ``rounding anomalies'' is only to force a few additional
\emph{minus} spins in~$\Innt$ and
\emph{plus} spins in~$\Exxt$ than would 
appear from a naive look at~$\Gamma$.

To make the aforementioned observations notationally apparent, we define~$\oInnt\subset\Innt$ 
to be the set of sites that can be flipped without changing~$\Gamma$ and similarly for~$\Exxt$. 
We thus have~$\sigma_x=-1$ for all~$x\in\Innt\setminus\oInnt$ and~$\sigma_x=+1$ 
for all~$x\in\Exxt\setminus\oExxt$.
Explicitly, there are a few more boundary spins than one might have thought, 
but they are always of the correct type.
Thus, clearly, although rather trivially, 
the measure~$P_L^{+,\beta}(\cdot|\Gamma_s^\ext(\sigma)=\Gamma)$ restricted to~$\Innt$ 
is simply the measure in~$\Innt$ with minus boundary conditions. 
The same measure on~$\Exxt$ is not quite the corresponding plus-measure due to the condition
that~$\Gamma$ constitutes \emph{all} the external contours visible on 
the scale~$s$. 
Thus, beyond the scale~$s$ in~$\Exxt$, we must see\dots{}no contours. 
But this is precisely the definition of the restricted ensemble.

We conclude that the conditional measure splits on~$\Innt$ and~$\Exxt$ 
into independent measures that are well understood. 
Explicitly, if~$\AA$ is an event depending only on the spins 
in~$\oInnt$ and~$\BB$ is an event depending only on 
the spins in~$\oExxt$, then
\begin{equation}
\label{Int-Ext}
P_L^{+,\beta}\bigl(\AA\cap\BB\big|\Gamma_s^\ext(\sigma)=\Gamma\bigr)=
P_{\oInnt}^{-,\beta}(\AA)P_{\oExxt}^{+,\beta,s}(\BB).
\end{equation}
This observation will be crucial for our estimates in the next section.

Next we will notice that the number of sites associated with the contours
can be easily bounded in terms of the total length of~$\Gamma$:

\begin{lemma}
\label{lemma-corners}
There exists a geometrical constant~$g_4<\infty$ such that the following is true: 
If~$\Gamma$ is a set of external contours and $\oInnt$ and $\oExxt$ are as defined above, then
\begin{equation}
\label{corners}
|\Lambda_L\setminus(\oInnt\cup\oExxt)|\le g_4\sum_{\gamma\in\Gamma}|\gamma|.
\end{equation}
\end{lemma}

\begin{proofsect}{Proof}
Each site from $\Lambda_L\setminus(\oInnt\cup\oExxt)$ is within 
some (Euclidean) distance from a
dual lattice site $x^*\in(\Z^2)^*$ such that some contour $\gamma\in\Gamma$ passes through~$x^*$.
On the other hand, the number of dual lattice sites~$x^*$ visited by contours from~$\Gamma$ 
does not exceed twice the total length of all contours in $\Gamma$. From here the existence of a $g_4$ satisfying
\eqref{corners} follows. 
\end{proofsect}

\smallskip
The definition of the event~$\RR^1_{\varkappa,s,L}$ gives us the following easy 
bounds:

\begin{lemma}
\label{lemma-volume}
Let $g_4$ be as in Lemma~\ref{lemma-corners}. Let~$\sigma\not\in\RR^1_{\varkappa,s,L}$ 
and let the sets~$\Innt=\Innt_{s,L}(\sigma)$,
$\oInnt=\oInnt_{s,L}(\sigma)$  and~$\oExxt=\oExxt_{s,L}(\sigma)$
be as above. Then 
we have the bounds
\begin{equation}
\label{bd-int-bd}
|\partial\oInnt|\le g_4\varkappa^{-1}s\sqrt{v_L}
\quad\text{and}\quad
|\partial\oExxt|\le g_4\varkappa^{-1}s\sqrt{v_L}
+4L
\end{equation}
and
\begin{equation}
\label{int-bd}
|\oInnt|\le|\Innt|\le g_4^2\varkappa^{-2}s^2 v_L.
\end{equation}
\end{lemma}

\begin{proofsect}{Proof}
Since $\partial\oInnt\subset\Lambda_L\setminus(\oExxt\cup\oInnt)$ which by Lemma~\ref{lemma-volume} implies $|\partial\oInnt|\le g_4\sum_{\gamma\in\Gamma_s(\sigma)}|\gamma|$, the first bound in \eqref{bd-int-bd} is an immediate consequence of the fact that~$\sigma\not\in\RR^1_{\varkappa,s,L}$.
Note that the same inequality is true for $|\partial\Innt|$.
The second bound in \eqref{bd-int-bd} then follows by the fact that $\partial\oExxt\subset\partial\Lambda_L\cup\Lambda_L\setminus(\oExxt\cup\oInnt)$.
The last bound, \eqref{int-bd}, is then implied by the first bound in \eqref{bd-int-bd} for $\partial\Innt$ instead of $\partial\oInnt$ and the isoperimetric inequality~$|\Lambda|\le\frac1{16}|\partial\Lambda|^2$ valid for any $\Lambda\subset\R^2$ that is a finite union of closed unit squares (see, e.g., Lemma~A.1 in \cite{Borgs-Kotecky}).
\end{proofsect}

\subsubsection{Volume of large contours}
The preceding lemma asserts that, for typical configurations, the interior of large contours is not too big.
Actually, one can be a bit more precise. Namely, introducing
\begin{equation}
%\label{}
\RR^2_{\varkappa,s,L}=\bigl\{\sigma\colon 
|V(\Gamma_s^\ext(\sigma))|\ge (1-\varkappa)v_L\bigr\},
\end{equation}
we will show in the next lemma that, whenever~$\varkappa$ is sufficiently small, 
the conditional probability of~$\RR^2_{\varkappa,s,L}$ given the~$M_L$'s of 
interest is still exponentially small in~$\sqrt{v_L}$. However, unlike in 
Lemma~\ref{lemma-length} (and Lemma~\ref{lemma-nomag} below), here the constant multiplying 
$\sqrt{v_L}$ in the exponent can no longer be made arbitrarily large.

\begin{lemma}
\label{lemma-vol}
There exist constants
$c_2>0$,~$\varkappa_0>0$,~$K_0<\infty$, and~$L_0<\infty$ such that
\begin{equation}
\label{R3bd} 
P_L^{+,\beta}\bigl(\RR^2_{\varkappa,s,L}
\big|M_L=\mstar\,|\Lambda_L|-2\mstar\, v_L\bigr)\le
e^{-c_2\sqrt{v_L}}
\end{equation}
holds for all~$K\ge K_0$,~$\varkappa\in(0,\varkappa_0]$,~$L\ge L_0$, and~$s=K\log L$.
\end{lemma}

\begin{proofsect}{Proof}
Let~$\Phi_\Delta^\star$ be as defined in
\eqref{Phi*-d}. Clearly, it suffices to prove the statement for \textit{some} 
$\varkappa>0$, so 
let~$\varkappa\in(0,1)$ be such that
\begin{equation}
\label{vrkrel}
c_2=w_1\bigl[(1-\varkappa)^2-
(\Phi_\Delta^\star+2\varkappa)\bigr]>0.
\end{equation} 
(This is possible because~$\Phi_\Delta^\star<1$ for all~$\Delta<\infty$.)
Let~$L_0$ be so large that 
$\epsilon_L$ from Theorem~\ref{lowerbound} satisfies~$\epsilon_L\le\varkappa$ for 
all~$L\ge L_0$.
Let~$K_0$ be chosen to exceed the quantity~$K_0(\varkappa,\beta)$ from
Lemma~\ref{lemma-Peierls}.

Fix~$K\ge K_0$,~$L\ge L_0$, and~$s=K\log L$. 
Let now~$\sigma\in\RR^2_{\varkappa,s,L}$ and let us temporarily 
abbreviate~$\Gamma=\Gamma_s(\sigma)$ and~$\Gamma'=\Gamma_s^\ext(\sigma)$.
Let~$\frakS$ be any~$s$-skeleton such that~$\frakS\sim\Gamma$, and let~$\frakS'$ be the set of skeletons in~$\frakS$ 
corresponding to~$\Gamma'$. 
First we note that we may as well assume that, for some fixed~$B>0$ to be specified 
later
\begin{equation}
\label{P(S)specialbd}
\sum_{S\in\frakS'}\bigl|\ssP(S)\bigr|\le 
\frac B{\tau_\mini}\sqrt{v_L}.
\end{equation}
Indeed, the contribution of 
the configurations violating this bound can be directly estimated, combining 
Corollary~\ref{cor-help} with~$\alpha=\varkappa$ and \eqref{triv-W},
by~$e^{-(1-\varkappa)B\sqrt{v_L}}$.
For configurations satisfying \eqref{P(S)specialbd}, Lemma~\ref{L:skeleti} in turn implies
\begin{equation}
\label{4.23}
\bigl|V(\frakS')\bigr|\ge\bigl|V(\Gamma')
\bigr|-g_3s\sum_{S\in\frakS'}\bigl|\ssP(S)\bigr|
\ge(1-\varkappa)^2v_L,
\end{equation}
provided~$L$ is sufficiently large to ensure that~$g_3 K \frac{\log L}{\sqrt{v_L}}\frac B{\tau_\mini}\ll 1$.
As a consequence of this and the Wulff variational problem, 
$\mathscr{W}_\beta(\frakS')\ge w_1(1-\varkappa)\sqrt{v_L}$. Since~$\frakS\supset\frakS'$, we have~$\mathscr{W}_\beta(\frakS)\ge 
\mathscr{W}_\beta(\frakS')$ and thus 
for every~$\sigma\in\RR^2_{\varkappa,s,L}$ satisfying \eqref{P(S)specialbd} there is a collection 
$\frakS$ of~$s$-skeletons such that $\frakS\sim\sigma$ and 
$\mathscr{W}_\beta(\frakS)\ge w_1(1-\varkappa)\sqrt{v_L}$.
Using, once more,  
Corollary~\ref{cor-help} with~$\alpha=\varkappa$ and our choice of~$K_0$, we 
have
\begin{equation}
P_L^{+,\beta}(\RR^2_{\varkappa,s,L})\le
e^{-(1-\varkappa)^2w_1\sqrt{v_L}} +e^{-(1-\varkappa)B\sqrt{v_L}} .
\end{equation}
Letting~$B=(1-\varkappa)w_1$, the right-hand side beats the lower bound
$P_L^{+,\beta}(M_L=\mstar\,|\Lambda_L|-2\mstar\, v_L)\ge
\exp\{-w_1\sqrt{v_L}(\Phi_\Delta^\star+\varkappa)\}$ from
Theorem~\ref{lowerbound} and our choice of~$L_0$ and~$\varkappa$
by exactly~$2e^{-(c_2+\varkappa w_1)\sqrt{v_L}}$. Using the 
leeway in the exponent to absorb the extra factor of~$2$
(which may require that we further increase~$L_0$), 
the estimate \eqref{R3bd} follows.
\end{proofsect}

\subsection{Magnetization deficit due to large contours}
In this section we will provide the necessary control over the
magnetization deficit inside and outside large contours. The
relevant statements come as Lemmas~\ref{lemma-nomag}-\ref{lemma-muchmag}.

\subsubsection{Magnetization inside}
Our next claim concerns the total magnetization
inside the large contours in~$\Lambda_L$. Recalling the definition of~$\oInnt$, 
we reintroduce the event
\begin{equation}
%\label{}
\RR^3_{\varkappa,s,L}=\bigl\{\sigma\colon
M_{\oInnt}\le-\mstar\,|\oInnt|
-\varkappa^{-1} sv_L^{3/4}\bigr\}.
\end{equation}
For the probability of~$\RR^3_{\varkappa,s,L}$ we have the
following bound:

\begin{lemma}
\label{lemma-nomag} For each~$c_3>0$
there exist~$\varkappa_0>0$,~$K_0<\infty$ and~$L_0<\infty$ such
that 
\begin{equation}
\label{4.11}
P_L^{+,\beta}\bigl(\RR^3_{\varkappa,s,L}
\big|M_L=\mstar\,|\Lambda_L|-2\mstar\, v_L\bigr)\le
e^{-c_3\sqrt{v_L}}
\end{equation}
for any~$\varkappa\le\varkappa_0$,~$K\ge K_0$,~$L\ge L_0$, and
$s=K\log L$.
\end{lemma}

\begin{proofsect}{Proof}
Fix a~$c_3>0$. By Lemma~\ref{lemma-length}, there are
$\vartheta<\infty$,~$K_0<\infty$ and~$L_0<\infty$ such that
$P_L^{+,\beta}(\RR^1_{\vartheta,s,L}|
M_L=\mstar\,|\Lambda_L|-2\mstar\, v_L)\le e^{-2c_3\sqrt{v_L}}$
whenever~$s=K\log L$ and~$L\ge L_0$. Let~$\boldsymbol\Gamma
=\{\Gamma_s^\ext(\sigma)\colon\sigma\not\in\RR^1_{\vartheta,s,L}\}$. Recalling
the lower bound in Theorem~\ref{lowerbound}, it is clearly sufficient to prove 
that for some~$c_3'>0$ large enough, 
\begin{equation}
\label{4.13} P_L^{+,\beta}\bigl(\RR^3_{\varkappa,s,L}\big|
\Gamma_s^\ext(\sigma)=\Gamma\bigr)\le 2e^{-c_3'\sqrt{v_L}}
\end{equation}
holds for all~$\Gamma\in\boldsymbol\Gamma$ and all~$L$ sufficiently 
large provided~$\varkappa$ is sufficiently small and that the~$K$ in~$s=K\log L$ is sufficiently large. 
(Note that,  for \eqref{4.13} to imply \eqref{4.11},~$c_3'$ will
have to exceed~$c_3$ by a~$\beta$-dependent factor. The factor of ``$2$'' was put in for later convenience.)

Pick a~$\Gamma\in\boldsymbol\Gamma$. 
Since~$\RR^3_{\varkappa,s,L}$ depends only 
on the configuration in~$\oInnt$, 
\eqref{Int-Ext} implies
\begin{equation}
\label{4.17}
P_L^{+,\beta}\bigl(\RR^3_{\varkappa,s,L}\big|
\Gamma_s^\ext(\sigma)=\Gamma\bigr)= P_{\oInnt}^{-,\beta}
\bigl(\RR^3_{\varkappa,s,L}\bigr).
\end{equation}
In order to apply Lemma~\ref{lemma-Gauss-positive}, we need to 
compare 
$-\mstar|\oInnt|$ with the actual average
magnetization of the Ising 
model in volume~$\oInnt$ with minus boundary condition. By \eqref{int-bd} and 
\eqref{bd-int-bd}, 
we have~$|\oInnt|\le g_4^2\vartheta^{-2}s^2v_L$ and 
$|\partial\oInnt|\le g_4\vartheta^{-1}s\sqrt{v_L}$. Then Lemma~\ref{L:M-fluctuations} 
and \eqref{sbd}
imply the existence of constants~$\alpha_1=\alpha_1(\beta)<\infty$ and 
$\alpha_2=\alpha_2(\beta)>0$ 
such that
\begin{equation}
%\label{}
\bigl|\langle 
M_{\oInnt}\rangle_{\oInnt}^{-,\beta}+
\mstar|\oInnt|\bigr|\le
\alpha_1\bigl(g_4\vartheta^{-1} s\sqrt{v_L}+
(g_4^2s^2\vartheta^{-2} v_L)^2e^{-\alpha_2s}\bigr).
\end{equation}
Now, since~$s=K\log L$, for~$K$ large the right-hand side 
is less than
$2\alpha_1g_4\vartheta^{-1}s\sqrt{v_L}$. Thus, if~$L$ is so large that
the latter does not exceed $\frac12\varkappa^{-1} s v_L^{3/4}$ 
(i.e., if $4\alpha_1g_4\vartheta^{-1}s\sqrt{v_L}\le\varkappa^{-1} sv_L^{3/4}$), 
then $\sigma\in\RR^3_{\varkappa,s,L}$ and~$\Gamma_s^\ext(\sigma)=\Gamma$ imply
\begin{equation}
%\label{}
M_{\oInnt}\le\langle M_{\oInnt}\rangle_{\oInnt}^{-,\beta,s}
-\frac12\varkappa^{-1}sv_L^{3/4}.
\end{equation}
Let now~$\varkappa_0>0$ be such that 
$c_3'\le\vartheta^2(8\varkappa_0^{2}\chi g_4^2)^{-1}$, 
where~$\chi=\chi(\beta)$ is the susceptibility, and let~$\varkappa\le\varkappa_0$. 
By equation \eqref{ubp} in Lemma~\ref{lemma-Gauss-positive} 
and the fact that~$|\oInnt|\le g_4^2\vartheta^{-2}s^2v_L$, 
the right-hand side of \eqref{4.17} is bounded by~$2e^{-c_3'\sqrt{v_L}}$.
The bound \eqref{4.13} is thus proved. 
\end{proofsect}

\subsubsection{Magnetization outside}
Recall the definition of~$\oExxt$. Our 
first concern here is an upper bound 
on the total magnetization in~$\oExxt$.
Let 
$\RR^4_{\varkappa,s,L}$ be the event
\begin{equation}
%\label{}
\RR^4_{\varkappa,s,L}=\bigl\{\sigma\colon M_{\oExxt}\ge
\mstar\,|\oExxt|-2\varkappa \mstar\, v_L\bigr\}.
\end{equation}

To bound the conditional probability of this event is easy;
we will actually show that it can be included into the preceding ones for configurations contained in
$\MM_L=\{\sigma\colon M_L=\mstar|\Lambda_L|-
2\mstar\,v_L\}$.

\begin{lemma}
\label{lemma-somemag} 
For any~$\varkappa>0$ and any $K<\infty$ there exists an $L_0<\infty$ such that
\begin{equation}
\label{A4bd} 
\RR^4_{\varkappa/2,s,L}\cap\MM_L\subset
\bigl(\RR^1_{\varkappa,s,L}\cup\RR^2_{\varkappa,s,L}\cup\RR^3_{\varkappa,s,L}\bigr)\cap\MM_L
\end{equation}
for any~$L\ge L_0$ and $s=K\log L$.
\end{lemma}

\begin{proofsect}{Proof}
Let~$\varkappa$ and~$K$ be fixed.
Let us abbreviate~$\oInnt=\oInnt_{s,L}(\sigma)$ and~$\oExxt=\oExxt_{s,L}(\sigma)$ for a
configuration~$\sigma$ which we will take to be 
in~$(\RR^1_{\varkappa,s,L})^{\text{\rm c}}\cap(\RR^2_{\varkappa,s,L})^{\text{\rm
c}}\cap(\RR^3_{\varkappa,s,L})^{\text{\rm c}}
\cap\MM_L$. First, we note that
if~$\sigma\not\in\RR^1_{\varkappa,s,L}$, we can use Lemmas~\ref{lemma-corners} and~\ref{lemma-volume} to get
\begin{equation}
%\label{}
|\Lambda_L|-\bigl(|\oExxt|+|\oInnt|\bigr)\le g_4\varkappa^{-1}s
\sqrt{v_L}
\end{equation}
and hence
\begin{equation}
%\label{}
|M_L-M_{\oExxt}-M_{\oInnt}|\le g_4\varkappa^{-1}s\sqrt{v_L}.
\end{equation}
Now, since the total magnetization is held fixed, i.e.,~$\sigma\in\MM_L$, 
we have~$M_L=\mstar\,|\Lambda_L|-2\mstar\, v_L$
and by a simple calculation we get
\ifroman
\begin{multline}
\label{VGbd}
M_{\oExxt}\le M_L-M_{\oInnt}+g_4\varkappa^{-1}s\sqrt{v_L}=\\=
\mstar\,(|\Lambda_L|-|\oInnt|)-M_{\oInnt}+\mstar\,|\oInnt|-2\mstar\, v_L 
+g_4\varkappa^{-1} s\sqrt{v_L}.
\end{multline}
\else
\begin{equation}
\begin{aligned}
\label{VGbd}
M_{\oExxt}&\le M_L-M_{\oInnt}+g_4\varkappa^{-1}s\sqrt{v_L}\\&=
\mstar\,(|\Lambda_L|-|\oInnt|)-M_{\oInnt}+\mstar\,|\oInnt|-2\mstar\, v_L 
+g_4\varkappa^{-1} s\sqrt{v_L}.
\end{aligned}
\end{equation}
\fi
At the expense of another factor of~$g_4\varkappa^{-1} s\sqrt{v_L}$, 
we can replace~$|\Lambda_L|-|\oInnt|$ with~$|\oExxt|$.
Finally, since~$\sigma\not\in\RR^2_{\varkappa,s,L}\cup\RR^3_{\varkappa,s,L}$ we can use the bounds
\begin{equation}
\label{zweite} 
M_{\oInnt}\ge-\mstar\,|\oInnt|-\varkappa^{-1} sv_L^{3/4}
\end{equation}
and
\begin{equation}
\label{erste} 
|\oInnt|\le|V(\Gamma_s^\ext(\sigma))|\le (1-\varkappa)v_L
\end{equation}
in succession to arrive at
\begin{equation}
\label{dritte}
M_{\oExxt}\le\mstar\,|\oExxt|-2\mstar\,\varkappa v_L+2g_4\varkappa^{-1} s\sqrt{v_L}+
\varkappa^{-1} sv_L^{3/4}.
\end{equation}
From here we see that~$\sigma\not\in\RR^4_{\varkappa/2,s,L}$ once~$L$ 
is so large  that the remaining terms on the right-hand side are  
swamped by $-\mstar\,\varkappa v_L$.
\end{proofsect}

Our second task concerning the magnetization outside the large external 
contours is to show that~$M_{\oExxt}-\mstar|\oExxt|$ will not get substantially below the deficit 
value forced in by the condition on overall magnetization. (Note, however, that we 
have to allow for the possibility that~$\oExxt=\Lambda_L$ in which case
the exterior takes the entire deficit.)
Let~$\varkappa>0$ and consider the event
\begin{equation}
%\label{}
\RR^5_{\varkappa,s,L}=\bigl\{\sigma\colon M_{\oExxt}\le
\mstar\,|\oExxt|-2\mstar\,(1+\varkappa^{-1}) v_L\bigr\}.
\end{equation}
The probability of~$\RR^5_{\varkappa,s,L}$ is bounded as follows:

\begin{lemma}
\label{lemma-muchmag} 
For any~$c_5>0$ there exist constants~$\varkappa_0>0$ ,  
$K_0<\infty$ and~$L_0<\infty$ such that
\begin{equation}
\label{A5bd} P_L^{+,\beta}\bigl(\RR^5_{\varkappa,s,L}
\big|M_L=\mstar\,|\Lambda_L|-2\mstar\, v_L\bigr)\le
e^{-c_5\sqrt{v_L}}
\end{equation}
for all~$K\ge K_0$,~$\varkappa\le\varkappa_0$ and~$L\ge L_0$, and~$s=K\log L$.
\end{lemma}

\begin{proofsect}{Proof}
With~$\Phi_\Delta^\star$  as in \eqref{Phi*-d} and~$c_5$ fixed,  choose~$\varkappa_0$ so 
that
\begin{equation}
\label{vrkrel4}
c_5\le\frac{w_1}2\Bigl[\Delta+\frac{\Delta}{3\varkappa_0}-
\Phi_\Delta^\star\Bigr].
\end{equation}
For this~$\varkappa_0>0$, let~$L_0$ be so large that
for all~$L\ge L_0$, the finite-$L$ expression on the right-hand side of \eqref{Delta-lim} 
exceeds~$\Delta(1+\frac1{2\varkappa_0})^{-1}$ and, at the same time,
$\epsilon_L$ from Theorem~\ref{lowerbound} is bounded by~$\Delta/(6\varkappa_0)$.

First, we can restrict ourselves to the 
complement of~$\RR^1_{\vartheta,s,L}$ with~$\vartheta$ so 
small that the corresponding~$c_1$ exceeds~$2c_5$. 
Once again using Lemma~\ref{L:M-fluctuations}, we get
\begin{equation}
\label{vierte}
\bigl|\langle 
M_{\oExxt}\rangle_{\oExxt}^{+,\beta}-
\mstar|\oExxt|\bigr|\le
\alpha_1\bigl(g_4\vartheta^{-1} s\sqrt{v_L}+4L
+L^4 e^{-\alpha_2s}).
\end{equation}
Now, since~$s=K\log L$ and~$v_L\sim L^{4/3}$, for~$K$ 
sufficiently large the right-hand side does not exceed 
$8\alpha_1L$. Thus, if~$L$ is so large that
the latter does not 
exceed~$\mstar\,v_L \varkappa_0^{-1}$, it
suffices to prove the corresponding bound for the event
\begin{equation}
%\label{}
\overline{\RR}=\bigl\{\sigma\colon M_{\oExxt}\le 
\langle 
M_{\oExxt}\rangle_{\oExxt}^{+,\beta}-\mstar\,(2+\varkappa_0^{-1}) v_L\bigr\}.
\end{equation}
Clearly,~$\overline{\RR}$ depends only on the configuration in~$\oExxt$, and thus \eqref{Int-Ext} makes
the estimates in Lemma~\ref{lemma-Gauss} available. We get
\begin{equation}
\label{4.29b}
\begin{aligned}
P_L^{+,\beta}\bigl(\overline{\RR}\big|\Gamma_s^\ext(\sigma)=\Gamma\bigr)
&\le 
C\exp\Bigl\{-2\frac{(\mstar v_L)^2}{\chi|\oExxt|}
\Bigl(1+\frac{1}{2\varkappa_0}\Bigr)^2\Bigr\}\\
&\qquad\qquad\le C\exp\Bigl\{-w_1\Delta\Bigl(1+\frac{1}{2\varkappa_0}\Bigr)\sqrt{v_L}\Bigr\}.
\end{aligned}
\end{equation}
Here 
$C=C(\beta)<\infty$ is independent of~$\Gamma$ and the second inequality follows from
our assumption about~$L_0$. 
Now, using \eqref{vrkrel4} and the fact that~$\epsilon_L\le\Delta/(6\varkappa_0)$,
we derive the bound
\begin{equation}
\label{4.29c}
P_L^{+,\beta}\bigl(\overline{\RR}\big|\Gamma_s^\ext(\sigma)=\Gamma\bigr)\le 
Ce^{-w_1\sqrt{v_L}(\Phi_\Delta^\star+\epsilon_L)-2c_5\sqrt{v_L}}.
\end{equation}
The claim then follows by multiplying both sides by $P_L^{+,\beta}(\Gamma_s^\ext(\sigma)=\Gamma)$, summing over all~$\Gamma$ with the 
above properties and comparing the right-hand side with the 
lower bound in Theorem~\ref{lowerbound}.
\end{proofsect}

\subsection{Proof of Theorem~\ref{upperbound}}
\label{sec4.4}\noindent
The ultimate goal of this section is to rule out the occurrence of
intermediate contours. As a first step we derive an upper bound on the 
probability of the occurrence of contours of intermediate sizes in a contour 
ensemble constrained to not contain contours with diameters larger than 
$\varkappa\sqrt{v_L}$. The relevant statement comes as Lemma~\ref{lemma-intermediate-help}. 
Once this lemma is established, we will
give a proof of Theorem~\ref{upperbound}.

\subsubsection{A lemma for the restricted ensemble}
Recall our notation 
$P_\Lambda^{+,\beta,s'}$ for the probability measure
in volume~$\Lambda\subset\Lambda_L$ 
conditioned on the event that the contour diameters do
not exceed~$s'$. We will show that
the occurrence of intermediate contours is improbable in~$P_\Lambda^{+,\beta,s'}$ 
with~$s'=\varkappa\sqrt{v_L}$ and magnetization restricted to
``reasonable'' values. 
For any~$\Lambda\subset\Lambda_L$ and any~$s>0$ and~$\varkappa>0$, let
\begin{equation}
%\label{}
\AA^{\text{\rm c}}_{\varkappa,s,\Lambda}=\bigl\{\sigma\colon 
\text{\rm there exists }\gamma\text{ \rm 
in }\Lambda\text{ such that } 
s\le\diam\gamma\le 
\varkappa\sqrt{v_L}\bigr\}.
\end{equation}
Then we have the following estimates:

\begin{lemma}
\label{lemma-intermediate-help} 
For any~$c_6>0$,~$\varphi_0>1$,
and~$\vartheta>1$, there exist~$\varkappa_0\in(0,1)$,~$K_0<\infty$, 
and~$L_0<\infty$, such that for~$s=K\log L$, 
all~$\varkappa\in(0,\varkappa_0]$,~$K\ge K_0$,~$L\ge L_0$, 
all~$\Lambda\subset\Lambda_L$ satisfying the bounds
\begin{equation}
\label{Vrestr}
|\Lambda|\ge\vartheta^{-
1}L^2\quad\text{and}\quad
|\partial \Lambda|\le\vartheta L,
\end{equation}
and all~$\varphi\in[\varkappa_0,\varphi_0]$ 
that make~$\mstar\,|\Lambda|-2\varphi\mstar\,  v_L$ an
allowed value of~$M_\Lambda$, we have
\begin{equation}
\label{Lc6}
P_\Lambda^{+,\beta,\varkappa\sqrt{v_L}}\bigr(\AA^{\textup{c}}_{\varkappa,s,\Lambda}\big|
M_\Lambda=\mstar\,|\Lambda|-2\varphi\mstar\,  v_L\bigr)\le L^{-c_6}.
\end{equation}
\end{lemma}

\begin{proofsect}{Proof}
Notice that the event~$\AA^{\textup{c}}_{\varkappa,s,\Lambda}$ is monotone in~$s=K\log L$
and thus it is sufficient to prove the claim for only a fixed~$K$ (chosen suitably large).
Let $\varkappa_0\in(0,1)$ be fixed and let $\varkappa\in(0,\varkappa_0]$.
(At the very end of the proof, we will have to assume that $\varkappa_0$ is sufficiently
small, see \eqref{4.54}.)
Fix a set~$\Lambda\subset\Z^2$ satisfying \eqref{Vrestr} and let 
\begin{equation}
\label{MMA}
\MM_\Lambda(\varphi)=\bigl\{\sigma\colon M_\Lambda=\mstar\,|\Lambda|-2\varphi\mstar\,  v_L\bigr\}.
\end{equation}  Let us define
\begin{equation}
\label{deltaV}
\delta_\Lambda=\langle 
M_\Lambda\rangle_\Lambda^{+,\beta,s}-\mstar|\Lambda|
\end{equation}
and note that, on~$\MM_\Lambda(\varphi)$, we 
have~$M_\Lambda=\langle M_\Lambda\rangle_\Lambda^{+,\beta,s}
-\delta_\Lambda-2\varphi\mstar v_L$.

The proof of \eqref{Lc6} will be performed by writing the conditional
probability as a quotient of two probabilities with unconstrained
contour sizes and estimating separately the numerator and the
denominator. Let
\begin{equation}
\label{4.28}
\EE=\bigl\{\sigma\colon\forall\gamma\in\Gamma_s(\sigma),
\,\diam\gamma\le\varkappa\sqrt{v_L}\bigr\}
\end{equation}
and, using the shorthand~$\AA=\AA_{\varkappa,s,\Lambda}$, write
\begin{equation}
P_\Lambda^{+,\beta,\varkappa\sqrt{v_L}}\bigr(\AA^{\text{\rm c}}\big|
\MM_\Lambda(\varphi)\bigr)=\frac{P_\Lambda^{+,\beta}(\AA^{\text{\rm c}}
\cap\MM_\Lambda(\varphi)\cap \EE)}
{P_\Lambda^{+,\beta}(\MM_\Lambda(\varphi)\cap\EE)}.
\end{equation}
As to the bound on the denominator, we restrict the contour
sizes in~$\Lambda$ to~$s=K\log L$ as in \eqref{MLlbd} and apply Lemmas~\ref{lemma-Gauss} 
and~\ref{lemma-nolog} with the result
\begin{equation}
\label{4.27a} 
P_\Lambda^{+,\beta}(\MM_\Lambda(\varphi)\cap\EE)\ge \frac
{C_1}{L^2}\exp\Bigl\{-2\frac{(\mstar\,v_L)^2}{\chi|\Lambda|}\varphi^2
-2\frac{\mstar\,\varphi\,v_L}{\chi|\Lambda|}\delta_\Lambda\Bigr\},
\end{equation}
where~$C_1=C_1(\beta,\vartheta, \varphi_0)>0$. 
Here, we note that two distinct terms were incorporated into the constant~$C_1$:
First, a term proportional to $\delta_\Lambda^2$ since, by Lemma~\ref{L:M-fluctuations} and 
\eqref{Vrestr},~$|\delta_\Lambda|\le 
2\alpha_1\vartheta L$ once~$K$ is sufficiently large 
and thus~$|\delta_\Lambda|^2/|\Lambda|$ is bounded by a constant independent of~$L$. 
Second, a term that comes from the bound \eqref{2.45} yielding
$|\Omega_\Lambda^s(\varphi v_L +\frac{\delta_\Lambda}{2\mstar})|\le 
C_2 \max\{K \frac{\log L}{L^{1/3}},1\}$
with some~$C_2=C_2(\beta,\vartheta, \varphi_0)<\infty$. 
(Notice that, to get a constant~$C_1$ independent 
of~$L$, we have to choose~$L_0$ after a choice of~$K$ is done.)
Although the second term on the right-hand side of \eqref{4.27a} 
is negligible compared to the first one, 
its exact form will be needed to cancel an inconvenient contribution of the 
complement of intermediate contours. 

In order to estimate the numerator, let
$\boldsymbol\Gamma=\{\Gamma_s(\sigma)\colon
\sigma\in\EE,\,\Gamma_s(\sigma)\ne\emptyset\}$ be the set of all collections of 
$s$-large contours that can possibly contribute to~$\EE$. (We also demand that 
$\Gamma_s(\sigma)\ne\emptyset$, because on~$\AA^{\text{\rm c}}$ there will be 
at least one~$s$-large contour.) Then we have
\begin{equation}
\label{4.29} 
P_\Lambda^{+,\beta}\bigl(\AA^{\text{\rm c}}\cap\MM_\Lambda(\varphi)\cap\EE\bigr)
\le\sum_{\Gamma\in\boldsymbol\Gamma}
P_\Lambda^{+,\beta}\bigl(\MM_\Lambda(\varphi)\big|\Gamma_s(\sigma)=\Gamma\bigr)
P_\Lambda^{+,\beta}\bigl(\Gamma_s(\sigma)=\Gamma\bigr).
\end{equation}
Our strategy is to derive a bound on 
$P_\Lambda^{+,\beta}(\MM_\Lambda(\varphi)|\Gamma_s(\sigma)=\Gamma)$ 
which is uniform in~$\Gamma\in\boldsymbol\Gamma$ and to estimate 
$P_\Lambda^{+,\beta}(\Gamma_s(\sigma)=\Gamma)$ using the
skeleton upper bound.

Let~$\Gamma\in\boldsymbol\Gamma$ and let~$\frakS$ be an~$s$-skeleton such 
that~$\frakS\sim\Gamma$. We claim that, for some 
$C'=C'(\beta,\vartheta)<\infty$ and some~$\eta_0=\eta_0(\beta,\vartheta)<\infty$, 
independent of~$\Gamma$,~$\frakS$,~$\varkappa_0$ and~$L$,
\begin{equation}
\label{relbd} 
\frac
{P_\Lambda^{+,\beta}(\MM_\Lambda(\varphi)|\Gamma_s(\sigma)=\Gamma)}
{P_\Lambda^{+,\beta}(\MM_\Lambda(\varphi)\cap\EE)}\le 
C'L^2 e^{\eta_0\sqrt
\varkappa_0
\mathscr{W}_\beta(\frakS)}
\end{equation}
holds true.
Indeed, let~$\Gamma'$ be the abbreviation for the set 
of external contours in~$\Gamma$ and let~$\frakS'$ be the set of skeletons in 
$\frakS$ corresponding to~$\Gamma'$. Recall the definition of~$\Innt$ and 
$\oInnt$ and note that~$\V(\Gamma')=\Innt$ and~$\mathscr{W}_\beta(\frakS)\ge\mathscr{W}_\beta(\frakS')$, since~$\frakS\supset\frakS'$. Also 
note that, by \eqref{triv-kappa} and \eqref{triv-W} and the fact that 
$\diam\gamma\le \varkappa\sqrt{v_L}$ for all~$\gamma\in\Gamma'$, we 
have
\begin{equation}
\label{VGammabd} 
|\Innt| \le g_2\varkappa\sqrt{v_L}\sum_{S\in\frakS'}\bigl|\ssP(S)\bigr|\le
g_2\varkappa_0\tau_\mini^{-1}\sqrt{v_L}\,\mathscr{W}_\beta(\frakS).
\end{equation}
This bound tells us that we might as well assume that~$|\Innt|\le 
\sqrt{\varkappa_0} v_L$. Indeed, in the opposite case, 
the bound \eqref{relbd} would directly follow by noting 
that \eqref{4.27a} implies 
$P_L^{+,\beta}(\MM_\Lambda(\varphi)\cap\EE)\ge 
C_1L^{-2}e^{-
\eta_1\sqrt{\varkappa_0}
\mathscr{W}_\beta(\frakS)}$ with~$\eta_1$ given by
\begin{equation}
%\label{}
\eta_1=2g_2
\Bigl[\frac{(\mstar\,\varphi)^2}
{\chi\tau_\mini}\frac{v_L^{3/2}}{|\Lambda|}+
\frac{\mstar\,\varphi}{\chi\tau_\mini}
\frac{\delta_\Lambda\sqrt{v_L}}{|\Lambda|}\Bigr].
\end{equation}
Notice that~$\eta_1$ is bounded uniformly in~$L$ and~$\Lambda$
by \eqref{Vrestr} and 
the facts that~$\Delta<\infty$ and $\delta_\Lambda\le2\alpha_1\vartheta L$. 
A similar bound, using  \eqref{triv-s} instead of \eqref{triv-kappa}, shows that also
$|\partial\Innt|\le s\sqrt{v_L}/\sqrt{\varkappa_0}$.
Indeed, if the opposite is true, then \twoeqref{triv-s}{triv-W} imply that $\sqrt{\varkappa_0}\mathscr{W}_\beta(\frakS)\ge\tau_\mini g_1^{-1}\sqrt{v_L}$ and we
can proceed as before.

Thus, let us assume that 
$|\Innt|\le\sqrt{\varkappa_0} v_L$ and~$|\partial\Innt|\le s\sqrt{v_L}/\sqrt{\varkappa_0}$ hold 
true. In order for $\MM_\Lambda(\varphi)$ to occur, the total
magnetization in~$\Lambda$ should deviate from~$\mstar\,|\Lambda|$ by~$-2\varphi\mstar\,  v_L$, 
while the volume~$\Innt$ can help the bulk only by at
most~$-|\Innt|$. More precisely,~$M_{\oExxt}$ is forced to deviate from its mean 
value~$\langle M_{\oExxt}
\rangle_{\oExxt}^{+,\beta,s}$ by 
at least~$-2\mstar u$ 
(and by not more than $-2\mstar u-2|\Innt|$) 
where~$u$ is defined by
\begin{equation}
%\label{}
-2\mstar u=-2\varphi\mstar\,  v_L-
\delta_{\oExxt}+2|\Innt|,
\end{equation} 
with~$\delta_{\oExxt}$ as in \eqref{deltaV}.
By the estimates~$|\Innt|\le \sqrt{\varkappa_0} v_L$,
$|\oExxt|\ge\frac12\vartheta^{-1}L^2$, 
$|\partial\oExxt|\le 2\vartheta  L$, and~$u\le C_3 L^{4/3}\ll L^2/\log L$, 
with~$C_3=C_3(\beta,\vartheta, \varphi_0)$ (all these bounds  hold for~$L$ 
sufficiently large---in particular, to ensure that~$K\sqrt{v_L}\log L\le \vartheta L$),
we now have, once more,  Lemma~\ref{lemma-Gauss} at our disposal. Thus,
\begin{equation}
\label{predchozi}
P_\Lambda^{+,\beta}\bigl(\MM_\Lambda(\varphi)\big|\Gamma_s(\sigma)=\Gamma\bigr)\le
C_4\exp\Bigl\{-2\frac{(\mstar\, v_L)^2}{\chi|\Lambda|}\varphi^2
-2\frac{\mstar\,\varphi v_L}{\chi|\Lambda|}\bigl(\delta_{\oExxt}-2|\Innt|\bigr)\Bigr\},
\end{equation}
where~$C_4=C_4(\beta,\vartheta,\varphi_0)<\infty$.
Similarly as in \eqref{4.27a}, the constant~$C_4$ incorporates also the error term
$\Omega_{\oExxt}^s(u)$.
To compare the right-hand side of 
\eqref{predchozi} and \eqref{4.27a}, we invoke 
the second part of Lemma~\ref{L:M-fluctuations} to 
note that, for~$K$ sufficiently large and some 
$\alpha_1=\alpha_1(\beta)<\infty$,
\begin{equation}
%\label{}
\delta_\Lambda-\delta_{\oExxt}
\le \alpha_1|\Lambda\setminus\oExxt|.
\end{equation}
Using \eqref{VGammabd} 
again,~$|\Innt|$ is bounded by a constant
times~$\varkappa_0\mathscr{W}_\beta(\frakS)\sqrt{v_L}$ 
and the same holds for~$|\Lambda\setminus\oExxt|$. Therefore, there is a constant 
$\eta_2=\eta_2(\beta,\vartheta)<\infty$,
independent of~$\varkappa_0$, such that
\begin{equation}
\label{2ndbd}
2\frac{\mstar\,\varphi v_L}{\chi|\Lambda|}
\bigl(\delta_\Lambda-\delta_{\oExxt}+2|\Innt|\bigr)\le 
\eta_2\varkappa_0\mathscr{W}_\beta(\frakS),
\end{equation}
holds true for all 
$\Gamma\in\boldsymbol\Gamma$ and their associated skeletons~$\frakS$. By 
combining this with \eqref{predchozi} and \eqref{4.27a},  the bound \eqref{relbd} is 
established with~$\eta_0=\max\{\eta_1,\eta_2\}$,
which we remind is independent of~$\varkappa_0$.

With \eqref{relbd}, the proof is easily concluded. Indeed, a straightforward 
application of the skeleton bound to the second
term on the right-hand side of \eqref{4.29} then shows that
\begin{equation}
\label{4.54}
P_\Lambda^{+,\beta,\varkappa\sqrt{v_L}}\bigl(\AA^{\text{\rm c}}\big|\MM_\Lambda(\varphi)\bigr)\le
\sum_{\frakS\ne\emptyset}C'L^2e^{-(1-\eta_0\sqrt{\varkappa_0})\mathscr{W}_\beta(\frakS)}.
\end{equation}
Now, choosing~$\varkappa_0$ sufficiently small, we have
$1-\eta_0\sqrt{\varkappa_0}>2/3$. Then we can extract the term
$C'e^{-\frac13\mathscr{W}_\beta(\frakS)}$ 
which, choosing the~$K$ in
$s=K\log L$ sufficiently large, can be made less than~$L^{-2-c_6}$, for any~$c_6$
initially prescribed. Invoking Lemma~\ref{lemma-Peierls}, the remaining sum is then estimated by one.
\end{proofsect}

\subsubsection{Absence of intermediate contours}
Lemmas~\ref{lemma-length} and~\ref{lemma-vol}-\ref{lemma-intermediate-help} 
finally put us in the position to rule out the intermediate contours altogether.

\begin{proofsect}{Proof of Theorem~\ref{upperbound}}
Recall that our goal is to prove \eqref{intermediate}, i.e., $P_L^{+,\beta}(\AA^{\text{c}}|\MM_L)\le L^{-c_0}$.
Pick any~$c_0>0$ and~$\varkappa_0<1$. Let~$K_0$ and~$L_0$
be chosen so that Lemmas~\ref{lemma-length}, \ref{lemma-vol}, 
\ref{lemma-nomag}, and \ref{lemma-muchmag} hold with \emph{some}
$c_1,c_2, c_3, c_5>0$ for all~$\varkappa\le2\varkappa_0$,~$K\ge K_0$ and~$L\ge L_0$. 
We also  assume that~$L_0$ is chosen so that Lemma~\ref{lemma-somemag} is valid 
for~$\varkappa=2\varkappa_0$.
We wish to restrict attention to configuration outside the sets $\RR^1_{\varkappa_0,s,L}$, 
$\RR^4_{\varkappa_0,s,L}$ and $\RR^5_{\varkappa_0,s,L}$, 
but since $\RR^4_{\varkappa_0,s,L}$ is essentially included in
$\RR^2_{\varkappa_0,s,L}$ and $\RR^3_{\varkappa_0,s,L}$, we might as well focus on the event
$\RR^{\text{c}}$, where
$\RR=\bigcup_{\ell=1}^5\RR^\ell_{\varkappa_0,s,L}$.
Fix any $\varkappa\le\varkappa_0$, let $s=K\log L$ and let us introduce 
the shorthand~$\AA=\AA_{\varkappa,s,L}$.
Appealing to the aforementioned lemmas, our goal will be achieved 
if we establish the bound~$P_L^{+,\beta}(\AA^{\text{c}}\cap\RR^{\text{c}}|\MM_L)\le L^{-2c_0}$.

Let us abbreviate~$q=\varkappa\sqrt{v_L}$ and 
let~$\boldsymbol\Gamma=\{\Gamma_q^\ext(\sigma)\colon \sigma\in\RR^{\text{c}}\}$ be the set of all
collections of external contours that can possibly arise from~$\RR^{\text{c}}$. 
Fix~$\Gamma\in\boldsymbol\Gamma$ and recall our notation~$\oExxt$ for
the exterior component of~$\Lambda_L$ induced by the contours in~$\Gamma$. 
To prove \eqref{intermediate}, it suffices to show that, for all $\Gamma\in\boldsymbol\Gamma$,
\begin{equation}
\label{frbd} P_L^{+,\beta}\bigl(\AA^{\text{c}}\cap\RR^{\text{c}}\cap\MM_L
\big|\Gamma_q^\ext(\sigma)=\Gamma\bigr) \le
L^{-2c_0}P_L^{+,\beta}\bigl(\MM_L
\big|\Gamma_q^\ext(\sigma)=\Gamma\bigr).
\end{equation}
Indeed, multiplying \eqref{frbd} by  $P_L^{+,\beta}(\Gamma_q^\ext(\sigma)=\Gamma)$ and
summing over all~$\Gamma\in\boldsymbol\Gamma$, we derive that
\begin{equation}
%\label{}
P_L^{+,\beta}\bigl(\AA^{\text{c}}\cap\RR^{\text{c}}\cap\MM_L\bigr)\le
L^{-2c_0}P_L^{+,\beta}(\MM_L).
\end{equation}
Thence, $P_L^{+,\beta}(\AA^{\text{c}}\cap\RR^{\text{c}}|\MM_L)\le L^{-2c_0}$ which, 
in light of the bound
$P_L^{+,\beta}(\RR|\MM_L)\le 4e^{-c\sqrt{v_L}}$ where $c=\min\{c_1,c_2,c_3,c_5\}$, implies
\eqref{intermediate} once~$L$ is sufficiently large.

It remains to prove \eqref{frbd} for all $\Gamma\in\boldsymbol\Gamma$.
Let~$\varphi\ge0$ be such that~$\mstar\,|\oExxt|-2\varphi\mstar\,v_L$ is an 
allowed value of~$M_{\oExxt}$ and consider the corresponding event~$\MM_{\oExxt}(\varphi)$ 
(cf. \eqref{MMA}).
Note that, by the restriction to the complements of~$\RR^4_{\varkappa_0,s,L}$ 
and~$\RR^5_{\varkappa_0,s,L}$, we only need to consider $\varphi\in[\varkappa_0,1+\varkappa_0^{-1}]$. 
We claim that, for all such allowed values of~$\varphi$, we have
\begin{equation}
\label{stepik}
P_L^{+,\beta}\bigl(\AA^{\text{c}}\big|\{\Gamma_q^\ext(\sigma)=\Gamma\}
\cap\MM_L\cap\MM_{\oExxt}(\varphi)\bigr)=
P_{\oExxt}^{+,\beta,q}\bigl(\AA^{\text{c}}
\big|\MM_{\oExxt}(\varphi)\bigr).
\end{equation}
Indeed, given 
that~$\Gamma_q^\ext(\sigma)=\Gamma$, the event~$\AA$
depends only on the configurations in~$\oExxt$.
Moreover,~$\MM_L\cap\MM_{\oExxt}(\varphi)$ can 
be written as an intersection of~$\MM_{\oExxt}(\varphi)$, which also 
depend only on~$\sigma$ in~$\oExxt$, and the event
$\{\sigma\colon M_{\Lambda_L\setminus\oExxt}=
\mstar\,(|\Lambda_L|-|\oExxt|)-2\mstar\, (1-\varphi) v_L\}$,
which depends only on the configuration in~$\oInnt$.
Thus, \eqref{stepik} follows from \eqref{Int-Ext} and some elementary manipulations.

By the restriction to the complement of~$\RR^1_{\varkappa_0,s,L}$, we have
$|\oExxt|\ge L^2/2$ and~$|\partial\oExxt|\le 8L$ for all 
$\Gamma\in\boldsymbol\Gamma$. 
Choosing now~$c_6=2c_0$ and then $K_0$ and~$L_0$ (if necessary, even bigger than before)
so that  Lemma~\ref{lemma-intermediate-help} can be applied, 
the right-hand side of \eqref{stepik} can be bounded by 
$L^{-c_6}= L^{-2c_0}$ uniformly in~$\Gamma\in\boldsymbol\Gamma$, provided~$\varkappa$ 
is sufficiently small and~$L\ge L_0$. Using \eqref{stepik}, we thus have
% ROMAN'S VERSION
\ifroman
\begin{multline}
P_L^{+,\beta}\bigl(\AA^{\text{c}}\cap\RR^{\text{c}}\cap\MM_L\cap\MM_{\oExxt}(\varphi)
\big|\Gamma_q(\sigma)=\Gamma\bigr)\le\\ 
\le
P_L^{+,\beta}\bigl(\AA^{\text{c}}\big|\{\Gamma_q^\ext(\sigma)=\Gamma\}
\cap\MM_L\cap\MM_{\oExxt}(\varphi)\bigr) P_L^{+,\beta}\bigl(\MM_L\cap\MM_{\oExxt}(\varphi)
\big|\Gamma_q(\sigma)=\Gamma\bigr)\le\\
\le L^{-2c_0}P_L^{+,\beta}\bigl(\MM_L\cap\MM_{\oExxt}(\varphi)
\big|\Gamma_q(\sigma)=\Gamma\bigr),
\end{multline}
\else
% EVERYBODY ELSE'S VERSION
\begin{equation}
\begin{aligned}
P_L^{+,\beta}\bigl(\AA^{\text{c}}\cap\RR^{\text{c}}\cap\MM_L&\cap\MM_{\oExxt}(\varphi)
\big|\Gamma_q(\sigma)=\Gamma\bigr)\\ 
&\le
P_L^{+,\beta}\bigl(\AA^{\text{c}}\big|\{\Gamma_q^\ext(\sigma)=\Gamma\}
\cap\MM_L\cap\MM_{\oExxt}(\varphi)\bigr)\\
&\qquad\qquad\qquad\qquad\times P_L^{+,\beta}\bigl(\MM_L\cap\MM_{\oExxt}(\varphi)
\big|\Gamma_q(\sigma)=\Gamma\bigr)\\
&\le L^{-2c_0}P_L^{+,\beta}\bigl(\MM_L\cap\MM_{\oExxt}(\varphi)
\big|\Gamma_q(\sigma)=\Gamma\bigr),
\end{aligned}
\end{equation}
\fi
%%%%
for \emph{all}~$\varphi$ for which~$\mstar\,|\oExxt|-2\varphi\mstar\,  v_L$
is an allowed value of~$M_{\oExxt}$. (In the cases when $\varphi\not\in[\varkappa_0,1+\varkappa_0^{-1}]$ 
we have $\RR^{\text{c}}\cap\MM_{\oExxt}(\varphi)=\emptyset$ and the left-hand side vanishes.) 
This implies \eqref{frbd} by summing over all allowed values of~$\varphi$. 
\end{proofsect}

\section{Proof of main results}
\label{sec5}\smallskip\noindent
Having established the absence of intermediate-size contours, we
are now in the position to prove our main results. 

\begin{proofsect}{Proof of Theorem~\ref{main-result}}
Fix a $\zeta>0$  and recall our notation~$\MM_L=\{\sigma\colon M_L=\mstar|\Lambda_L|-
2\mstar\,v_L\}$. Our goal is to estimate the conditional probability 
$P_L^{+,\beta}(\AA_{\varkappa,s,L}^{\text{c}}\cup\BB_{\epsilon,s,L}^{\text{c}}|\MM_L)$
by $L^{-\zeta}$. Let $c_0>\zeta$ and note that, by Theorem~\ref{upperbound}, we have
\begin{equation}
%\label{}
P_L^{+,\beta}(\AA_{\varkappa,s,L}^{\text{c}}|\MM_L)\le L^{-c_0},
\end{equation}
provided~$\varkappa$ is sufficiently small and~$L$ 
sufficiently large.
This means we can restrict our attention to the event 
$\BB_{\epsilon,s,L}^{\text{c}}\setminus\AA_{\varkappa,s,L}^{\text{c}}$. 
Furthermore, we can use Lemmas~\ref{lemma-length},~\ref{lemma-vol},~\ref{lemma-nomag}, 
and~\ref{lemma-somemag} 
to exclude the events~$\RR^1_{\vartheta,s,L}$,~$\RR^2_{\vartheta,s,L}$,
$\RR^3_{\vartheta,s,L}$, and~$\RR^4_{\vartheta,s,L}$, provided~$\vartheta$ is 
sufficiently small. We therefore introduce the event 
$\EE_{\epsilon,\varkappa,\vartheta}$ defined 
by
\begin{equation}
%\label{}
\EE_{\epsilon,\varkappa,\vartheta}=
\BB_{\epsilon,s,L}^{\text{c}}\setminus(\AA_{\varkappa,s,L}^{\text{c}}
\cup 
\RR^1_{\vartheta,s,L}\cup \RR^2_{\vartheta,s,L}\cup 
\RR^3_{\vartheta,s,L}\cup \RR^4_{\vartheta,s,L}),
\end{equation}
where we have suppressed~$s=K\log L$ and~$L$ from the notation.

On the basis of the aforementioned Lemmas, the proof of Theorem~\ref{main-result} 
will follow if we can establish that for each~$\varkappa>0$ and each 
$\epsilon>0$ there are~$K_0<\infty$,~$\vartheta>0$ and~$c_7>0$ such that
\begin{equation}
\label{expbd}
P_L^{+,\beta}(\EE_{\epsilon,\varkappa,\vartheta}|\MM_L)\le
e^{-c_7\sqrt{v_L}}
\end{equation}
whenever~$L$ is sufficiently large.
The proof of \eqref{expbd} will be performed by conditioning on the set of $s$-large 
exterior contours and applying separately the Gaussian 
estimates and the skeleton upper bound. The argument will be split into several cases, 
depending on which of the bounds \twoeqref{hausbd}{magbd} 
constituting the event~$\BB_{\epsilon,s,L}$ fail to hold.

Let~us write $\EE_{\epsilon,\varkappa,\vartheta}$ as the disjoint union
$\EE_{\epsilon,\varkappa,\vartheta}^1
\cup\EE_{\epsilon,\varkappa,\vartheta}^2$, where 
$\EE_{\epsilon,\varkappa,\vartheta}^1$ is the set of all configurations on which 
one of \eqref{hausbd} or \eqref{volbd} fail and where 
$\EE_{\epsilon,\varkappa,\vartheta}^2
=\EE_{\epsilon,\varkappa,\vartheta}
\setminus\EE_{\epsilon,\varkappa,\vartheta}^1$.
Let $\boldsymbol\Gamma=\{\Gamma_s^\ext(\sigma)\colon 
\sigma\in\EE_{\epsilon,\varkappa,\vartheta}\}$ be the set of all collections of exterior contours 
allowed by~$\EE_{\epsilon,\varkappa,\vartheta}$. (Here $s=K\log L$.) 
Since $\Gamma_s(\sigma)$ is non-empty for all~$\sigma$ contributing to 
$\BB_{\epsilon,s,L}^{\text{c}}$, we have~$\Gamma\ne\emptyset$ for all $\Gamma\in\boldsymbol\Gamma$. 
Let
\begin{equation}
%\label{}
\lambda_\Gamma=v_L^{-1}|V(\Gamma)|.
\end{equation}
To apply the Gaussian estimate, we need the following \emph{upper} bound on the magnetization in~$\oExxt$.

\begin{lemma}
\label{lemma5.1}
Let $\epsilon>0$, $\varkappa>0$ and $\vartheta>0$ and let the~$K$ in $s=K\log L$ be sufficiently large. 
Then there exists a sequence $(\kappa_L)$ with $\lim_{L\to\infty}\kappa_L=0$ 
such that for both $i=1,2$, all $\Gamma\in\boldsymbol\Gamma$ and all $\sigma\in\MM_L\cap
\EE_{\epsilon,\varkappa,\vartheta}^i\cap\{\Gamma_s^\ext(\sigma)=\Gamma\}$, the magnetization
$M_{\oExxt}=M_{\oExxt_{s,L}(\sigma)}(\sigma)$ obeys the bound
\begin{equation}
\label{MExtbd}
M_{\oExxt}\le\langle M_{\oExxt}\rangle_{\oExxt}^{+,\beta,s}
-2\mstar\,v_L(1-\lambda_\Gamma+\epsilon_i-\kappa_L).
\end{equation}
Here $\epsilon_1=0$  and $\epsilon_2=\epsilon/(2\mstar)$.
\end{lemma}

\begin{proofsect}{Proof}
Recall the exact definition of $\oExxt$. The proof is similar in spirit to the reasoning 
\twoeqref{erste}{dritte}.
First we will address
the case of configurations in~$\EE_{\epsilon,\varkappa,\vartheta}^1$. Using the equality
$M_L=\mstar|\Lambda_L|-2\mstar\,v_L$ and our restriction to the complement of $\RR^1_{\vartheta,s,L}$,
we have
\begin{equation}
%\label{}
M_L\le \mstar|\oExxt|+\mstar|V(\Gamma)|-2\mstar v_L+g_4\vartheta^{-1}s\sqrt{v_L},
\end{equation}
where $g_4\vartheta^{-1}s\sqrt{v_L}$ bounds the volume of~$\Exxt\setminus\oExxt$
according to Lemma~\ref{lemma-corners}. 
Next, in view of the restriction to $(\RR^3_{\vartheta,s,L})^{\text{\rm c}}$, we have
\begin{equation}
\label{vG1}
M_{\V(\Gamma)}\ge -\mstar|V(\Gamma)|-\vartheta^{-1}sv_L^{3/4}-g_4\vartheta^{-1}s\sqrt{v_L}.
\end{equation}
Finally, since $M_{\oExxt}\le M_L-M_{\V(\Gamma)}+g_4\vartheta^{-1}s\sqrt{v_L}$ and since \eqref{vierte}
implies that $\mstar|\oExxt|-\langle M_{\oExxt}\rangle_{\oExxt}^{+,\beta,s}$ 
can be bounded by~$8\alpha_1L$ once~$K$ is sufficiently large, we have \eqref{MExtbd} with~$\kappa_L$
given by
\begin{equation}
2\mstar\kappa_L=
\vartheta^{-1}sv_L^{-1/4}+3g_4\vartheta^{-1}sv_L^{-1/2}+8\alpha_1L v_L^{-1}.
\end{equation}
Since $v_L\sim L^{4/3}$, we have $\lim_{L\to\infty}\kappa_L=0$ as claimed.

Next we will attend to the  case of configurations from $\EE_{\epsilon,\varkappa,\vartheta}^2$, 
for which the bound \eqref{magbd} must fail. 
Since $\EE_{\epsilon,\varkappa,\vartheta}^2$ is still a subset  of
$(\RR^3_{\vartheta,s,L})^{\text{\rm c}}$, 
we still have the bound  \eqref{vG1} at our disposal implying that $M_{\V(\Gamma)}\ge
-\mstar|V(\Gamma)|-\epsilon v_L$ once~$L$ is sufficiently large.  
However, this means that the only way \eqref{magbd} can fail is
that, in fact, the lower bound
\begin{equation}
\label{vG2}
M_{\V(\Gamma)}\ge -\mstar|V(\Gamma)|+\epsilon v_L
\end{equation}
holds. Substituting this 
stronger bound in the above derivation 
in the place of \eqref{vG1}, the desired estimate follows.
\end{proofsect}

With Lemma~\ref{lemma5.1} in the hand, we are ready to start proving the bound \eqref{expbd}. 
We begin with the Gaussian estimate.
By the restriction to 
the complement of~$\RR^2_{\vartheta,s,L}$, 
we have the bound~$\lambda_\Gamma\le 1-\vartheta$ and thus 
$1-\lambda_\Gamma+\epsilon_i-\kappa_L\ge0$ once~$L$ is sufficiently large.
Moreover, since we also discarded~$\RR^1_{\vartheta,s,L}$,
 Lemma~\ref{lemma-Gauss} for $A=\oExxt$ applies. 
Combining this with the observation \eqref{Int-Ext} and the bound \eqref{MExtbd}, 
there exists a constant $C<\infty$ such that
\begin{equation}
\label{vodhad}
P_L^{+,\beta}\bigl(\MM_L\cap\EE_{\epsilon,\varkappa,\vartheta}^i\big|
\Gamma_s^\ext(\sigma)=\Gamma\bigr)
\le C\exp\biggl\{-
2\frac{(\mstar\,v_L)^2}
{\chi|\Lambda_L|}(1-\lambda_\Gamma+\epsilon_i-\kappa_L)^2\biggr\}
\end{equation}
holds for all $\Gamma\in\boldsymbol\Gamma$.
Next we will estimate the probability that $\Gamma_s^\ext(\sigma)=\Gamma$.
Let~$\frakS$ be a collection of skeletons corresponding to 
$\Gamma$. The skeleton upper bound in Lemma~\ref{lemma-skeletonUB} along with the 
estimates featured in Lemma~\ref{lemma-Peierls} then yields
\begin{equation}
\label{druhej}
P_L^{+,\beta}\bigl(\Gamma_s^\ext(\sigma)=\Gamma\bigr)
\le \sum_{\frakS'\supseteq\frakS}e^{-\mathscr{W}_\beta(\frakS')}\le
C'e^{-\mathscr{W}_\beta(\frakS)},
\end{equation}
where~$C'<\infty$ and 
where~$\frakS'$ corresponds to the skeleton of 
a full set $\Gamma_s(\sigma)$ with $\Gamma_s^\ext(\sigma)=\Gamma$.

To estimate the probability of~$\MM_L\cap 
\EE_{\epsilon,\varkappa,\vartheta}^i\cap\{\Gamma_s^\ext(\sigma)=\Gamma\}$, 
we will write~$\boldsymbol\Gamma$ as the union of two disjoint sets, 
$\boldsymbol\Gamma=\boldsymbol\Gamma_1\cup\boldsymbol\Gamma_2$.
Here
\begin{equation}
\label{G1}
\boldsymbol\Gamma_1=
\bigl\{\Gamma\in\boldsymbol\Gamma\colon\exists\frakS\sim\Gamma,\,
\mathscr{W}_\beta(\frakS)\le 
w_1\sqrt{\lambda_\Gamma v_L}(1+\epsilon c^{-2})\bigr\},
\end{equation}
where~$c$ is the constant from Lemma~\ref{lemma-onecontour}, and
$\boldsymbol\Gamma_2=\boldsymbol\Gamma\setminus\boldsymbol\Gamma_1$.
First we will 
study the cases when $\Gamma\in\boldsymbol\Gamma_1$. 
By the restriction to the 
event~$\AA_{\varkappa,s,L}$, we know that 
$\diam\gamma\ge\varkappa\sqrt{v_L}$ for all~$\gamma\in\Gamma$.
Using that 
$\lambda_\Gamma\le1-\vartheta$---recall that we are in the complement of $\RR^2_{\vartheta,s,L}$---we have
$\diam\gamma\ge c(\epsilon c^{-2})\sqrt{|V(\Gamma)|}$ whenever~$\varkappa\ge 
\epsilon/c$. Moreover, the upper bound on~$\mathscr{W}_\beta(\frakS)$ from \eqref{G1} 
along with the estimate~$\mathscr{W}_\beta(\frakS)\ge\tau_\mini\varkappa\sqrt{v_L}$ imply
that~$\lambda_\Gamma$ is bounded away from zero and thus 
$\epsilon\sqrt{|V(\Gamma)|}=\epsilon\sqrt{\lambda_\Gamma v_L}\ge s$ for~$L$ 
sufficiently large. This verifies the assumptions of Lemma~\ref{lemma-onecontour} 
with~$\epsilon$ replaced by~$\epsilon c^{-2}$, which then guarantees that~$\Gamma$ 
is a singleton,~$\Gamma=\{\gamma_0\}$, and~that
\begin{equation}
%\label{}
\inf_{z\in\R^2}\dH\bigl(V(\gamma_0),\sqrt{|V(\gamma_0)|}W+z\bigr) \le
\sqrt{\epsilon}\sqrt{|V(\gamma_0)|}.
\end{equation}
Now, $|V(\gamma_0)|=\lambda_\Gamma v_L\le v_L$ (because, as noted before,~$\lambda_\Gamma\le1$), which 
means that the right-hand side is less than~$\sqrt{\epsilon v_L}$ and 
\eqref{hausbd} holds. But on~$\EE_{\epsilon,\varkappa,\vartheta}^i$ the event 
$\BB_{\epsilon,s,L}$ must fail, so we must have 
either that $\Phi_\Delta(\lambda_\Gamma)>\Phi_\Delta^\star+\epsilon$, 
which only applies when $i=1$, or that
\eqref{magbd} fails, which only applies when~$i=2$.

We claim that, in both cases, there exists an $\epsilon'>0$ 
and an~$\alpha>0$---both proportional to~$\epsilon$---such that for some $\frakS\sim\Gamma$
and~$L$ sufficiently large, we have
\begin{equation}
\label{exponentbd}
(1-\alpha)\mathscr{W}_\beta(\frakS)+
2\frac{(\mstar\,v_L)^2}
{\chi|\Lambda_L|}(1-\lambda_\Gamma+\epsilon_i-\kappa_L)^2\ge 
w_1\sqrt{v_L}\bigl(\Phi_\Delta^\star+\epsilon'\bigr).
\end{equation}
Indeed, the Wulff variational problem in conjunction with Lemma~\ref{L:skeleti}, 
the restriction to $(\RR^1_{\vartheta,s,L})^{\text{\rm c}}$ and the bound $(1-x)^{1/2}\ge 1-x$ for $x\in[0,1]$
imply that
\begin{equation}
\label{WS36}
\begin{aligned}
\mathscr{W}_\beta(\frakS)&\ge w_1|\V(\frakS)|^{1/2}\ge
w_1\Bigl(|V(\gamma_0)|-g_3\vartheta^{-1}s^2\sqrt{v_L}\Bigr)^{1/2}\\
&\ge w_1\sqrt{\lambda_\Gamma 
v_L}-g_3w_1\bigl(\vartheta\sqrt{\lambda_\Gamma}\bigr)^{-1}s^2.
\end{aligned}
\end{equation}
Observing also that the difference $2(\mstar)^2v_L^{3/2}/(\chi|\Lambda_L|)-w_1\Delta\to 0$ 
as $L\to \infty$,  the left
hand side of \eqref{exponentbd} can be bounded from below by
\begin{equation}
\label{tteerrmm}
w_1 \sqrt{v_L}\Phi_\Delta(\lambda_\Gamma)-\alpha w_1\sqrt{\lambda_\Gamma v_L}
-\delta_L\sqrt{v_L}
+2w_1\Delta\sqrt{v_L} (\epsilon_i-\kappa_L)\vartheta,
\end{equation}
where $\delta_L\to 0$ (as well as $\kappa_L\to 0$) with $L\to\infty$.
(Here we again used that $1-\lambda_\Gamma\ge\vartheta$.)
Now, for $i=1$ we have $\Phi_\Delta(\lambda_\Gamma)>\Phi_\Delta^\star+\epsilon$ 
from which \eqref{exponentbd} follows once $\alpha<\epsilon$ and $L$ is sufficiently large. 
For $i=2$, we use $\Phi_\Delta(\lambda_\Gamma)\ge \Phi_\Delta^\star$
and get the same conclusion since  \eqref{tteerrmm} 
now contains the positive term $2 w_1\Delta \epsilon_2\sqrt{v_L}
\propto \epsilon\sqrt{v_L}$.

By putting \eqref{vodhad} and \eqref{druhej} together, applying \eqref{exponentbd}, 
choosing~$K\ge K_0(\alpha,\beta)$ and invoking Lemma~\ref{lemma-Peierls} 
to bound the sum over all skeletons~$\frakS$, we find that
\begin{equation}
\label{1eq}
P_L^{+,\beta}\bigl(\MM_L\cap\EE_{\epsilon,\varkappa,\vartheta}
\cap\{\Gamma_s^\ext(\sigma)\in\boldsymbol\Gamma_1\}\bigr)\le 
2CC'\exp\bigl\{-
w_1\sqrt{v_L}\bigl(\Phi_\Delta^\star+\epsilon'\bigr)\bigr\}.
\end{equation}
whenever~$L$ is sufficiently large. 
(Here the embarrassing factor ``$2$'' comes from combining the corresponding estimates for $i=1$ and $i=2$.)

Thus, we are down to the cases~$\Gamma\in\boldsymbol\Gamma_2$, 
which means that for every skeleton $\frakS\sim\Gamma$, 
we have $\mathscr{W}_\beta(\frakS)> 
w_1\sqrt{\lambda_\Gamma v_L}(1+\epsilon c^{-2})$.
Moreover, since $\EE_{\epsilon,\vartheta,\varkappa}\subset\AA_{\varkappa,s,L}$, all $s$-large contours 
that we have to consider
actually satisfy that $\diam\gamma\ge\varkappa\sqrt{v_L}$. 
In particular, we
also have that $\mathscr{W}_\beta(\frakS)\ge \tau_\mini\varkappa\sqrt{v_L}$. 
Combining these bounds we derive that,
for some $c'>0$ and regardless of the value of $\lambda_\Gamma$,
\begin{equation}
%\label{}
\mathscr{W}_\beta(\frakS)\ge w_1\bigl(\sqrt{\lambda_\Gamma}+c'\bigr)\sqrt{v_L}.
\end{equation} 
Disregarding the factor $\epsilon_i$ in \eqref{vodhad} and performing
similar estimates as in the derivation of \eqref{1eq}, we find that \eqref{exponentbd} 
holds again for some~$\alpha>0$. Hence an analogue of \eqref{1eq} is valid also for all
$\Gamma\in\boldsymbol\Gamma_2$. A combination of these 
estimates in conjunction with Theorem~\ref{lowerbound} show that, indeed, 
\eqref{expbd} is true with a~$c_7$ proportional to~$\epsilon$. This finishes the 
proof.
\end{proofsect}

The previous proof immediately provides us with the proof of the other main 
results:

\begin{proofsect}{Proof of Theorem~\ref{LDP-thm}}
In light of Theorem~\ref{lowerbound}, we need to prove an appropriate upper bound on 
$P_L^{+,\beta}(\MM_L)$, where~$\MM_L=\{\sigma\colon M_L=\mstar|\Lambda_L|-
2\mstar\,v_L\}$. First we note that for~$L$ sufficiently large,
the probability 
$P_L^{+,\beta}(\MM_L)$ is comparable with~$P_L^{+,\beta}(\FF_L)$, 
where~$\FF_L$ is the event
\begin{equation}
%\label{}
\FF_L=\MM_L\cap\AA_{\varkappa,s,L}\cap\BB_{\epsilon,s,L}
\cap\bigl(\RR^1_{\vartheta,s,L}\cup\RR^3_{\vartheta,s,L}
\cup\RR^4_{\vartheta,s,L}\bigr)^{\text{\rm c}}
\end{equation}
with~$\epsilon$,~$\varkappa$,~$\vartheta$ as in the proof of 
Theorem~\ref{main-result}. 
But on~$\FF_L$, we have at most one large contour and the 
skeleton and Gaussian upper bounds readily give us that
\begin{equation}
%\label{}
P_L^{+,\beta}(\FF_L)\le C e^{-w_1\sqrt{v_L}(\Phi_\Delta^\star-\epsilon')}.
\end{equation}
for some $C<\infty$ and some $\epsilon'>0$ proportional to~$\epsilon$. 
From here and Theorem~\ref{lowerbound}, the claim \eqref{LDP} follows 
by letting $L\to\infty$ and $\epsilon\downarrow0$.
\end{proofsect}

Our last task is to prove Corollary~\ref{cor-main-result}.

\begin{proofsect}{Proof of Corollary~\ref{cor-main-result}}
By Proposition~\ref{prop-Phi}, if~$\Delta<\Deltac$, the unique minimizer of 
$\Phi_\Delta(\lambda)$ is~$\lambda=0$. Thus, for~$\epsilon>0$ sufficiently small 
and~$L$ large enough, the contour volumes are restricted to a small number times~$v_L$. 
Since \eqref{hausbd} says that the contour volume is proportional to the square of its diameter, this
(eventually) forces~$\diam\gamma<\varkappa\sqrt{v_L}$ for any fixed $\varkappa>0$. 
But that contradicts
the fact that~$\AA_{\varkappa,s,L}$ holds for a~$\varkappa$ sufficiently small.  Hence, no such
intermediate~$\gamma$ exists and all contours have a diameter smaller than~$K\log L$.

In the cases~$\Delta>\Deltac$, the function 
$\Phi_\Delta(\lambda)$ is minimized only by a non-zero~$\lambda$ (which is, in fact, 
larger than~$2/3$)
and so the scenarios without large contours are exponentially suppressed. 
Since, again,~$\diam\gamma>\varkappa\sqrt{v_L}$ for all 
potential contours, Theorem~\ref{main-result} guarantees that there is only one such 
contour and it obeys the bounds \eqref{hausbd} and \eqref{volbd}. 
All the other contours have diameter less than $K\log L$.
\end{proofsect}

\section*{Acknowledgments}
\noindent
The research of R.K. was supported by
the grants GA\v{C}R~201/00/1149 and MSM~110000001.
The research of L.C.~was supported by the NSF under the grant DMS-9971016
and by the NSA under the grant NSA-MDA~904-00-1-0050.
R.K. would also like to thank 
the UCLA Department of Mathematics  and
the Max-Planck Institute for Mathematics in Leipzig for their hospitality
as well as the A.~von Humboldt Foundation whose Award made the stay in Leipzig possible.

\end{document}